\documentclass[journal,twoside,web]{ieeecolor}
\usepackage{epsfig,color,amsmath,cite}
\usepackage{amsfonts}
\usepackage[svgnames]{xcolor}
\usepackage{pgfplots,caption}
\usepackage{subfigure}

\usepackage{xcolor}
\usepackage{epstopdf}
\usepackage{amssymb}
\usepackage{url}
\usepackage{mathtools}
\usepackage{multirow}
\usepackage{ifthen}
\usepackage{stackrel}
\usepackage{booktabs}
\usepackage[flushleft]{threeparttable}
\usepackage{makecell}
\usepackage[linesnumbered,lined,boxed,commentsnumbered,ruled,longend]{algorithm2e}


\DeclareMathOperator{\diag}{diag}

\DeclareMathOperator{\rank}{rank}

\makeatother

\usepackage{amsthm}
\newtheorem{remark}{Remark}
\newtheorem{thm}{Theorem}

\newtheorem{definition}{Definition}
\newtheorem{corollary}{Corollary}
\newtheorem{assumption}{Assumption}
\newtheorem{lemma}{Lemma}

\newcommand{\prob}{\mathbb{P}}

\newcommand{\m}{\boldsymbol}
\allowdisplaybreaks[4]
\newcommand{\mbb}[1]{\mathbb{#1}}
\newcommand{\mr}[1]{\mathrm{#1}}
\usepackage{graphicx} %Loading the package 
\usepackage{tikz}
\usetikzlibrary{matrix,positioning,decorations.pathreplacing}
\usepackage{array}
\usepackage{filecontents}
\usepackage[framemethod=TikZ]{mdframed}
\usepackage{bm}

\usepackage{multicol}

\newcommand{\E}{\mathbb{E}}

\newcommand{\F}{\mathbb{F}}

\newcommand{\Var}{\mr{Var}}
\newcommand{\Cov}{\mr{Cov}}

\usepackage{graphicx}
\usepackage{generic}
\usepackage[export]{adjustbox}
\usepackage{textcomp}

\def\BibTeX{{\rm B\kern-.05em{\sc i\kern-.025em b}\kern-.08em
    T\kern-.1667em\lower.7ex\hbox{E}\kern-.125emX}}
\usepackage{hyperref}
\hypersetup{
	colorlinks=true,
	linkcolor=black,
	filecolor=black,      
	urlcolor=cyan,
	citecolor=black
}
\begin{document}
\title{Probabilistic State Estimation in Water Networks}

\author{Shen Wan$\text{g}^\dagger$, Ahmad F. Tah$\text{a}^{\dagger, *}$, Nikolaos Gatsi$\text{s}^{\dagger}$,  Lina Sel$\text{a}^{\ddagger}$, and Marcio H. Giacomon$\text{i}^{**}$	% \vspace{-0.85cm}
	\thanks{
$^{*}$Corresponding author. $^\dagger$Department of Electrical and Computer Engineering, The University of Texas at San Antonio. $^{**}$Department of Civil and Environmental Engineering, The University of Texas at San Antonio. $^{\ddagger}$Department of Civil, Architectural and Environmental Engineering, Cockrell School of Engineering, The University of Texas at Austin. Emails: mvy292@my.utsa.edu, \{ahmad.taha, nikolaos.gatsis, marcio.giacomoni\}@utsa.edu, linasela@utexas.edu. This material is based upon work supported by the National Science Foundation under Grants 1728629 and 2015671. }}

%\author{First A. Author, \IEEEmembership{Fellow, IEEE}, Second B. Author, and Third C. Author, Jr., \IEEEmembership{Member, IEEE}
%\thanks{This paragraph of the first footnote will contain the date on 
%which you submitted your paper for review. It will also contain support 
%information, including sponsor and financial support acknowledgment. For 
%example, ``This work was supported in part by the U.S. Department of 
%Commerce under Grant BS123456.'' }
%\thanks{The next few paragraphs should contain 
%the authors' current affiliations, including current address and e-mail. For 
%example, F. A. Author is with the National Institute of Standards and 
%Technology, Boulder, CO 80305 USA (e-mail: author@boulder.nist.gov). }
%\thanks{S. B. Author, Jr., was with Rice University, Houston, TX 77005 USA. He is 
%now with the Department of Physics, Colorado State University, Fort Collins, 
%CO 80523 USA (e-mail: author@lamar.colostate.edu).}
%\thanks{T. C. Author is with 
%the Electrical Engineering Department, University of Colorado, Boulder, CO 
%80309 USA, on leave from the National Research Institute for Metals, 
%Tsukuba, Japan (e-mail: author@nrim.go.jp).}}

\maketitle

\begin{abstract}
State estimation in water distribution networks (WDN), the problem of estimating all unknown network heads and flows given select measurements, is challenging due to the nonconvexity of hydraulic models and significant uncertainty from water demands, network parameters, and measurements. To this end, a probabilistic modeling for state estimation (PSE) in WDNs is proposed. After linearizing the nonlinear hydraulic WDN model, the proposed PSE shows that the covariance matrix of unknown system states (unmeasured heads and flows) can be linearly expressed by the covariance matrix of three uncertainty sources (i.e., measurement noise, network parameters, and water demands). Instead of providing deterministic results  for unknown states, the proposed PSE approach  \textit{(i)} regards the system states and uncertainty sources as random variables and yields variances of individual unknown states,  \textit{(ii)} considers thorough modeling of various types of valves and measurement scenarios in WDNs, and \textit{(iii)} is also useful for uncertainty quantification, extended period simulations, and confidence limit analysis.  The effectiveness and scalability of the proposed approach is tested using several WDN case studies.
\end{abstract}

\begin{IEEEkeywords}
Water distribution networks, probabilistic state estimation, confidence limit analysis, uncertainty quantification.
\end{IEEEkeywords}
\section{Introduction and Literature Review}~\label{sec:literature}
\IEEEPARstart{S}{}tate estimation (SE) is a technique used to estimate the unknown state variables based on a set of measurements, a mathematical model linking the measurable and unmeasurable (unknown) variables, and model parameters~\cite{tshehla2017state}. SE plays an important role in a water distribution network (WDN) with a variety of applications such as detecting water loss~\cite{diaz2018probabilistic}. System states in a WDN include the collections of heads and flow rates in our paper. These variables depend on water demands, which are characterized by frequent temporal changes and inherent variability and uncertainty; hence the state variables require frequent estimation. The parameter estimation or calibration problem~\cite{Kumar2010,Diaz2017}  is closely related to SE problem, in which unknown or unmeasurable model parameters need to be estimated based on measurements and the mathematical model.

However, the  distinction between parameter and state estimation is manifested in the frequency of calibration or the time scales---and their subsequent applications. For example, the studies~\cite{Kumar2010,Diaz2017} investigate estimation of fast-changing variables (for states such as heads, flow rates, settings, and demands) and slow-changing parameters (pipe roughness coefficients). In our paper, the model parameters, i.e., pipe roughness coefficients, statuses or settings of pumps and valves, are assumed to be known~\cite{diaz2018topological}.

The measurements in a WDN typically include \textit{(i)} a small (relative to the system size) subset of flows in pipes and pumps, heads in junctions and storage tanks which can be measured by a supervisory control and data acquisition (SCADA) system, and \textit{(ii)} consumer demands which are viewed as \textit{pseudo-measurements}. 

The WDN mathematical model linking system states, measurements, and parameters is built on the principles of \textit{conservation of mass} and \textit{energy}. The former implies the continuity of flow at each node, and the latter states that the energy difference between two connected nodes is equal to the energy losses or gains due to pipe friction or pumping~\cite{puig2017real}. With such a mathematical model, the unknown and unmeasurable system states can be computed or estimated with known measurements from sensors located in key network locations and pseudo-measurements.
 
In particular, the SE problem is referred to as deterministic state estimation (DSE) or {point state estimation} when all measurements, parameters, and variables used in that mathematical model are considered to be deterministic.   However, measurements and parameters in practice contain a significant amount of uncertainty~\cite{kang2009approximate} which might lead to inaccurate SE---and mislead the planning and management of the WDN.  To that end, the paper's \textit{\textbf{objective}} is three-fold: \textit{(i)} to propose probabilistic modeling for SE (PSE) via viewing all system states, network parameters, and uncertainty as random variables,  \textit{(ii)} to investigate scalable computational algorithms to solve the probabilistic SE problem in WDN considering the various sources of uncertainty, and \textit{(iii)} to explore the impact on system states from the uncertainty sources under general or specific statistical distributions.

In this work, the typical sources of uncertainty are considered including sensor measurement noise, demand estimation errors, and WDN pipe parameters.  Other sources of uncertainty such as cyber-attacks to SCADA systems~\cite{tcst1,taormina2018battle} or pressure-deficient conditions that could incur from pump failures, pipe leaks, and fire demands, are not considered in our study. The literature of SE in WDN is rich and is summarized as follows. 
\subsection{Literature review}
The DSE problem is typically formulated as an inverse optimization problem in which the objective is minimizing the error between the mathematical model and measurements, constrained by the network hydraulics and parameter bounds/limits~\cite{tshehla2017state}. DSE formulations treat each variable as a deterministic one, and uncertainty is often not considered yielding a deterministic solution. The studies in~\cite{andersen2001constrained,tshehla2017state,diaz2018topological,wang2019state,preis2011} make up the bulk of recent DSE literature. 

DSE formulations can result in unconstrained or  constrained problems depending on the specific optimization approach and problem formulation. The unconstrained DSE, can be divided into the following groups~\cite{arsene2014mixed}: \textit{(i)} least squares (LS) where the sum of squared differences between measurement and model is minimized, \textit{(ii)} least absolute value (LAV) where the sum of absolute differences is minimized, and \textit{(iii)} minimax where the maximum difference is minimized. The sensitivity to errors of above methods varies.
%, and the corresponding weighted versions are proposed, e.g., weighted least square (WLS), weighted least absolute value (WLAV). The authors in~\cite{andersen2001constrained} propose a constrained DSE to overcome the downside of a weighted based approach via adding parts of objective functions to the constraints. 
Recently, we propose a novel geometric programming (GP)-based method~\cite{wang2019state} to solve the  DSE problem which relies on derivative-free and scalable formulations.

The aforementioned DSE studies do not comprehensively consider the inherent uncertainty embedded in WDN. As previously mentioned, the main sources of uncertainty stem from: \textit{(i)} measurement noise of sensors, \textit{(ii)} demands, and \textit{(iii)}  pipe roughness coefficients~\cite{bargiela1989pressure,arsene2011confidence,du2018direct,Vrachimis2018RealtimeHI,kang2010demand,diaz2016uncertainty}.
In order to overcome the limitations of DSE, uncertainty quantification in WDN---also termed as confidence limit analysis (CLA)---is first proposed in~\cite{bargiela1989pressure}. Specifically, the authors present three CLA techniques which can calculate the inaccuracy of heads and flows caused by the uncertainty of measurements and network parameters. The techniques are based on Monte Carlo Simulation (MCS) method, an optimization method, and a sensitivity matrix technique. The study~\cite{arsene2011confidence} presents two CLA-based techniques based on a least squares loop flows state estimation. This approach computes the confidence limits (lower and upper limit) for the system state variables (i.e., flows and heads at unmeasured locations). Another study~\cite{du2018direct} presents a method quantifying uncertainties propagated from measurement noises, demand uncertainties, and model simplification errors to the states in their uncertainty quantification section.
 
In~\cite{kang2010demand}, a two-step sequential method for estimation of demand and pipe roughness coefficient is presented based on a weighted LS scheme. The uncertainty in estimated variables and resulting nodal head predictions are quantified in terms of confidence limits using first-order second moment (FOSM) method. The authors in~\cite{diaz2016uncertainty} set out an alternative approach with respect to the weighted LS problem to determine the upper and lower limits of states. This enables general quantification of SE uncertainty for all state variables by applying the FOSM method. Moreover, it enables the computation of the covariance matrix of state variables. Another branch of CLA is interval state estimation (ISE), introduced recently in~\cite{vrachimis2018iterative,Vrachimis2018RealtimeHI}, based on interval arithmetic which aims to find the DSE solution region, rather than finding a specific point-based solution.
\subsection{Research gaps, paper contributions and organization}
The research gaps are summarized as follows. First, the deterministic SE studies~\cite{andersen2001constrained,tshehla2017state,preis2011,diaz2018topological,wang2019state} do not consider critical uncertainty (demand and pipe coefficient uncertainty) as random variables in the network with different types of valves and various measurement scenarios that could lead to significantly different estimates for unmeasured state variables.  Second, the methods based on  CLA or uncertainty quantification~\cite{bargiela1989pressure,du2018direct,arsene2011confidence,diaz2016uncertainty,kang2010demand,Vrachimis2018RealtimeHI} suffer from the following limitations. The first one is that there is an absence of deriving or studying {all three uncertainty sources simultaneously and their impact on state estimates} as well as the covariance of uncertain variables and unknown system states variables. We note that FOSM only produces a mapping from some uncertainty sources to specific systems states (particular heads and flows). 
Furthermore, the CLA approach focuses on how the state variables are impacted by uncertainty instead of  studying the relationship between probability distribution functions (PDF) of state variables and PDF of uncertainties when variables follow a certain type of distribution, i.e., normal distribution and uniform distribution. 

Third, we note that the probabilistic modeling proposed in~\cite{xu1998probabilistic}, built for water distribution reliability analysis instead of state estimation,  only considers the SE of heads. Furthermore, it is unclear how the method in~\cite{xu1998probabilistic} can solve the over-determined measurement scenario, i.e., when there are more observations or measurements than unknown variables~\cite{powell1999state}. A novel contribution in~\cite{xu1998probabilistic} is the authors' correct conjecture that 
\begin{quote}
\textit{If the uncertainty from nodal demands, reservoir levels, and pipe roughness coefficients are normally distributed, the linearized nodal heads also follow a normal distribution.}
\end{quote} However, no formal proof is given, and the PDF produced in~\cite{xu1998probabilistic} only depicts heads in the network. Specifically, it remains unclear how changes in the uncertainty distribution (i.e., if the uncertainty follows a non-Gaussian distribution) impacts PDF of heads and flows in the network. Finally, it is also unclear how the aforementioned studies can be extended to various types of valves. 

The paper's objective is to address the aforementioned research gaps by  proposing probabilistic modeling for SE that generalizes the approaches in the literature while considering various sources of uncertainty with arbitrary statistical distributions, types of valves, and sufficient/over-determined measurement scenarios.  The specific paper contributions are:
\begin{itemize}
\item The presented PSE  views system states, measurements, and network parameters as random variables, and offers a general method to connect unknown system states with uncertain variables. This contribution generalizes various types of distributions (i.e., uniform, normal, and Laplace distributions) of the uncertainty random variables as long as their expectation and covariance exist. 

\item A scalable algorithm to find the covariance matrix of system states is proposed given the covariance matrix of uncertainty through solving a linear system of equations. In particular, we prove that the linear systems of equations have a solution under mild conditions. Various types of valves and measurement scenarios are considered in this algorithm which makes it applicable to various WDN state estimation scenarios.

\item We also show that if  the uncertainty follows the normal or uniform distribution (motivated by empirical evidence), then system state variables also follow the corresponding distribution.  Uncertainty propagation in water networks is studied and the importance or impact of each uncertainty source is given and verified by case studies.

\item  {Thorough case studies are presented showcasing the scalability and effectiveness of the proposed PSE formulations in dealing with uncertainty.}
\end{itemize}

The rest of the paper is organized as follows. Section~\ref{sec:ModelAssumption} introduces WDN modeling and assumptions.  {PSE} formulation is given first, then methods and the corresponding algorithm to solve the {PSE} are presented in Section~\ref{sec:formulation}.  The generalization of flows and heads distributions are discussed in Section~\ref{sec:distribution}. Section~\ref{sec:Discussion} presents important discussions and insights related to the proposed PSE algorithm and different measurement scenarios. Section~\ref{sec:Test} presents case studies to corroborate the paper's theoretical findings. Mathematical proofs and extensive details about the presented models are all given in the appendices.  The notation for this paper is introduced next.

\noindent {\textbf{Notation.}} 
\noindent Italicized, boldface upper and lower case characters represent matrices and column vectors: $a$ is a scalar, $\m a$ is a vector, and $\m A$ is a matrix. The notation $\mathbb{R}^n$   denotes the sets of column vectors with $n$  real  numbers.
For $\m x \in\mbb{R}^m$, $\m y \in \mbb{R}^n$, a compact column vector in $\mathbb{R}^{m+n}$  is defined as $\{\m x\, , \m y\} = [\m x^\top \ \m y^\top]^\top$.  Similarly, for matrices $\m A$ and $\m B$ with same number of columns, $\{\m A\, , \m B\}$ stands for $[\m A^\top \ \m B^\top]^\top$. For any random vector $\m x, \m b \in \mathbb{R}^{n}$, and ${ \m x =[x_{1},x_{2},...,x_{n}]^{\top}}$, the $\E{(\m x)} = [\E(x_{1}),\E(x_{2}),...,\E(x_{n})]^{\top}$, $\Var{(\m x)}= [\Var(x_{1}),\Var(x_{2}),...,\Var(x_{n})]^{\top}$, and ${ \operatorname {K}_{\m x \m x }}$ are the expectation, variance, and  covariance matrix of the  vector $\m x$. We also define the operation of covariance over a linear system of equations $\m A \m x = \m b$, the notation  $\Cov(\m A  \m x = \m b)$ stands for $\Cov(\m A \m x, \m A \m x) = \Cov(\m b, \m b)$, which is to apply covariance operator to each side of $\m A  \m x = \m b$. By this notation, another equivalent formulation of $\m A  \m x = \m b$ is $\m A \m x -  \m b = \m 0$, and the operation of covariance over it results in $\Cov((\m A \m x - \m b),(\m A \m x - \m b))=\m 0$. We note that $\Cov(\m A \m x, \m A \m x) = \Cov(\m b, \m b)$  and $\Cov((\m A \m x - \m b),(\m A \m x - \m b))=\m 0$  are equivalent; see Appendix~\ref{app:covform}.

The variables with upper case characters $\m \cdot^{\mathrm{J}}$, $\m \cdot^{\mathrm{R}}$, $\m \cdot^{\mathrm{TK}}$,  $\m \cdot^{\mathrm{P}}$, $\m \cdot^{\mathrm{M}}$, and  $\m \cdot^{\mathrm{L}}$ represent the variables related to junctions, reservoirs, tanks, pipes, pumps, and valves.  
\section{WDN Modeling and Assumptions}~\label{sec:ModelAssumption}
We define the column vectors $\m h^{\mathrm{J}} \in \mathbb{R}^{n_j}$  , $\m h^{\mathrm{R}} \in \mathbb{R}^{n_r}$, and $\m h^{\mathrm{TK}} \in \mathbb{R}^{n_t}$ to collect the heads at $n_j$ junctions, $n_r$  reservoirs, and $n_t$  tanks, respectively; the column vectors $\m q^{\mathrm{P}} \in \mathbb{R}^{n_p} $, $\m q^{\mathrm{M}} \in \mathbb{R}^{n_m}$ and $\m q^{\mathrm{L}} \in \mathbb{R}^{n_l}$ to collect the flow through $n_p$ pipes, $n_m$ pumps, and $n_l$ valves, respectively. Then, the compact column vectors $\m h \in \mathbb{R}^{n_{h}}$ and $\m q \in \mathbb{R}^{n_{q}}$ collect all heads and flows are defined as $\m h  \triangleq   \left\lbrace \m h^{\mathrm{J}}, \m h^{\mathrm{R}}, \m h^{\mathrm{TK}}\right\rbrace$, $\m q  \triangleq  \left\lbrace \m q^{\mathrm{P}}, \m q^{\mathrm{M}},\m q^{\mathrm{L}}\right\rbrace .$
Vector $\m x \in \mathbb{R}^{n_{x}}$ at time $k$, collects all heads and flows, is defined as
\begin{align} \label{equ:compactxi}
{\m x(k) \triangleq  \left\lbrace \m h(k), \m q(k)  \right\rbrace }.
\end{align}
Note that $n_{h}  = n_j + n_r + n_t$, $n_{q} = n_p + n_m + n_l$, and $n_{x} = n_{h} + n_{q}$.

\subsubsection{DAE form of WDN} 
The modeling of a WDN can be written in the form of difference algebraic equation (DAE)
\begin{subequations} ~\label{equ:dae-abstract}
	\begin{align}
	\mathrm{DAE}:\;\;\;	\m h^{\mathrm{TK}}(k + 1) &= \m A_{\m h} \m h^{\mathrm{TK}}(k) + \m B_{\m q} \m q(k)  ~\label{equ:tankhead-abcstracted} \\
	\hspace{-10pt}	\m 0&=\m E_{\m q}^\mathrm{J} \m q(k) + \m d(k) ~\label{equ:nodes-abcstracted}\\
	\hspace{-10pt}\m 0 &= \m E_{\m h} \m h(k) + \m \Phi(\m q(k), \m c(k))~\label{equ:PumpPipe-abstract},
	\end{align}
\end{subequations}
where $n_t \times n_t$ matrix $\m A_{\m h}$, $n_t \times n_{q}$ matrix $\m B_{\m q}$, $n_j \times n_{q}$ matrix  $\m E_{\m q}^\mathrm{J}$, and $n_{q} \times n_{h}$ matrix  $\m E_{h}$ are constant matrices that depend on the topology and hydraulic properties of the underlying WDN.  Equation~\eqref{equ:tankhead-abcstracted} collects the dynamic equations of tanks~\eqref{equ:tankhead}; Equation~\eqref{equ:nodes-abcstracted} collects the mass balance equations for all junctions~\eqref{equ:nodes}; the nonlinear function $\m \Phi(\cdot)$ in Equation~\eqref{equ:PumpPipe-abstract}  includes  the nonlinear pipe~\eqref{equ:head-flow-pipe}, pump~\eqref{equ:head-flow-pump}, and valve~\eqref{equ:head-fcv-valve},~\eqref{equ:head-prv-valve} models in Appendix~\ref{sec:Modeling}. Note that the statuses and settings of pressure reducing and flow control valves are assumed to be known in~\eqref{equ:head-fcv-valve},~\eqref{equ:head-prv-valve}.
The roughness coefficients for pipes are collected in vector $\m c(k)$. 
\subsubsection{Operational limits}  
The state of operational conditions in  WDN is bounded by physical and operational constraints, hence, the overall operational limits of hydraulic elements can be expressed as
\begin{subequations}\label{equ:constraints}
	\begin{align}
	\hspace{-10em}\mathrm{Limits} \hspace{2em}
	\m x_{i}(k) &= \m x^{\mathrm{set}}_{i}(k) \label{equ:equality}\\
	\m x_{j}^{\mathrm{min}} \leq  &\,\,\m x_{j}(k) \leq \m x_{j}^{\mathrm{max}}, \label{equ:inequality}
	\end{align}
\end{subequations} 
where $i \neq j$ indicates that each element in $\m x$ is either limited by equality or inequality expression. For example, we assume that the head at  Reservoir $i$ is constant for simplicity, but they can be seamlessly modeled as uncertain similar to water levels in tanks. The relative speeds of pumps are assumed to be fixed as $1$ (full speed) and the flow in pipe $j$ is limited by $ q_{j}^{\mathrm{min}} \leq   q_{j}(k) \leq  q_{j}^{\mathrm{max}}$.

\subsubsection{Measurement modeling} 
The model for WDN measurements (select heads and flows at certain nodes in the network) can be expressed as
\begin{align}\label{equ:measure}   
\mathrm{Measurements}:\;\;\;	\m y(k) =  \m C \m  x(k) + \m v(k),
\end{align}
where $\m y(k) \in \mathbb{R}^{n_{y}}$ is measurement vector, $\m C$ is the binary selection matrix with $n_{y} \times n_{x}$ depicting where sensors are installed, and $\m v(k) \in \mathbb{R}^{n_{y}} $ is  measurement noise vector; see~\cite{tcst4} for a study on sensor placement in water networks.

Before we propose PSE formulation next, an assumption and an definition are introduced. Common statistical assumptions~\cite{xu1998probabilistic,kang2009approximate,andersen2001constrained} is given next. 
\begin{assumption}\label{assumtion:all}
	%	Demand $\m d$, measurement noise $\m v$, and  pipe roughness coefficient $\m c$ are assumed to be independent of each other. 
	The entries of the demand vector $\m d$, measurement noise vector $\m v$, and pipe roughness coefficient vector $\m c$ are mutually independent.
\end{assumption}
The assumption is reasonable seeing there is no connection between the three uncertain variables. The definition of the cross-covariance matrix is given as follows.
\begin{definition}~\label{def:cc}
	For  column vectors
	${ \m x =[x_{1},x_{2},...,x_{m}]^{\top}}$ and ${ \m y =[y_{1},y_{2},...,y_{n}]^{\top}}$ consisting of random variables, then the cross-covariance matrix ${ \operatorname {K}_{\m x \m y }}$~\cite{Park2018} is the matrix whose ${ (i,j)}$ entry is the covariance between $x_{i}$ and $y_{j}$. That is
	\begin{equation}
	\hspace{-1em}\small{{  \operatorname {K}_{\m x \m y }} = \begin{bmatrix}
		\Cov(x_1,y_1) & \Cov(x_1,y_2) & \ldots & \Cov(x_1,y_n) \\
		\Cov(x_2,y_1) & \Cov(x_2,y_2) & \ldots & \Cov(x_2,y_n) \\
		\vdots & \vdots & \ddots & \vdots \\
		\Cov(x_m,y_1)&\Cov(x_m,y_2)&\ldots& \Cov(x_m,y_n)\\
		\end{bmatrix}}.
	\end{equation}
\end{definition}
Note that \textit{(i)} $\Cov(\m x, \m y) =  \operatorname {K}_{\m x \m y }$, \textit{(ii)} when $\m x= \m y$, the cross-covariance matrix ${ \operatorname {K}_{\m x \m y }}$ turns into ${ \operatorname {K}_{\m x \m x }} $, which is the covariance matrix,  and \textit{(iii)}  according to the definition and notation of $\Var(\m x)$,  each entry in $\Var(\m x)$ is on the diagonal of ${ \operatorname {K}_{\m x \m x }}$. 
\section{Formulating and Solving the PSE Problem}\label{sec:formulation}
As mentioned in Section~\ref{sec:literature}, the uncertainty of SE lies in demand in the mass balance equations $\m d(k)$, pipe roughness coefficients $\m c(k)$, and measurement noise $\m v(k)$, and these uncertainty sources would impact the system state $\m x (k)$. The objective of the presented probabilistic state estimation (PSE) is to find the covariance of the unknown variable $\m x (k)$, or $\operatorname{K}_{\m x\m x}$, over a time-horizon $k$ to $k+T$ thereby producing the variance of all variables in vector $\m x$. This contrasts point-based estimation which focuses on generating estimates for $\m x(k)$.  To formulate the PSE problem, we first obtain an operating point $\m x^0 := \{\m x^0 (k),\ldots, \m x^0(k+T) \}$ through solving the existing and well-developed deterministic state estimation routine (see~\cite{wang2019state} or a general DSE formulation in~\cite{tshehla2017state}) given the means of all uncertainty for all time-steps and given the measurement vector $\m y^0$. For example, we consider that a demand prediction is given from $k$ to $k+T$. The PSE can be depicted as the following \textit{high-level optimization problem}
\begin{subequations}\label{equ:pse-formulation-variance}
	\begin{align}  
	& \mathrm{find}  \hspace{8em} \operatorname {K}_{\m x \m x }  \label{equ:pse-formulation-obj1}  \\
	&\mathrm{s.t.}\; \;\Cov(\mathrm{DAE}~\eqref{equ:dae-abstract}, \mathrm{Limits}~\eqref{equ:equality},\mathrm{Measurements}~\eqref{equ:measure}) \label{equ:pse-constraint-variance}.
	\end{align}
\end{subequations}
Problem~\eqref{equ:pse-formulation-variance} finds a feasible covariance matrix of $\m x$ while satisfying the constraints~\eqref{equ:pse-constraint-variance}. The constraint set is an implicit function of the covariance matrix. The set also physically defines the probabilistic propagation of uncertainty to the system states. Before proceeding to the paper's approach, we emphasize the following traits of~\eqref{equ:pse-formulation-variance}:\\
\noindent \textit{(i)} Problem~\eqref{equ:pse-formulation-variance} is nonconvex and extremely difficult to solve for large networks, due to the nonlinear, nonconvex hydraulic constraints and the corresponding covariance operator. 

\noindent \textit{(ii)} The covariance operation of equality limits~\eqref{equ:equality} has clear physical meaning: Reservoir $i$ has a fixed head $h_i$, the deviation of head is $\Cov(h_i,h_i) = 0$ implies that the head at Reservoir $i$ does not change and is deterministic. However, the covariance operation over the inequality constraint~\eqref{equ:inequality} is not meaningful and thus not included in~\eqref{equ:pse-constraint-variance}. With that in mind, the operating point $\m x^0$  satisfies~\eqref{equ:inequality}.

\noindent \textit{(iii)} The optimization variable is $\operatorname {K}_{\m x \m x}$ which is encoded in~\eqref{equ:pse-constraint-variance} after performing the covariance operation on $\mathrm{DAE}$, $\mathrm{Limits}$, and $\mathrm{Measurements}$ models. 

The objective of this paper is to solve a simplified version of~\eqref{equ:pse-formulation-variance} through a scalable computational method. 
\subsection{Linear modeling}
We regard the elements in $\m {x}$ as random variables from a statistical perspective, apply covariance operation $\Cov$ on the constraints in~\eqref{equ:pse-formulation-variance}, and use the law of covariance of linear combination to build the relationship between uncertainty sources and system state $\m {x}$. However, the DAE model~\eqref{equ:dae-abstract} remains nonlinear and the law of covariance of linear combination can not be applied. Hence, the first step to solve~\eqref{equ:pse-formulation-variance} is to linearize the nonlinear hydraulic model around an operating point $\m x^0$. Note that this linearization procedure does not compromise modeling integrity due to slow hydraulic time constants in WDN.

From the nonlinear DAE model~\eqref{equ:dae-abstract} for various network components and hydraulic models, we obtain a linearized DAE model around $\m x^0$.  
%see that flow control valve (FCV) and pressure reducing valve (PRV) are piece-wise linear since minor head loss of valves are not considered. Hence, only the nonlinearity from pipes~\eqref{equ:head-flow-pipe} and pumps~\eqref{equ:head-flow-pump} need to be linearized, the results are presented as~\eqref{equ:head-flow-pipe-linear} and~\eqref{equ:head-flow-pump-linear} in Tab.~\ref{tab:models} in 
Appendices~\ref{sec:Modeling} and~\ref{sec:LinearDAE} include the complete derivation for the nonlinear and linearized DAE models. The linearized, compact form for all pipes and pumps can be written as 
\begin{subequations} \label{equ:linearPipesPumps}  
	\begin{align}
	\m \Delta \m h^\mathrm{P} (k)  &= \m K_q^{\mathrm{P}}  \m q(k)   + \m K_c^{\mathrm{P}}  \m c(k)  + \m b^{\mathrm{P}} \label{equ:linearPipe}\\
	\m \Delta \m h^\mathrm{\mathrm{M}} (k)  &=  -\m K_q^{\mathrm{M}} \m q (k)  + \m b^{\mathrm{M}}.  \label{equ:linearPump}
	\end{align}
\end{subequations}
where $\m K_q^{\mathrm{P}}$, $\m K_c^{\mathrm{P}}$, and $\m K_q^{\mathrm{M}}$ are diagonal slope matrices with size $n_p \times n_{p}$, $n_p \times n_{p}$, and $n_m \times n_{m}$; vectors $\m b^{\mathrm{P}}$ and $\m b^{\mathrm{M}}$ are intercepts with size $n_p \times 1$ and $n_m \times 1$. These matrices/vectors are all known and can be calculated around operating point $\m x^0$ efficiently. All pipe roughness coefficients are collected in $ \m c(k)$. Hence, the linear DAE can be expressed as 
\begin{subequations} ~\label{equ:dae-linear}
	\begin{align}
	\mathrm{DAE_{linear}}:\;\;\;	\m h^{\mathrm{TK}}(k + 1) &= \m A_{\m h} \m h^{\mathrm{TK}}(k) + \m B_{\m q} \m q(k)  ~\label{equ:tankhead-linear} \\
 	\m 0&=\m E(k) \m x(k) - \m z(k).  ~\label{equ:linearHD}
	\end{align} 
\end{subequations}
where $\m z(k) \triangleq \{\m d(k), \m K_c^{\mathrm{P}} \m c(k) + \m b^{\mathrm{P}}(k) , \m  b^{\mathrm{M}}(k), \m l(k)\}$, and $\m l(k)$ is an $n_l \times 1$ vector. Equation~\eqref{equ:nodes-abcstracted}, collecting the mass balance equations for all junctions, is included in~\eqref{equ:linearHD}. Matrix $\m E(k)$ in~\eqref{equ:linearHD} depends on the topology and parameters in the linearized model. Details of the linearized model are all given in Appendix~\ref{sec:LinearDAE}. 
\subsection{Reformulating the PSE problem}
The PSE can now be written as
	\begin{align} 
	\mathrm{find}& \hspace{8em} \operatorname {K}_{\m x \m x }  \label{equ:pse-linear}   \\
	\mathrm{s.t.}& \hspace{2pt} \Cov(\mathrm{Limits}~\eqref{equ:equality},\,\mathrm{Measurements}~\eqref{equ:measure},\,\mathrm{DAE_{linear}}~\eqref{equ:dae-linear}). \notag
	\end{align}
Interestingly, all constraints in~\eqref{equ:pse-linear} can be expressed as a compact linear system  of equations $\m A[k] \m x [k] = \m b[k]$, where 
\begin{align*}
\m x[k] &\triangleq \{\m x(k), \m x(k+1),\ldots, \m x(k+T) \} \in \mbb{R}^{Tn_x}  \\
\m A[k] &\triangleq  \{\m A^\mathrm{s}[k], \m A^\mathrm{TK}[k] \}, \hspace{8pt}\m b[k] \triangleq  \{\m b^\mathrm{s}[k], \m b^\mathrm{TK}[k] \}.
\end{align*}
We now show the structure of the linear system of equations $\m A[k] \m x [k] = \m b[k]$. This linear system of equations can be obtained as
\begin{equation*}
\resizebox{\columnwidth}{!}
{ 
\begin{tikzpicture}[
baseline={([yshift=-.5ex]current bounding box.center)},
style0/.style={
	matrix of math nodes,
	every node/.append style={text width=#1,align=center,minimum height=3.5ex},
	nodes in empty cells,
		inner sep=0pt,
	left delimiter=.,
	right delimiter=.,
},
style1/.style={
	matrix of math nodes,
	every node/.append style={text width=#1,align=center,minimum height=3.5ex},
	nodes in empty cells,
	inner sep=0.1pt,
	left delimiter=[,
	right delimiter=],
},
style2/.style={
	matrix of math nodes,
	every node/.append style={text width=#1,align=center,minimum height=3.5ex},
	nodes in empty cells,
		inner sep=0.1pt,
	left delimiter=\lbrace,
	right delimiter=\rbrace,
}
]

\matrix[style1=1.6cm] (1mat)
{
	& & & & \\
	& & & &  \\
	& & & &  \\
	& & & &  \\
	& & & &  \\
	& & & &  \\
	& & & &  \\
	& & & &  \\
	& & & &  \\
};
\draw[dashed]
(1mat-5-1.south west) -- (1mat-5-5.south east);
\draw[dashed]
(1mat-5-1.south west) -- (1mat-5-2.south east);
\draw[dashed]
(1mat-6-1.south west) -- (1mat-6-3.south east);
\draw[dashed]
(1mat-7-2.south west) -- (1mat-7-3.south east);
\draw[dashed]
(1mat-8-4.south west) -- (1mat-8-5.south east);
\draw[dashed]
(1mat-9-4.south west) -- (1mat-9-5.south east);
\draw[dashed]
(1mat-6-1.north west) -- (1mat-6-1.south west);
\draw[dashed]
(1mat-6-2.north east) -- (1mat-6-2.south east);
\draw[dashed]
(1mat-7-1.north east) -- (1mat-7-1.south east);
\draw[dashed]
(1mat-7-3.north east) -- (1mat-7-3.south east);
\draw[dashed]
(1mat-9-3.north east) -- (1mat-9-3.south east);
\draw[dashed]
(1mat-9-5.north east) -- (1mat-9-5.south east);
\node[font=\large]
at (1mat-1-1) {$\m A^\mathrm{s}(k)$};
\node[font=\large]
at (1mat-2-2) {$\m A^\mathrm{s}(k\hspace{-2pt}+\hspace{-2pt}1)$};
\node[font=\large]
at (1mat-3-3) {$\m A^\mathrm{s}(k\hspace{-2pt}+\hspace{-2pt}2)$};
\node[font=\large]
at (1mat-4-4) {$\ddots$};
\node[font=\large]
at ([xshift=-12pt]1mat-5-5) {$\m A^\mathrm{s}(k\hspace{-2pt}+\hspace{-2pt}T)$};
\node[font=\large]
at ([xshift=-20pt]1mat-6-2) {$\m A^\mathrm{TK}(k)$};
\node[font=\large]
at ([xshift=-20pt]1mat-7-3) {$\m A^\mathrm{TK}(k\hspace{-2pt}+\hspace{-2pt}1)$};
\node[font=\large]
at ([xshift=-20pt]1mat-8-4) {$\ddots$};
\node[font=\large]
at ([xshift=-22pt]1mat-9-5) {$\m A^\mathrm{TK}(k\hspace{-2.5pt}+\hspace{-2.5pt}T\hspace{-3.5pt}-\hspace{-3.5pt}1)$};
\draw [decorate, decoration={brace,mirror,amplitude=10pt,raise=0.35cm}] (1mat-9-1.west) -- (1mat-9-5.east)  node[midway,yshift=-0.9cm] {\large $\m A[k]$};
%\matrix[style1=0.6cm,right=15pt of 1mat] (2mat)
%{
%	\hspace{-5pt}{\large \m x(k)} \\
%	\hspace{-8pt}\m x(k\hspace{-2pt}+\hspace{-2pt}1) \\
%	\hspace{-8pt}\m x(k\hspace{-2pt}+\hspace{-2pt}2) \\
%	\vdots \\
%	\hspace{-8pt}\m x(k\hspace{-2pt}+\hspace{-2pt}T) \\
%};
%\draw [decorate, decoration={brace,mirror,amplitude=5pt,raise=0.35cm}] (2mat-5-1.west) -- (2mat-5-1.east)  node[midway,yshift=-0.9cm] {$\m x[k]$};

\node at ([xshift=30pt,yshift=-1.2pt]1mat.east) {\Large${\m x[k]=}$};

\matrix[style1=1.3cm,right=55pt of 1mat] (3mat)
{
	\\
	\\
	\\
	\\
	\\
	\\
	\\
	\\
	\\
};
\node[font=\large] at (3mat-1-1) {${\m b^{\mathrm{s}}(k)} $};
\node[font=\large] at (3mat-2-1) {${\m b^{\mathrm{s}}(k\hspace{-3pt}+\hspace{-3pt}1)} $};
\node[font=\large] at (3mat-3-1) {${\m b^{\mathrm{s}}(k\hspace{-3pt}+\hspace{-3pt}2)} $};
\node[font=\large] at (3mat-4-1) {\vdots};
\node[font=\large] at (3mat-5-1) {${\m b^{\mathrm{s}}(k\hspace{-3pt}+\hspace{-3pt}T)} $};
\node[font=\large] at (3mat-6-1) {$\m 0_{n_t\times1} $};
\node[font=\large] at (3mat-7-1) {$\m 0_{n_t\times1}$};
\node[font=\large] at (3mat-8-1) {\vdots};
\node[font=\large] at (3mat-9-1) {$\m 0_{n_t\times1}$};
\draw[dashed]
(3mat-5-1.south west) -- (3mat-5-1.south east);
\draw [decorate, decoration={brace,mirror,amplitude=5pt,raise=0.35cm}] (3mat-9-1.west) -- (3mat-9-1.east)  node[midway,yshift=-0.9cm] {\Large $\m b[k]$};
\end{tikzpicture}
,
}
\end{equation*}
where the submatrices/vectors in $\m A[k]$ and $\m b[k]$ are discussed next. Note that only the current $k$ is considered since submatrices from $k+1$ to $k+T$ have similar structure.

First, the matrices $\m A_{\m h}$, $\m B_{\m q}$, and $\m I_{n_t \times n_t}$ from~\eqref{equ:tankhead-linear} are collected in $\m A^{\mathrm{TK}}(k)$, and the right-hand side for these dynamic constraints is the zero vector of dimension $n_t$. That is, each subvector in $\m b^\mathrm{TK}[k]$ is $\m 0_{n_t \times 1}$.

Second, the equality constraints given in~\eqref{equ:equality} can be collected in~\eqref{equ:measure}, i.e.,  the elevations at reservoirs can be viewed as measurements without any noise. Thus, we only need to deal with~\eqref{equ:measure} and~\eqref{equ:linearHD} in the constraints~\eqref{equ:pse-linear}. After merging~\eqref{equ:measure} and~\eqref{equ:linearHD}, we obtain
%\begin{bmatrix}
%\end{bmatrix}
\setlength\arraycolsep{0pt}
\begin{align}~\label{equ:combin}
&\underbrace{
	\renewcommand{\arraystretch}{0.9}
	\begin{bmatrix}
	\m E(k) \\
	\m C 
	\end{bmatrix}}_{\text{ \normalsize $\m A^\mathrm{s}(k) $}}  \m x(k)\hspace{-3pt}= \hspace{-3pt}    
\underbrace{
	\renewcommand{\arraystretch}{0.8}
	\begin{bmatrix}
	\m z(k) \\
	\m y(k)-\m v(k)
	\end{bmatrix},}_{\text{ \normalsize $\m b^\mathrm{s}(k) $}}
\end{align} 
where $\m A^\mathrm{s}(k) \in \mbb{R}^{(n_y+n_x-n_t-n_r) \times n_x}$ is a matrix collecting and modeling the linear DAEs~\eqref{equ:linearHD}, equality limits~\eqref{equ:equality}, and measurement modeling; vector  $\m b^\mathrm{s}(k)$ is a vector collecting uncertainty source from  demand, pipe, and measurement noise. Appendix~\ref{app:inver} contains more details about this derivation. 
In short, we have now mapped all the constraints inside the covariance operation in~\eqref{equ:pse-linear} into the following problem
\begin{subequations}\label{equ:pse-linear-system-equation}  
	\begin{align} 
	&\mathrm{find} \hspace{3em} \operatorname {K}_{\m x \m x }  \label{equ:pse-linear-obj1} \\
	&\mathrm{s.t.} \hspace{2em}  \Cov(\m A[k]  \m x[k] = \m b[k] ) \label{equ:covaraince-pse-l1}.
	\end{align}
\end{subequations}
The solution to  Problem~\eqref{equ:pse-linear-system-equation} provides the covariance of $\m x[k]$, and the variance of each entry of $\m x[k]$ is located in the diagonal element of $\operatorname {K}_{\m x \m x }$. Two main questions need to be investigated for~\eqref{equ:pse-linear-system-equation}:

\noindent  \textit{(\textbf{Q1})} How is the constraint~\eqref{equ:covaraince-pse-l1} constructed given $\m A[k]$ and $\m b[k]$? \\
\noindent \textit{(\textbf{Q2})} Is Problem \eqref{equ:pse-linear-system-equation} convex and is there a closed-form, unique solution in terms of the problem data?
%{def:cc}

For \textit{(\textbf{Q1})}, the construction of constraint~\eqref{equ:covaraince-pse-l1} is obtained by the  cross-covariance properties $\Cov(\m A \m x, \m A \m x) \hspace{-2pt}= \hspace{-2pt} \m A { \operatorname {K}_{\m x \m x }} \m A^\top$ and  $\Cov(\m b,\m b) = \operatorname {K}_{\m b \m b }$; see Definition~\ref{def:cc} and~\cite{gubner2006probability}. Hence, we can now write 
$$\mathrm{Constraint}~\eqref{equ:covaraince-pse-l1}  \Leftrightarrow \m A[k] {\operatorname {K}_{\m x \m x}} (\m A[k])^\top   = { \operatorname {K}_{\m b \m b }},
$$
where $\operatorname{K}_{\m b \m b}$ is the covariance matrix of all sources of uncertainty which is given next.

Recall that  $\m z(k) = \{\m d(k),  \m K_c^{\mathrm{P}} \m c(k) + \m b^{\mathrm{P}}(k), \m  b^{\mathrm{M}}(k), \m l(k)\}$, and measurement vector ${\m y}(k)$ in~\eqref{equ:combin} is known from sensors, and is replaced by ${\m y^0}(k)$ around the given operating point $\m x^0(k)$.  Vectors $ \m b^{\mathrm{P}}(k)$, $ \m b^{\mathrm{M}}(k)$, and $\m l(k)$ are not random variables, and they are constant since they are all computed around $\m x^0 (k)$, while $\m x(k)$, $\m c(k)$, $\m d(k)$, and $\m v(k)$ are all random. Given Assumption~\ref{assumtion:all}, define the covariance matrix of uncertainty ${\operatorname {K}_{\m b^\mathrm{s} \m b^\mathrm{s} }}$ around $\m x^0 (k)$ as 
	\begin{align}~\label{equ:kbsbs}
	&{ \operatorname {K}_{\m b^\mathrm{s} \m b^\mathrm{s} }(k)} = \Cov\left(	\begin{bmatrix}
	-\m d(k)\\ 
	\m b^{\mathrm{P}}(k) \hspace{-2pt}+\hspace{-2pt} \m  K_c^{\mathrm{P}}  \m c(k)    \\ 
	\m  b^{\mathrm{M}}(k) \\ 
	\m l(k) \\
	\m y^0(k) \hspace{-2pt}-\hspace{-2pt} \m v(k)
	\end{bmatrix}, \begin{bmatrix}
	-\m d(k) \\ 
	\m b^{\mathrm{P}}(k)  \hspace{-2pt}+\hspace{-2pt}  \m  K_c^{\mathrm{P}}  \m c(k)    \\ 
	\m  b^{\mathrm{M}}(k) \\ 
	\m l(k) \\
	\m y^0(k) \hspace{-2pt}-\hspace{-2pt} \m v(k)
	\end{bmatrix}\right) \notag \\
	%	& = 
	%	\Cov\left(	\begin{bmatrix}
	%	\blue{-\m d(k)} \\ 
	%	\m  K_c^{\mathrm{P}}  \m c(k)   \\ 
	%	\m 0_{n_m\times1}  \\ 
	%	\m 0_{n_l\times1} \\
	%	- \m v(k)
	%	\end{bmatrix},	\begin{bmatrix}
	%	\blue{-\m d(k)} \\ 
	%	\m  K_c^{\mathrm{P}}  \m c(k)   \\ 
	%	\m 0_{n_m\times1}  \\ 
	%	\m 0_{n_l\times1} \\
	%	- \m v(k)
	%	\end{bmatrix}\right) \notag \\ 
	& = \begin{array}{c@{\!\!\!}l}
	\diag\left(\left[ \begin{array}[c]{c}
	\Var(\m d(k)) \\ 
	(\m  K_c^{\mathrm{P}})^2  \Var(\m c(k))  \\ 
	\m 0_{n_m\times1}  \\ 
	\m 0_{n_l\times1} \\
	\Var(\m v(k))
	\end{array}  \right]\right)
	&
	\begin{array}[c]{@{}l@{\,}l}
	\hspace{2pt} \leftarrow \text{demand uncertainty} \\
	\hspace{2pt} \leftarrow \text{pipe uncertainty} \\
	\\
	\\
	\hspace{2pt} \leftarrow \text{noise uncertainty.} \\
	\end{array}
	\end{array}
	\end{align}
We note that 
\newline $\Cov(\m  K_c^{\mathrm{P}}  \m c(k),\m  K_c^{\mathrm{P}}  \m c(k)  ) =  \m  K_c^{\mathrm{P}}  \Cov(\m c(k),\m c(k)) (\m  K_c^{\mathrm{P}})^\top $ 
\newline $\hspace*{3.6cm}= \diag{\left((\m  K_c^{\mathrm{P}})^2  \Var(\m c(k))\right)}$  
\newline since $\m  K_c^{\mathrm{P}}$ and $\Cov(\m c(k),\m c(k))$ are diagonal. Note that the covariance of $\m b^\mathrm{TK}[k]$ is a zero matrix. If we consider all time-steps then we can obtain the aggregate covariance matrix of all uncertainty sources as
\begin{align}~\label{equ:kbb}
{ \operatorname {K}_{\m b \m b }} &=\diag{\left(\operatorname {K}_{\m b^\mathrm{s} \m b^\mathrm{s} }[k],\operatorname {K}_{\m b^\mathrm{TK} \m b^\mathrm{TK} }[k]\right)} \notag \\
&=\diag{\left(\operatorname {K}_{\m b^\mathrm{s} \m b^\mathrm{s} }(k),\dots,\operatorname {K}_{\m b^\mathrm{s} \m b^\mathrm{s}} (k+T),\m 0\right)}.
\end{align}
This shows how the covariance operator of the linear equality constraints can be constructed. 
\subsection{Explicit solution for {PSE} and a real-time algorithm} \label{sec:linmatrix}
As for \textit{(\textbf{Q2})}, Problem~\eqref{equ:pse-linear-system-equation} now becomes a feasibility problem and is convex since the constraint only consists of linear equality constraints. Here, we present explicit solutions to the PSE problem, thereby addressing the second part of \textit{(\textbf{Q2})}.
To solve~\eqref{equ:pse-linear-system-equation},  we can compute the covariance matrix through solving
\begin{align}\label{equ:pse-matrix}
\m A[k] { \operatorname {K}_{\m x \m x }} (\m A[k])^\top   = { \operatorname {K}_{\m b \m b }}.
\end{align}
The next result provides a solution to~\eqref{equ:pse-matrix} with a guarantee on an existence of a unique solution.
\begin{thm}~\label{thm:inver}
	Matrix $\m A[k]:=\m A$ is full rank, and the solution to~\eqref{equ:pse-matrix} is uniquely given by 
	\begin{align}\label{equ:pse-solution}
	{ \operatorname {K}_{\m x \m x }} = (\m A^\top \m A)^{-1} { \m A^\top \operatorname {K}_{\m b \m b }} \m A (\m A^\top \m A)^{-1}.
	\end{align}
\end{thm}
The proof of Theorem~\ref{thm:inver} is given in Appendix~\ref{app:inver}. The physical meaning of~\eqref{equ:pse-solution} indicates that uncertainty of system state $\m x(k)$ can be expressed through a linear combination of the uncertainty from demand ${\m d(k)}$,  pipe parameters $ \m c(k)$, and measurement noise $\m v(k)$. That is,  the covariance matrix of system states ${\operatorname {K}_{\m x \m x }}$ is expressed by the linear combination of covariance matrix of uncertainty $\operatorname {K}_{\m b \m b}$. Note that the above discussion on \textit{(\textbf{Q2})} is similar as the observability analysis identifying if a set of available measurements is sufficient to estimate the system state; see~\cite{Pruneda2010,Diaz2016a} for details. 
\begin{algorithm}[t]
	\small	\DontPrintSemicolon
	\KwIn{ WDN topology, expectation of measurement  $\m y(k)$, demand  $\{\m d(k)\}_{k=1}^{T}$, variance of measurement noise $\{\Var{(\m v(k))}\}_{k=1}^{T}$,   demand $\{\Var{(\m d(k))}\}_{k=1}^{T}$, pipe parameter $\{\Var{(\m c(k))}\}_{k=1}^{T}$}
	\KwOut{The estimated covariance matrix of $\m x$: $\{ { \operatorname {K}_{\m x \m x }} \}_{k=1}^{T}$ }
	%Convert the variable $\langle{{\m x}}\rangle_{0}$  to its GP form $\langle{\hat{\m x}}\rangle_{0}$\;
	Set  $k=1$ \;
	\While {  $k \leq T$ }{
		Obtain an operating point ${\m x^0}$\;
		Linearize the nonlinear DAE around $\m x^0$ via~\eqref{equ:linearPipesPumps}\;
		Compute ${ \operatorname {K}_{\m x \m x }}$ through ~\eqref{equ:pse-solution}\;
		Extract the variances of heads and flows from the diagonal elements of  ${ \operatorname {K}_{\m x \m x }}$\;
		Update $k = k + 1 $
	}
	\caption{ Estimating the covariance matrix of $\m x$.}
	\label{alg:alg2}
\end{algorithm}

Given the above discussions, Algorithm~\ref{alg:alg2} presents a near real-time algorithm to solve the PSE in water networks.  That is, the algorithm solves the PSE for the expectation and covariance of the heads and flows assuming demand, pipe roughness coefficients, and measurement follow any statistical distribution, while knowing the expectation and variance for the uncertain variables.
% Vector ${\m x^0}$ is viewed as the operating point for~\eqref{equ:pse-linear-system-equation} or its matrix form~\eqref{equ:pse-matrix}. 

We  note that the \textit{system matrix} $\m A[k]$ can have more equations than variables---in the case of over-determined measurement scenario. Under such scenario, some equations in~\eqref{equ:pse-matrix} might be more \textit{reliable} than their counterparts. For example, the demand predictions that are encoded in~\eqref{equ:pse-matrix} might be less reliable than the flow and head measurements.  As a result, the weighted version can be rewritten as
\begin{align}\label{equ:weighted-hy-matrix}
\underbrace{\m A^\top \m W \m A}_{\text{ \normalsize $\m A^{\m W}[k]$}} \m x[k] = 	\underbrace{\m A^\top \m W \m b[k],}_{\text{ \normalsize $\m b_1$[k]}}
\end{align}
where $\m W$ is a diagonal, positive-definite weight matrix. The corresponding weighted version of~\eqref{equ:pse-matrix} is
\begin{align}\label{equ:weighted-pse-matrix}
\m A^{\m W}[k] { \operatorname {K}_{\m x \m x }} (\m A^{\m W}[k])^\top   = { \operatorname {K}_{\m b_1 \m b_1 }}.
\end{align}
\begin{corollary}~\label{cor:inver} 
The weighted formulation~\eqref{equ:weighted-pse-matrix} with full-rank matrix $\m A^{\m W}$ has a unique solution given by
\begin{equation} \label{equ:cce1}
\operatorname {K}_{\m x \m x}= ({\m A^{\m W}}^\top \m A^{\m W})^{-1}  {\m A^{\m W}}^\top \operatorname {K}_{\m b \m b } \m A^{\m W} ({\m A^{\m W}}^\top \m A^{\m W})^{-1}
\end{equation}
\end{corollary}

Corollary~\ref{cor:inver} builds on Theorem~\ref{thm:inver} with nearly identical proof; see Appendix~\ref{app:inver}. We note that~\eqref{equ:weighted-pse-matrix} is a simple yet powerful formula to calculate the covariance matrix of system states. Intuitively, we know that the state uncertainty in a system is reduced when more measurements on system states are performed. For example, if all states are measured and corresponding equations are given larger weights, then the uncertainty is only introduced from  the noise of sensors, and the other uncertainty sources are suppressed. Specifically, when more head and flow measurements are obtained, the impacts are not only on the measured states themselves but also on the other unmeasured states due to the system equation. 

Algorithm~\ref{alg:alg2} suits both sufficient and over-determined scenarios with various types of valves when the corresponding formulations are adopted; see Section~\ref{sec:SMSODS}.
We note here that~\eqref{equ:pse-solution} and~\eqref{equ:cce1} do \textit{not} rely on any statistical distribution of the uncertainty. Furthermore, this approach for PSE is scalable to networks with thousands of nodes due to the fact that  \textit{system matrix} $\m A[k]$ is highly sparse. Section~\ref{sec:Test} produces examples on the scalability of this algorithm.
\begin{remark}
	Equation~\eqref{equ:pse-matrix} is not based on the first-order second moment (FOSM) method, since  FOSM uses the \textit{sensitivity matrix} which is the Jacobian matrix and only focuses on the relationship between part of system states and uncertainty source~\cite{kang2009approximate,xu1998probabilistic,lansey2001calibration}. {Moreover, FOSM needs to be applied two or more times to find the mapping between all system states and uncertainty}, whereas~\eqref{equ:pse-matrix} uses \textit{system matrix} $\m A[k]$ and our formulation connects between all system states and all uncertainty sources in one shot.
\end{remark}

\section{Generalization of Flows and Heads Distributions}~\label{sec:distribution}
In Section~\ref{sec:formulation}, we showcase that~\eqref{equ:pse-solution} computes the variance of each system state efficiently, but it still cannot explain how uncertainty propagates to the system states clearly and intuitively for the linearized or the original, nonlinear hydraulic modeling. That is, \textit{how would the distributions of unmeasured flows and heads change with the change in the uncertainty distribution?} Furthermore, and albeit useful,  \eqref{equ:pse-solution} cannot generate insights related to which uncertainty source is more critical than the other in terms of SE performance. 
In this section, in order to answer the above questions, we need a further assumption that builds on Assumption~\ref{assumtion:all}.
\begin{assumption}\label{assumtion:all2}
	Demand $\m d$, measurement noise $\m v$, and  pipe roughness coefficient $\m c$ follow either normal {or} uniform distributions. 
\end{assumption}
This assumption is not needed to run Algorithm~\ref{alg:alg2}; rather, it is given to investigate how uncertainty propagates for specific distributions. With that in mind,
Section~\ref{sec:uncertainty} presents a thorough discussion on the practicality of Assumption~\ref{assumtion:all2}. The next theorem and corollary explain how distributions of flows and heads change with the uncertainty distribution for linearized or nonlinear hydraulic models.
\begin{lemma}\label{lemma:treeNonlinear}
	For tree networks with nonlinear hydraulic models, flows are normally (uniformly) distributed as long as the nodal demands follow the normal (uniform) distribution, whereas the heads are never normally (uniformly) distributed even if all uncertainty sources are normally (uniformly) distributed. Furthermore, for looped networks with nonlinear hydraulics, neither heads nor flows are normally (uniformly) distributed. 
\end{lemma}
\begin{thm}~\label{thm:looped}
	For tree or looped water networks with linearized hydraulics, all heads and flows follow the normal (uniform) distribution  when uncertainty is (uniformly) distributed.
\end{thm}

The proof of Lemma~\ref{lemma:treeNonlinear} and Theorem~\ref{thm:looped}  are  given in Appendix~\ref{sec:proof} along with the corresponding PDFs for flows and heads for linearized head loss models. The proof explains the propagation of uncertainty and the importance of each uncertainty source. The next discussion summarizes and explains these findings.

\begin{figure}[t]
	\centering
	\includegraphics[width=0.65\linewidth]{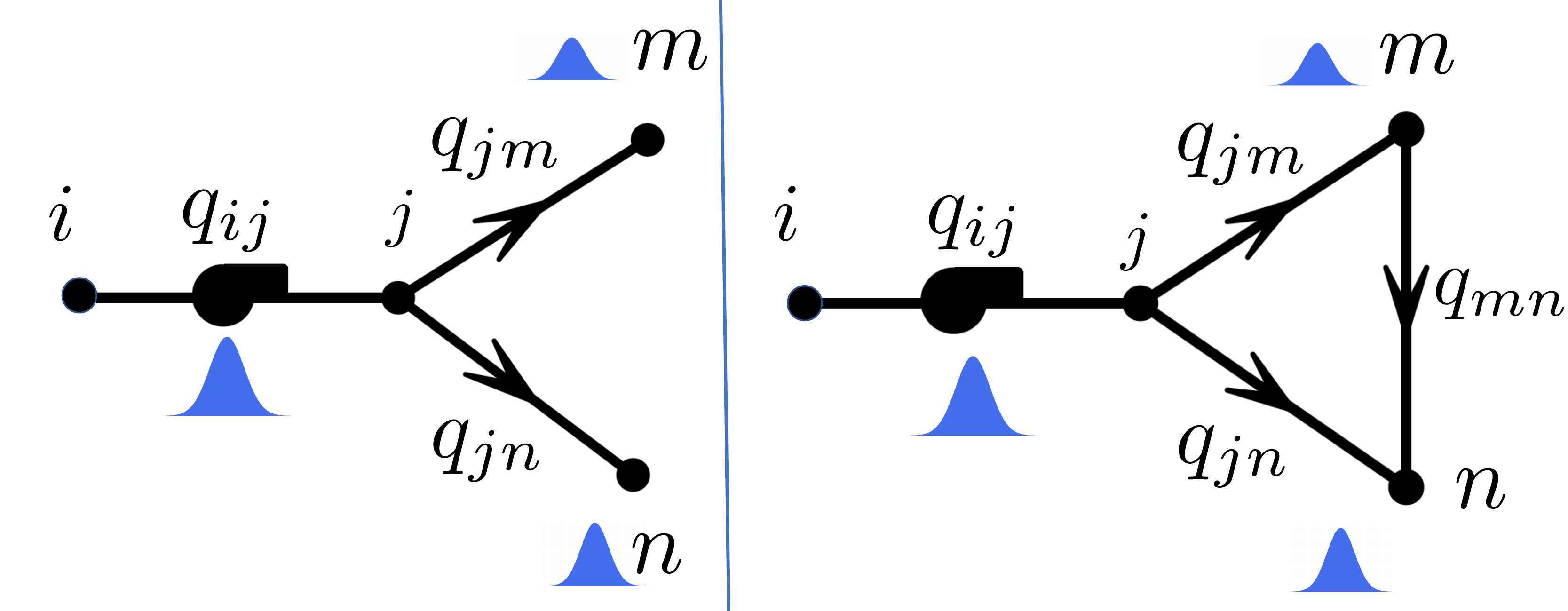}
	\caption{A tree network (left), and a looped network (right). }
	\label{fig:tree}
\end{figure}

\begin{figure}[t]
	\centering
	\includegraphics[width=\linewidth]{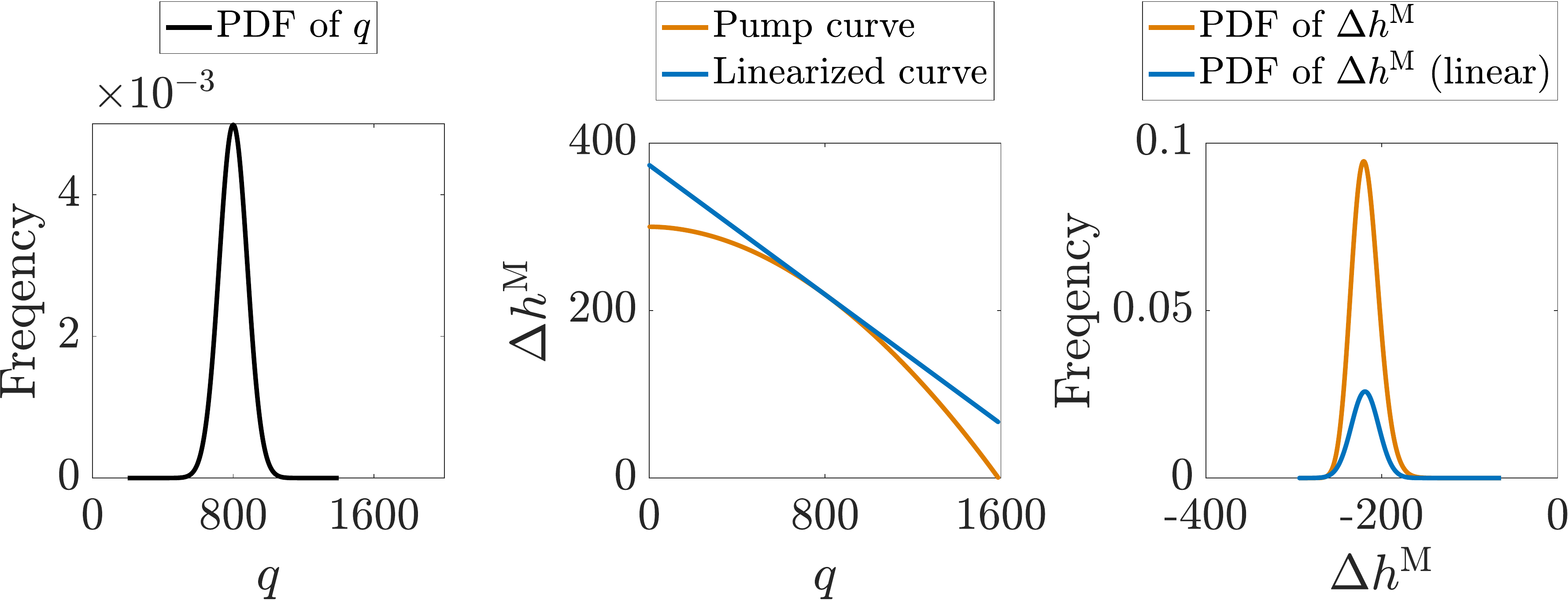}
	\caption{ PDF of the flow $q$ through a pump (left), original head increase curve of the pump and its linear form (middle), and  PDF of original and linearized head increase (right).}
	\label{fig:pump}
\end{figure}

\begin{figure}[t]
	\centering
	\includegraphics[width=\linewidth]{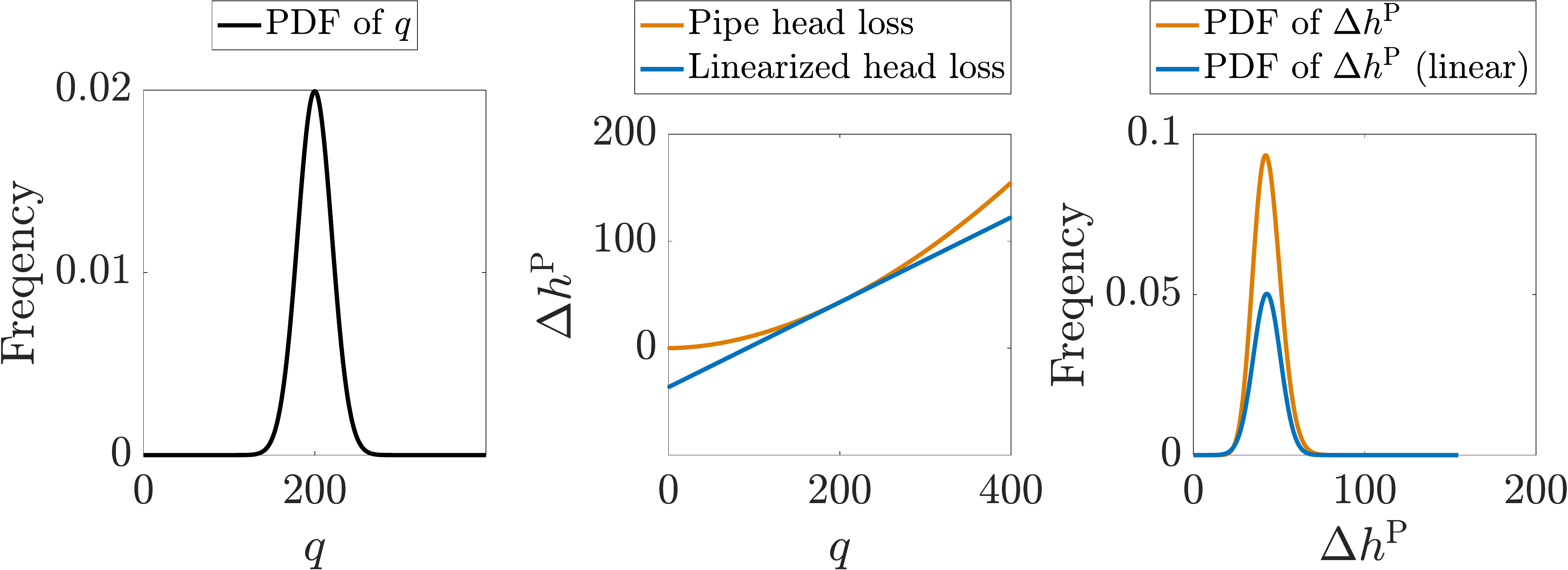}
	\caption{PDF of flow $q$ through a pipe (left), original head loss curve of the pipe and its linear form based only on demand (middle), and PDF of original and linearized head loss base only on demand (right).}
	\label{fig:pipe}
	% \vspace{-10pt}
\end{figure}

In the proof of Lemma~\ref{lemma:treeNonlinear}, we show the PDF of head gain of a pump $\Delta h_{ij}^{\mathrm{M}}$and head loss of a pipe $\Delta h_{ij}^{\mathrm{P}}$ in tree network with nonlinear hydraulic modeling, given by $f_{\Delta H_{ij}^{\mathrm{M}}}(\Delta h_{ij}^{\mathrm{M}})$~\eqref{equ:pdf_pump_new}  and $f_{\Delta H_{ij}^{\mathrm{P}}}(\Delta h_{ij}^{\mathrm{P}})$~\eqref{equ:pdf_pipe_new}.  By comparing them with the  formula of normal distribution function,  and we know that they are not normally distributed. We also depict these conclusions in Figs.~\ref{fig:pump} and~\ref{fig:pipe}, which are based on simple numerical tests for a network shown in Fig.~\ref{fig:tree}. The nonlinear pump curve and pipe head loss curve are the yellow lines in the middle of  Figs.~\ref{fig:pump} and~\ref{fig:pipe}, whereas the shape of $f_{\Delta H_{ij}^{\mathrm{M}}}(\Delta h_{ij}^{\mathrm{M}})$ and $f_{\Delta H_{ij}^{\mathrm{P}}}(\Delta h_{ij}^{\mathrm{P}})$ are  yellow lines on Figs.~\ref{fig:pump} and~\ref{fig:pipe} (right) when the PDF of flow $q$ is  normally distributed.

\begin{figure}[t]
	\centering
	\includegraphics[width=0.80\linewidth]{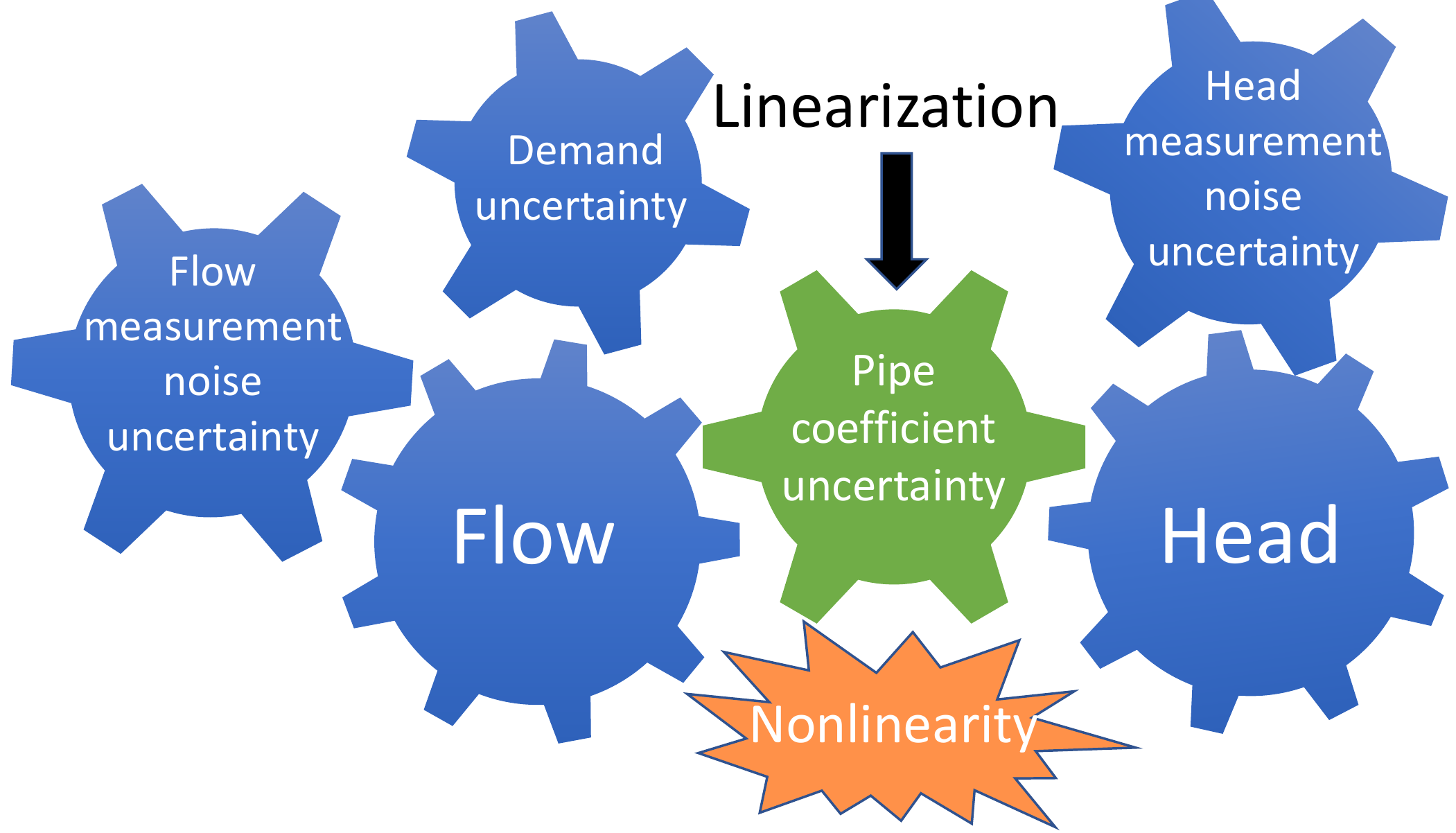}
	\caption{Propagation of uncertainty.}
	\label{fig:poof}
\end{figure}

Based on these results, we conclude that \textit{(i)} the normality of demand uncertainty results in the normality of flows in tree networks regardless of nonlinear or linear hydraulic modelings. Furthermore, the uncertainty from pipe roughness coefficients prevents the propagation of normality from flows to heads. \textit{(ii)} The normality of noise uncertainty only impacts the variables it measures, i.e., the flow (head) measurement noise uncertainty only impacts flows (heads). The linearization bridges the gap between the normality of flows and heads, which makes the system states follow normal distributions. These discussions are abstractly depicted in Fig.~\ref{fig:poof}. From the discussions, the importance or impact of each uncertainty source is given in the next remark.
\begin{remark}\label{remark:impact}
	Impact of pipe roughness parameter uncertainty  $>$  Impact of demand uncertainty  $>$  Impact of measurement noise uncertainty.% 
\end{remark}
Remark~\ref{remark:impact} is discussed in breadth in the case studies section.
\section{Discussions and Insights }~\label{sec:Discussion}
This section presents some discussions that are relevant to uncertainty and its distribution,  different measurement scenarios, and relevance and connections to the literature.
\subsection{Dealing with uncertainty}~\label{sec:uncertainty}
The rationale behind Assumption~\ref{assumtion:all2} is discussed here even though Algorithm~\ref{alg:alg2} can be implemented regardless of the uncertainty distribution type.
For the uncertainty from  measurement noise,  we simply assume that they follow a known normal distribution~\cite{andersen2001constrained}. A uniform distribution can also be applicable in our paper as we mentioned. As for the uncertainty for pipe roughness coefficients, and specifically for the Hazen-Williams coefficients, typical values lie in $[75,130]$~\cite{mohapatra2014distribution}, which is narrow enough, and can be fully covered by a normal distribution $\mathcal{N}(100,135.5)$. This means the worse case $[75,130]$ is covered by $\mathcal{N}(100,135.5)$ with $99\%$ confidence level. If the material of pipes is known, the range can be further narrowed down. In this way, we can deal with the uncertainty from roughness coefficient. In fact, the network calibration research~\cite{ormsbee1989implicit} is able to estimate the roughness coefficient, and recently, the standard deviation errors for the estimated pipe roughness are available~\cite{ostfeld2011battle}.

For demand uncertainty, various methods~\cite{guercio1970instantaneous,lopez2018multi,guo2018short} can be adopted to estimate the demand (mean and variance). For example,  the authors in~\cite{guo2018short} propose a gated recurrent unit network model to predict short-term water demand, and the histogram of relative errors indicates that $95\%$ of the forecast relative errors fall within the range of $\pm$12.65\% for 24-hour forecasts. In addition to these prediction studies, demand estimation studies based on the historical and real-time measurements (a calibration process using sensors), such as~\cite{shang2008real,andersen2000implicit}, the knowledge of population densities, and billing data from water utilities are also helpful to deal with demand uncertainty.

\begin{figure}[t]
	\centering
	\includegraphics[width=0.75\linewidth]{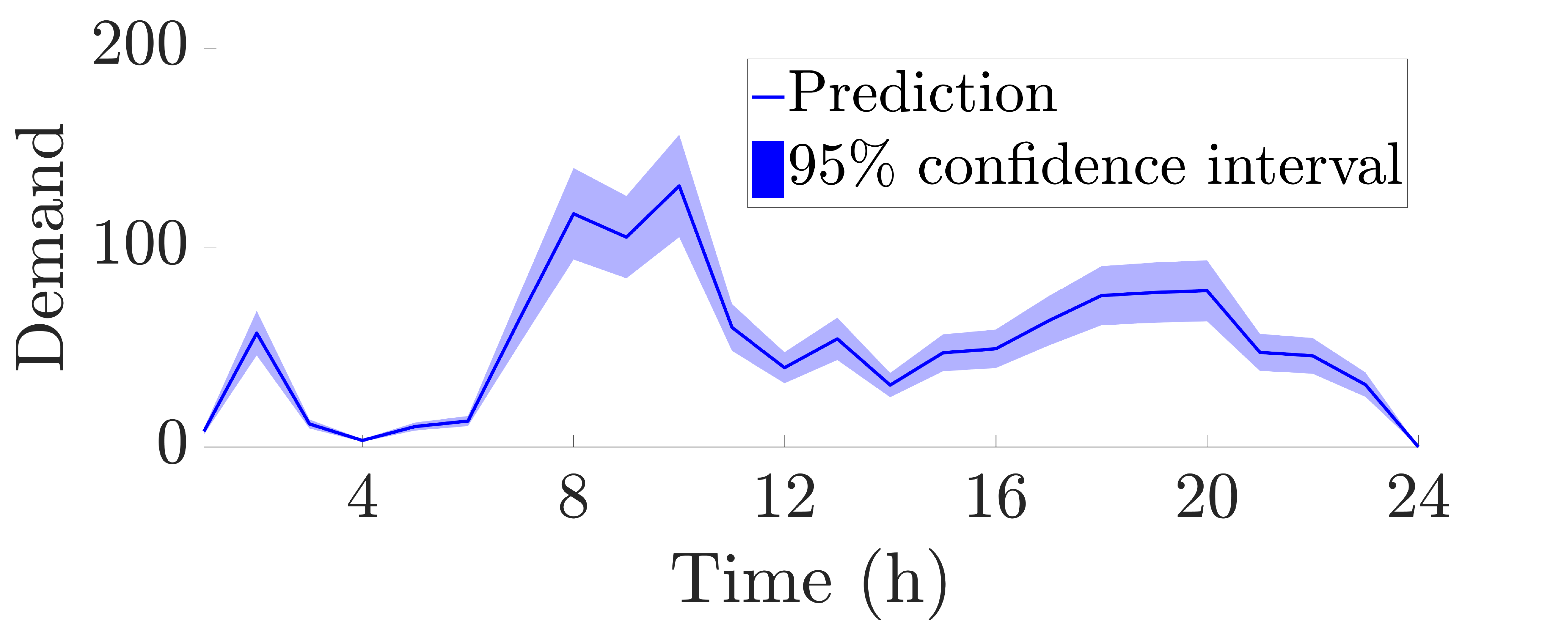}
	\caption{Demand prediction with 95\% confidence interval.}
	\label{fig:confidentlevel}
\end{figure}

Thus, given above studies, prediction or estimation error of demand can be assumed to follow a certain type of distributions, i.e., a uniform distribution or a normal distribution and vary around a prediction value (mean) with a standard deviation~\cite{andersen2001constrained}. Fig.~\ref{fig:confidentlevel}  illustrates a 24-hour estimated demand and the associated $95\%$ confidence interval. 
\subsection{Sufficient and over-determined scenarios}\label{sec:SMSODS}
There are two scenarios for SE in water networks~\cite{powell1999state,wang2019state}: \textit{(i)} Sufficient scenario is described by having the same number of equations and unknowns.  When measurements include the heads at tanks and reservoirs, and demand pseudo-measurements are available for state estimation, we consider it as a sufficient scenario. In fact, the DSE problem under sufficient scenario is similar to the water flow problem (WFP)~\cite{singh2019flow} or reliability  analysis~\cite{xu1998probabilistic}. \textit{(ii)} Over-determined scenario has more equations than unknowns. For example,  additional sets of heads are measured at several key nodes besides the head at tanks and reservoirs. For more details of solving deterministic SE (DSE) under both scenarios, please refer to the examples in~\cite{powell1999state,wang2019state}. For PSE, we note that as the number of extra measurements increases, the uncertainty of system variables gradually disappears, and SE accuracy is enhanced. 
\section{Case Studies}~\label{sec:Test}
We present several simulation examples (illustrative three-node network, 8-node network~\cite{rossman2000epanet}, Anytown, BAK, PESCARA, OBCL,  and D-Town~\cite{Eliades2016}) to illustrate the applicability of our approach. The first three-node network is used to illustrate the details  of proposed method. Then we test the 8-node network to illustrate that our approach can deal with the looped topology and various types of valves. The rest of testcases are used to test the scalability and the efficacy of the proposed approach for non-Gaussian uncertainty distributions. All test cases are simulated using MCS methods via EPANET Matlab Toolkit~\cite{Eliades2016} on Ubuntu 16.04.4 LTS with an Intel(R) Xeon(R) CPU E5-1620 v3 @ 3.50 GHz, and results of MCS are used to verify the accuracy of proposed approach.  
All codes, parameters, tested networks, and results are available on Github~\cite{shenwangSE}.

In order to compare our solution to MCS (with 1000 randomization), we need to define the criteria at first. The absolute error between $\m \sigma_{\mathrm{MCS}}$ and $\sigma_{\mathrm{PSE}}$ is defined as $\m {\mathrm{AE}} = |\m \sigma_{\mathrm{MCS}} - \m \sigma_{\mathrm{PSE}}|$, and  the corresponding relative error is $\m {\mathrm{RE}}=\frac{\m {\mathrm{AE}}}{|\m \sigma_{\mathrm{MCS}}|} \times 100\%$, where the entries in $\m \sigma$ are standard deviations.

Assumption~\ref{assumtion:all2} (only normality) are used in case studies, and there are three equivalent ways to express the uncertainty of a normally distributed random variable $x$. First, the direct way is using distribution $x \sim \mathcal{N}(\mu_x,\,\sigma_x^2)$ which indicates the $\E(x) = \mu_x$ and $\Var(x) = \sigma_x^2$. The second one is given by $\E(x) = \mu_x$ and the corresponding margin of error (ME) which is the  percentage of relative changes deviating from the $\E(x)$ under a certain confidence interval. That is, the $x_{a\%}$ notation defines that $\Delta x$ falls into the range of $\pm$a\% around value $\E(x)$ under  a 99\%  confidence interval, and ME of $x$ is calculated by $a = 2.576\frac{\sigma_x}{\mu_x} \times 100\%$. For 95\% confidence interval, $a = 1.96\frac{\sigma_x}{\mu_x} \times 100\%$. Quantity $x_{0\%}$ indicates that $x$ is not a random variable. Third, the uncertainty can also be described via a range given by $\E(x)$ and $\Delta x$ under a certain confidence interval, that is $x \in [\E(x) - \Delta x, \E(x) + \Delta x]$.

\subsection{Three-node network}\label{sec:3-node}
\begin{figure}[t]
	\centering
	\includegraphics[width=0.34\linewidth]{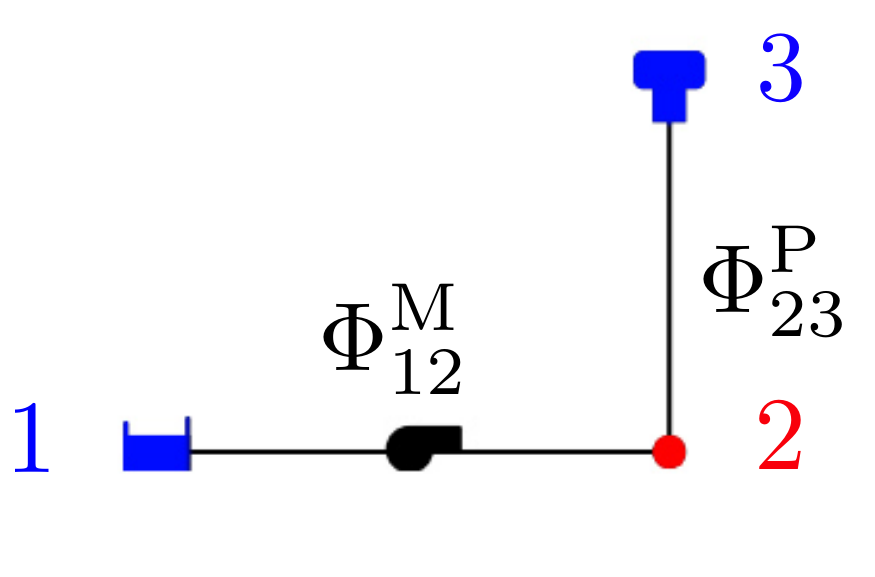}
	\hspace{5em}
	\includegraphics[width=0.34\linewidth]{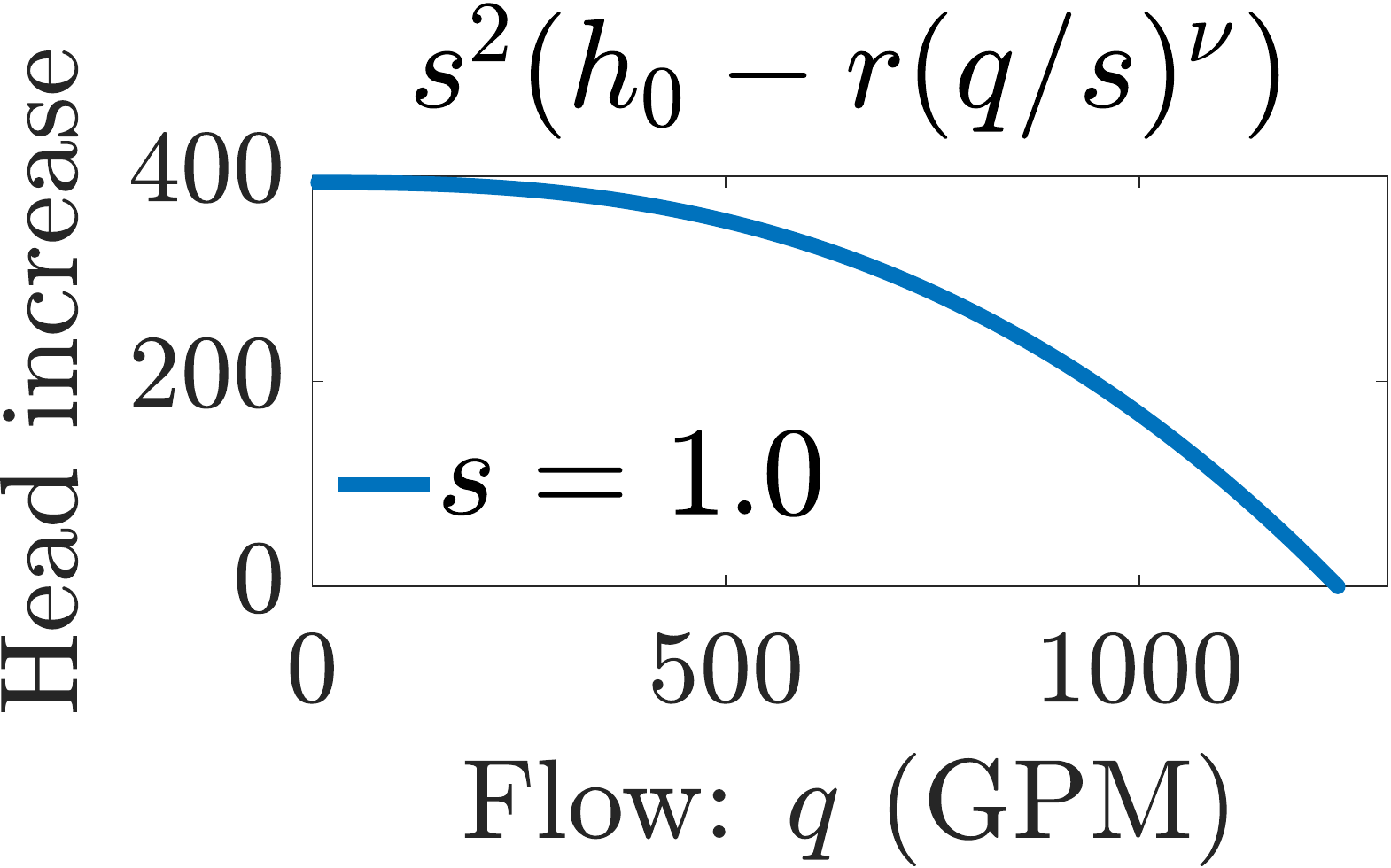}
	\caption{Three-node network(left) and its variable-speed pump curve (right). The $\m \Phi(\cdot)$ function denotes  nonlinear head loss or head gain. }
	\label{fig:setup}
\end{figure}

The three-node network includes one junction, one tank, one pipe, one pump, and one reservoir, and is shown in Fig.~\ref{fig:setup} (left). The corresponding pump curve is  shown in Fig.~\ref{fig:setup} (right). Only Junction 2 consumes water, and we assume that estimated demand is $\mu_{d_2} = 100\ \mathrm{GPM}$, and the uncertainty sources and corresponding values are summarized in Table~\ref{tab:interval} under 99\% confidence interval. For example, 99\% of the estimated demand relative errors fall within the range of $\pm$20\% around estimated value $\mu_{d_2} = 100 \mathrm{GPM}$ (denoted as $d_{20\%}$), or $d_2\sim \mathcal{N}(\mu_{d_2},\,\sigma_{d_2}^2)$ and $\sigma_{d_2} = \frac{20\% \mu_{d_2}}{2.576} = 7.7$.  Similarly, we assume the uncertainty of Hazen-Williams coefficient $c$ for Pipe 23 and measurement noise $v$ for Tank 3 are in $c_{20\%}$ and $v_{1\%}$ situation when measurements of head $\textcolor{black}{y_{h_3} = 908\ \mathrm{ft}}$ with $\sigma^2_{{h_3}} = 0.0241$.  The head at Reservoir 1 equals to its elevation which is  treated as a constant without any uncertainty. that is, a equality limit $y_{h_1} = 700\ \mathrm{ft}$ with $\sigma^2_{{h_1}} = 0$. 

As mentioned, state estimation contains two typical scenarios: sufficient and over-determined scenarios. The solution of DSE under sufficient (over-determined) scenario is served as operating point $\m x^0$ for PSE  under sufficient (over-determined) scenario. We note that the sufficient measurement scenario corresponds to the traditional water flow problem~\cite{singh2019flow} which can be solved by the standard water system simulation software \textit{EPANET}. That is, we can use EPANET to verify our solutions under the sufficient scenario even though EPANET cannot solve state estimation problems.

\begin{table}[t]
	\caption{Setup of uncertainty (99\% confidence interval).}
	\renewcommand{\arraystretch}{1.12}
	\resizebox{\linewidth}{!}{%
		\begin{tabular}{c|c|c|c|c}
			\hline
			\textit{\makecell{Uncertainty\\source}} & \textit{\makecell{Estimated \\ value }} & \textit{\makecell{ME\\ (\%)}}  & \textit{\makecell{Distribution\\ $\mathcal{N}(\mu,\,\sigma^2)$}} & \textit{\makecell{Range}} \\ \hline
			\textit{Demand} $d$& 100 $\mathrm{{\small GPM}}$  &  $20\%$ &   $\mathcal{N}(100,60.28)$ & $[80,120]$\\ \hline
			\textit{Roughness} $c$&  100 &  $20\%$ & $\mathcal{N}(100,60.28)$ &    $[80,120]$\\ \hline
			\textit{\makecell{Noise}} $v$ &  0 $\mathrm{ft}$ &  $1\%$& $\mathcal{N}(0,0.0241)$ & $[-0.04,0.04] $  \\ \hline \hline
		\end{tabular}%
	}
	\label{tab:interval}
	% \vspace{-0.5cm}
\end{table}
%\begin{table}[t]
%	\caption{Formulations of three-node network  when $T = 1$.}
%	\setlength\arraycolsep{1pt}
%	\setlength{\tabcolsep}{0.5pt}
%	\renewcommand{\arraystretch}{0.8}
%	\resizebox{\linewidth}{!}{%
%	\begin{tabular}{c}
%		\hline
%		\parbox{8.2cm}{
%				% \vspace{-0.7em}
%				\begin{align}
%				{\small \text{DSE:}} \hspace{1em} &   {\large \min}\;\; \hspace{6em} f(\m x(k)) = \m v^\top \m W \m v \label{equ:dse-classical}  \\
%				&\mathrm{s.t.}\;\; q_{12} - q_{23} =d_2;\ \  h_{1} - h_{2} =\Phi_{12}^{\mathrm{M}};\ \  h_{2} - h_{3} = \Phi_{23}^{\mathrm{P}}\notag
%				\end{align}
%				% \vspace{-1.0em}
%		}     \\ \hline
%\hspace{-10pt}
%		\parbox{8.2cm}{	
%				% \vspace{0.2em}
%			\hspace{-0em} {\small \text{System Equation:}}
%			% \vspace{-0.4em}
%			\begingroup
%			
%			% \vspace{-0.7em}
%			\endgroup
%		}          \\ 
%   \\ \hline
%		\hline
%		\multicolumn{1}{l}{}
%	\end{tabular}%
%	\label{tab:3node-model}
%	}
%	% \vspace{-0.5em}
%\end{table}
\subsubsection{Sufficient measurements scenario} \label{sec:sufficient} Supposing that we only know the measurements of head $\textcolor{black}{y_{h_1}}$ and $\textcolor{black}{y_{h_3}}$  for Reservoir 1 and Tank 3.  The operating point $\m x^0$ is obtained first using the authors' approach in~\cite{wang2019state}. The approach is based on linear approximations of the nonlinear state estimation problem under no uncertainty. This linear approximation returns nearly identical solutions to {EPANET}.  
  
Next, the  corresponding PSE  is solved to obtain  covariance matrix ${ \operatorname {K}_{\m x \m x }}$. The linear system of equations for three-node network when $T=1$ based on~\eqref{equ:combin} can be expressed as\footnote{Measurement formulation under sufficient scenario is in \textcolor{blue}{blue};  Formulation of extra measurement is in \textcolor{red}{red},  and over-determined scenario consider both \textcolor{blue}{blue} and  \textcolor{red}{red} formulation.The tank dynamic equations are not listed due to $T = 1$, hence $\m A = \m A^{\mathrm{s}}$ and $\m b = \m b^{\mathrm{s}}$ in this case.}
{\small \begin{align} ~\label{eq:3node-LP-matrix}
	\setlength\arraycolsep{3pt}
	\renewcommand{\arraystretch}{1.05}
	\underbrace{\begin{bmatrix}
		0 & 0 &0 &-1& 1\\
		1 &0& -1& \tiny{-k_{23}^{\mathrm{P}}}& 0\\
		-1& 1& 0& 0& \tiny{k_{12}^{\mathrm{M}}}\\
		\textcolor{blue}{0} &\textcolor{blue}{0}& \textcolor{blue}{1}& \textcolor{blue}{0}& \textcolor{blue}{0} \\
		\textcolor{blue}{0} &\textcolor{blue}{1}& \textcolor{blue}{0}& \textcolor{blue}{0}& \textcolor{blue}{0} \\
		\textcolor{red}{1} &\textcolor{red}{0}& \textcolor{red}{0}& \textcolor{red}{0}& \textcolor{red}{0}
		\end{bmatrix}}_{\textstyle \m A} \hspace{-3.5pt}
	\underbrace{\begin{bmatrix}
		h_2\\
		h_1 \\
		h_3\\
		q_{23}\\
		q_{12}
		\end{bmatrix}}_{\textstyle \m x}& =
	\renewcommand{\arraystretch}{1.05}
	\underbrace{\begin{bmatrix}
		d_2\\
		{b}_{12}^{\mathrm{P}}\hspace{-1pt} + \hspace{-1pt}  k_{c}^{\mathrm{P}} c_{23}\\
		{b}_{12}^{\mathrm{M}} \\
		\textcolor{blue}{y_{h_3}} - \textcolor{blue}{v_{h_3}}\\
		\textcolor{blue}{y_{h_1}} - \textcolor{blue}{v_{h_1}}\\
		\textcolor{red}{y_{h_2}} - \textcolor{red}{v_{h_2}}
		\end{bmatrix}.}_{\textstyle \m b} 
	\end{align}
}
Based on~\eqref{equ:kbb} and Table~\ref{tab:interval}, the uncertainty of demand, pipe roughness coefficient, and measurement noise is $\operatorname {K}_{ {\m b}  {\m b}} = \diag([60.28, (k_{c}^{\mathrm{P}})^2 60.28, 0, 0.00241,0])$.
 Thus, we can solve the covariance matrix $\operatorname {K}_{ {\m x}  {\m x}}$ by~\eqref{equ:pse-solution} of Theorem~\ref{thm:inver}, and the solution under sufficient scenario is presented in Table~\ref{tab:3node-result-covariance} and compared with MCS. Notice that only the entries for $q_{12}$, $q_{23}$, and $h_2$ in covariance matrix $\operatorname {K}_{ {\m x}  {\m x}}$  are listed since $h_1$ and $h_3$ are measurements.

The variance of  $q_{12}$, $q_{23}$, and $h_2$ are $0.16$, $55.44$, and $0.6$ from Table~\ref{tab:3node-result-covariance}. We note that \textit{(i)} the uncertainty introduced from demand $d_2$, pipe roughness coefficient $c_{12}$, and measurement noise $v_{h3}$ are mainly passed to $q_{23}$, \textit{(ii)} the solution from PSE are close to the one from MCS which confirms the effectiveness of proposed approach.
\begin{table}[t]
	\setlength{\tabcolsep}{3pt}
	\renewcommand{\arraystretch}{1.12}
	\centering
	\caption{Covariance matrix for $q_{12}$, $q_{23}$, and $h_2^*$.}
	\begin{tabular}{c|ccc|ccc|ccc}
		\hline
		& \multicolumn{3}{c|}{\textit{\makecell{MCS}}} & \multicolumn{3}{c|}{\textit{\makecell{Sufficient Scenario \\  by Theorem~\ref{thm:inver}}}} & \multicolumn{3}{c}{\textit{\makecell{Over-determined Sce- \\nario by Corollary~\ref{cor:inver}}}} \\ \hline
		& $q_{12}$     & $q_{23}$   & $h_{2}$  & $q_{12}$     & $q_{23}$   & $h_{2}$ & $q_{12}$     & $q_{23}$   & $h_{2}$   \\ \hline
		$q_{12}$   & 0.17 &	0.18 & 	-0.34 &  0.16 &	0.20 &	-0.31 & 0.00 &	-0.07&	0.00 \\ 
		$q_{23}$  & 0.18 & 	55.46 &	-0.35  & 0.20 &	55.44 &	-0.39 & -0.07&	11.50 &	0.03 \\
		$h_{2}$ &-0.34 &	-0.35& 0.65 & -0.31	& -0.39 &	0.60& 0.00&	0.03&	0.00 \\
		\hline \hline
		\multicolumn{10}{l}{\footnotesize{ \makecell{$*$Only result in sufficient scenario is compared with MCS since MCS \\cannot provide solution for over-determined scenario.}}}
	\end{tabular}
	\label{tab:3node-result-covariance}
\end{table}

\begin{figure}[t]
	\centering
	\includegraphics[width=0.42\linewidth]{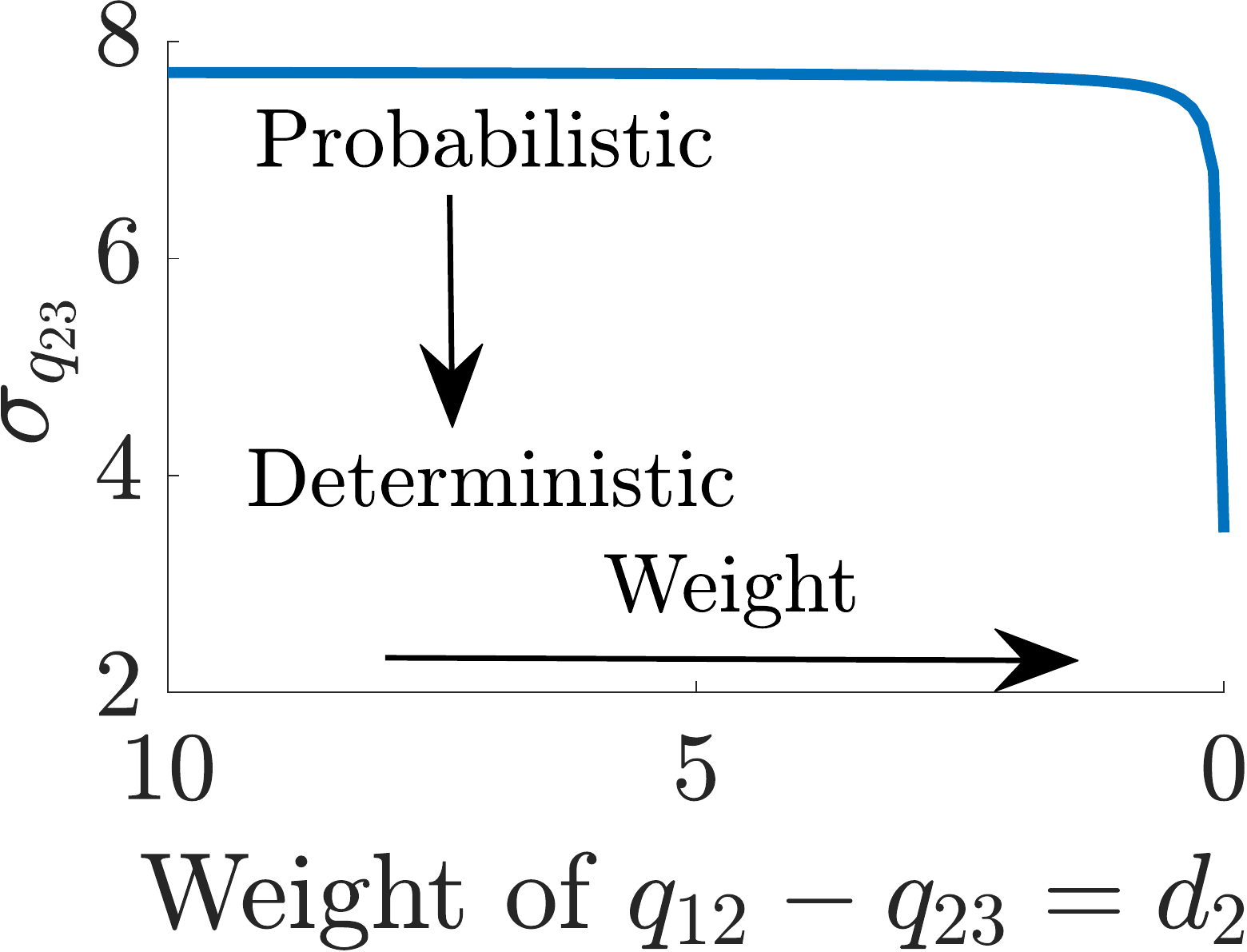}
	\caption{Changes of $\sigma_{q_{23}}$ when weight of $q_{12} - q_{23} = d_2$ decreases.}
	\label{fig:sigma}
	% \vspace{-0.35cm}
\end{figure}
\subsubsection{Over-determined measurements scenario} Under the over-determined scenario, an extra head measurement at Junction 2 (red in Fig.~\ref{fig:setup}) is provided with $y_{h_2} = 910\ \mathrm{ft}$ with $\sigma^2_{{h_2}} = 0.0241$. Its corresponding formulation is also marked as red in~\eqref{eq:3node-LP-matrix}, and results found by~\eqref{equ:cce1} in Corollary~\ref{cor:inver}. 
As we mentioned, a diagonal $6 \times 6$ weight matrix $\m W$ can be  introduced to reflect the importance of each equation of~\eqref{eq:3node-LP-matrix}. We note that  the weight matrix $\m W$ chosen has a significant impact on the state variances as we discussed after presenting Corollary~\ref{cor:inver}. 

\begin{figure}[t]
	%	% \vspace{-0.75cm}
	\centering
	\includegraphics[width=0.7\linewidth]{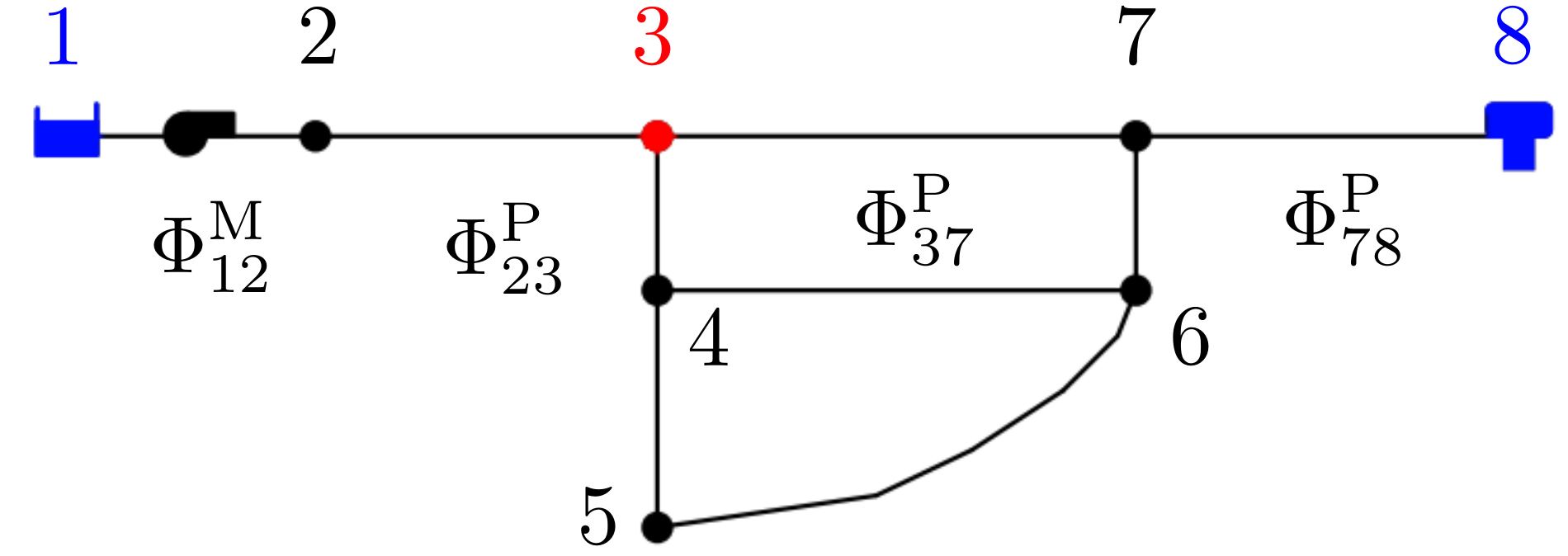}
	\caption{8-node network. The $\m \Phi(\cdot)$ function denotes the nonlinear head loss or head gain.}
	\label{fig:8node}
	% \vspace{-0.25cm}
\end{figure}
Suppose that the weight of the equation $q_{12} - q_{23} = d_2$ is $10$, and we decrease the weight due to lack of confidence in the prediction of demand, then the role it plays is less important while solving for ${\operatorname {K}_{ {\m x}  {\m x} }}$. Fig.~\ref{fig:sigma} demonstrates the changes in $\sigma_{q_{23}} = \sqrt{\Var(q_{23})}$ as a function of decreasing weight of the continuity constraint at Junction 2. We note that  $\sigma_{q_{23}}$ changes from $7.7$ to $3.39$ ($\Var(q_{23}) =11.50$ in Table~\ref{tab:3node-result-covariance}, and $\sqrt{11.50} = 3.39$) when weight of $q_{12} - q_{23} = d_2$ decreases, and it means the uncertainty of $q_{23}$ is reduced with extra measurement at Junction 2, and it relatively becomes more deterministic. This result is reasonable because uncertainty from demand is less valued. The covariance matrix under over-determined scenario with weighted equation is presented as the last three columns in Table~\ref{tab:3node-result-covariance}. Note that the result is not compared with MCS since MCS cannot be applied to over-determined scenario. From the result, we can see that uncertainty is reduced not only on the measured head $h_2$ but also on the other variables compared with the sufficient scenario. 

% \vspace{-0.5cm}

\subsection{Eight-node network}
The 8-node network, adapted from EPANET~\cite{rossman2000epanet}, is a looped network shown in~Fig.~\ref{fig:8node}. The measurement noise, demand, and pipe coefficients are assumed to be in the range of $\pm$20\% around the average values with 99\% confidence level, similar to the setup in previous three-node network.

For this case study, we only show the results under sufficient scenario. That is, the head at reservoir $y_{h_1} = 700\ \mathrm{ft}$ and tank $y_{h_8}= 834\ \mathrm{ft}$ are known.  
\begin{table}[t]
	\caption{Results of 8-node network under sufficient scenario ($T = 1$).}
	\centering
	\renewcommand{\arraystretch}{1.12}
	\begin{tabular}{c|c|c|c|c|c}
		\hline
		\multirow{2}{*}{} & \multirow{2}{*}{\textit{ID}} & \textit{MCS} & \textit{PSE} & \multirow{2}{*}{$\m {\mathrm{AE}}$} & \multirow{2}{*}{$\m {\mathrm{RE}} (\%)$} \\ \cline{3-4}
		&  & \multicolumn{2}{c|}{\textit{Standard deviation} $\sigma$} &  &  \\ \hline
		\multirow{6}{*}{\textit{Head}} 
		& J2 & 4.276 & 4.183 & 0.094 & 2.191\% \\ 
		& J3 &  4.281 & 4.189 & 0.092 & 2.151\% \\  
		& J4 & 4.265 & 4.169 & 0.096 & 2.241\% \\  
		& J5 & 4.365 & 4.267 & 0.098 & 2.241\% \\  
		& J6 & 4.26 & 4.161 & 0.099 & 2.325\% \\  
		& J7 & 4.245 & 4.145 & 0.099 & 2.343\% \\  
		& T8 &  0.015 & 0.015 & 0.00& 0.00\% \\ \hline 
		\multirow{9}{*}{\textit{Flow}}  & P23 & 7.625 & 7.434 & 0.191 & 2.506\% \\  
		& P34 & 14.17 & 13.92 & 0.256 & 1.806\% \\  
		& P45& 4.943 & 4.862 & 0.081 & 1.642\% \\  
		& P37 &16.7 & 16.4 & 0.305 & 1.825\% \\  
		& P46 &12.69 & 12.51 & 0.181 & 1.429\% \\  
		& P76 &17.04 & 16.78 & 0.258 & 1.515\% \\  
		& P65 &5.746 & 5.622 & 0.125 & 2.171\% \\  
		& P78 &13.23 & 13.11 & 0.117 & 0.888\% \\  
		& PU12 &7.625 & 7.434 & 0.191 & 2.506\% \\ \hline \hline
	\end{tabular}
	\label{tab:8noderesultcompare}
	% \vspace{-0.5cm}
\end{table}
We have 17 variables in $\m x(k)$ with 17 equations since we are considering the sufficient measurement scenario and $T$ is set as 1. The final solution and comparison with MCS are presented in~Table~\ref{tab:8noderesultcompare}, and the accuracy of proposed method is guaranteed by  the small relative error $\m {\mathrm{RE}}$. 

\subsubsection{Extended period simulation}
\begin{figure}[t]
	\centering
	\includegraphics[width= 0.9\linewidth]{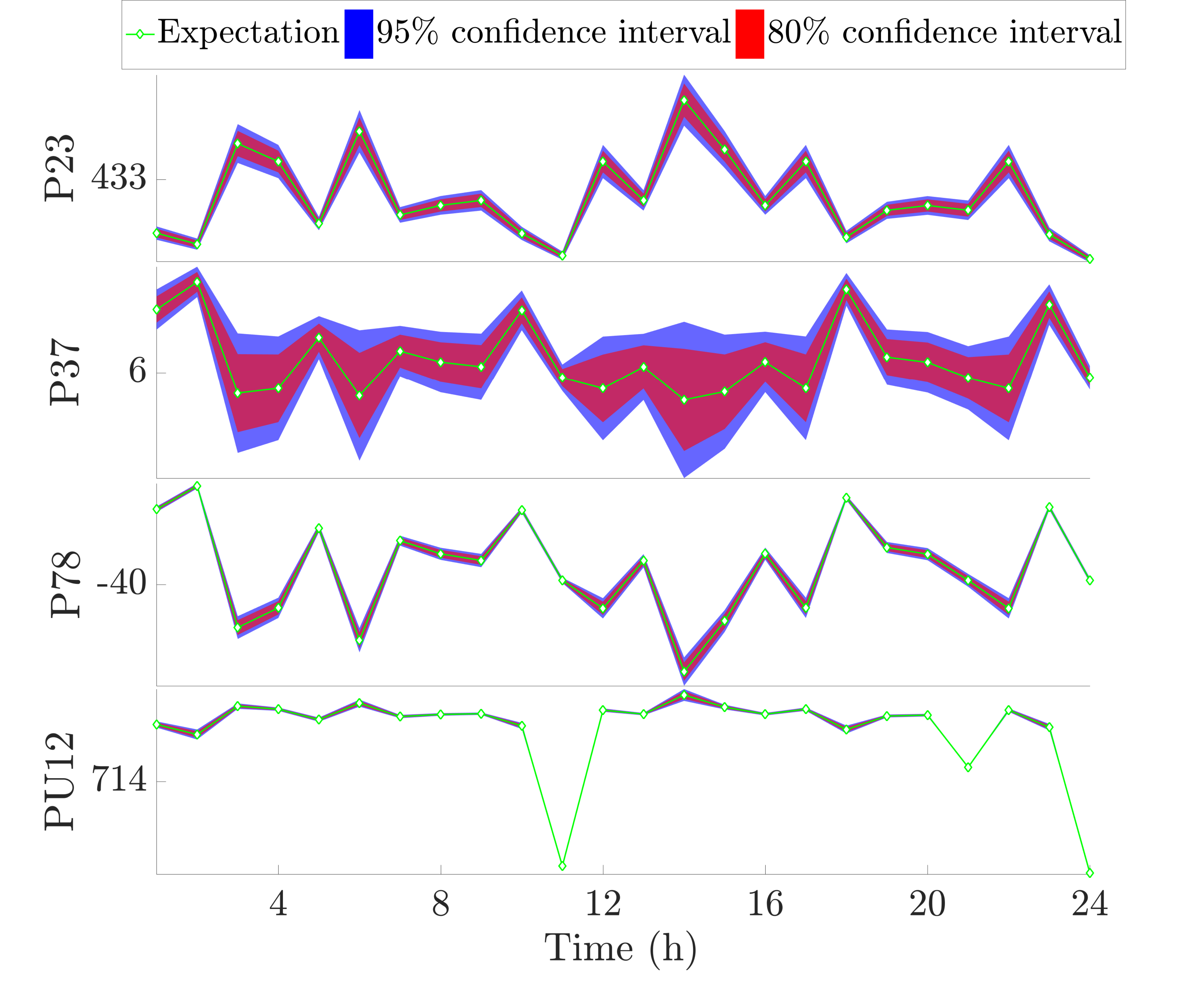}
	\caption{Confidence intervals for flow ($\mathrm{GPM}$)  of Pipes 23, 37, and 78 and Pump 12.}
	\label{fig:flowinterval}
	% \vspace{-10pt}
\end{figure}

\begin{figure}[t]
	\centering
	\includegraphics[width= 0.9\linewidth]{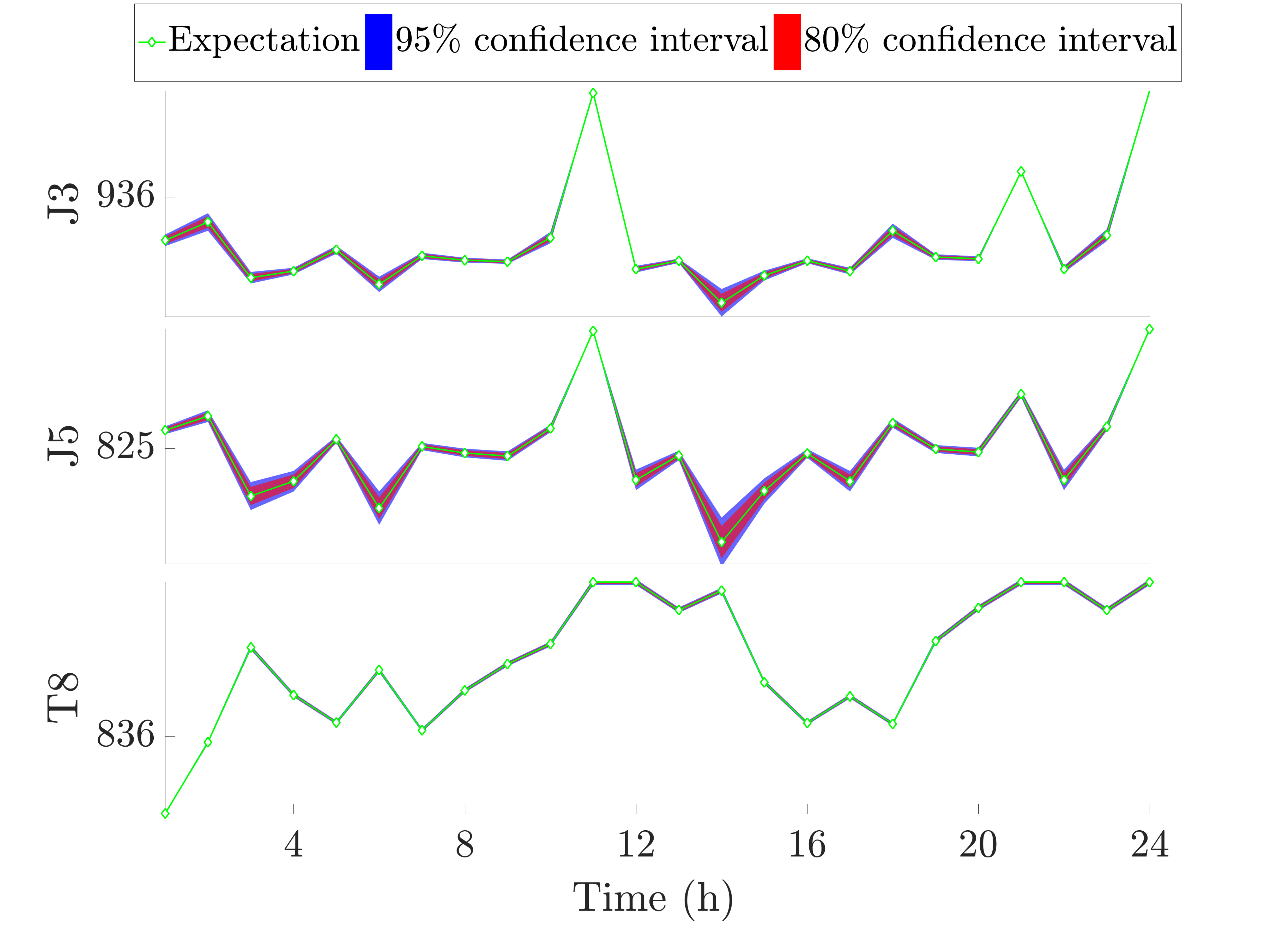}
	\caption{Confidence intervals for head ($\mathrm{ft}$) of Junctions 3 and 5 and Tank 8 under sufficient scenario.}
	\label{fig:headinterval}
	% \vspace{-0.5cm}
\end{figure}

We present the results of an extended period simulation (EPS) for $T = 24$ hours after applying Algorithm~\ref{alg:alg2}. We select three nodes (Junctions 3 and 5, Tank 8) and fours links (Pipes 23, 37, and 78 and Pump 12) (see Fig.~\ref{fig:8node}), and show the $80\%$ and $95\%$ confidence intervals after solving the PSE for 24 hours. The confidence intervals for the estimated heads and flows for the selected nodes and links are presented in Figs.~\ref{fig:flowinterval} and~\ref{fig:headinterval}. 

The green lines in Figs.~\ref{fig:flowinterval} and~\ref{fig:headinterval}  represent the operating point or expectation, and the red interval is in $80\%$ confidence level, whereas the blue one is in $95\%$ confidence level. Both intervals are calculated according to the standard deviation (square root of variance) obtained by PSE in Theorem~\ref{thm:inver}. In  Fig.~\ref{fig:flowinterval}, the expectation of flow in Pipe 37 $\mu_{q_{37}}$ is relatively small, whereas the corresponding variance $\sigma_{q_{37}}$ is larger than the other three links. This is because the flow direction in Pipe 37 changes frequently along with water consumption at junctions. Compared with Fig.~\ref{fig:headinterval},  the flow fluctuates more than the head. The head intervals of $80\%$ and $95\%$ confidence level at Tank 8 are small and almost overlapping with each other due to sensors accuracy and the small variance of measurement noise. The head at Junctions 3 and 7 are mainly decided by the overall network modeling: they tend to fluctuate more compared to Tank 8, hence reflecting the uncertainty in demand and pipe coefficients as well. 


The pumped water is either consumed by junctions or injected into tanks. In particular, the water pumped by Pump 12 is the sum of user demand and water injected into Tank 8, that is $q_{12} = \sum d_{i} + q_{78}$. While it might be intuitive to consider that $\sigma_{q_{12}}$ should be large because all uncertainty of $\sigma_{d_{i}}$ are accumulated in $q_{12}$, this case study shows that  $\sigma_{q_{12}}$ is much smaller than we expected (see flow through Pump 12 in  Fig.~\ref{fig:flowinterval}). This is because Tank 8 acts as a buffer for the network, thereby providing sufficient pressure and flow if demand in network is huge. Otherwise, Tank 8 performs as a junction to consume water when demand---and its corresponding uncertainty---is small. The uncertainty from demand is actually handled by the buffer mechanism from tanks which leads to small flow fluctuation in pumps. 
\begin{figure}[t]
	\centering
	\includegraphics[width=1.2\linewidth]{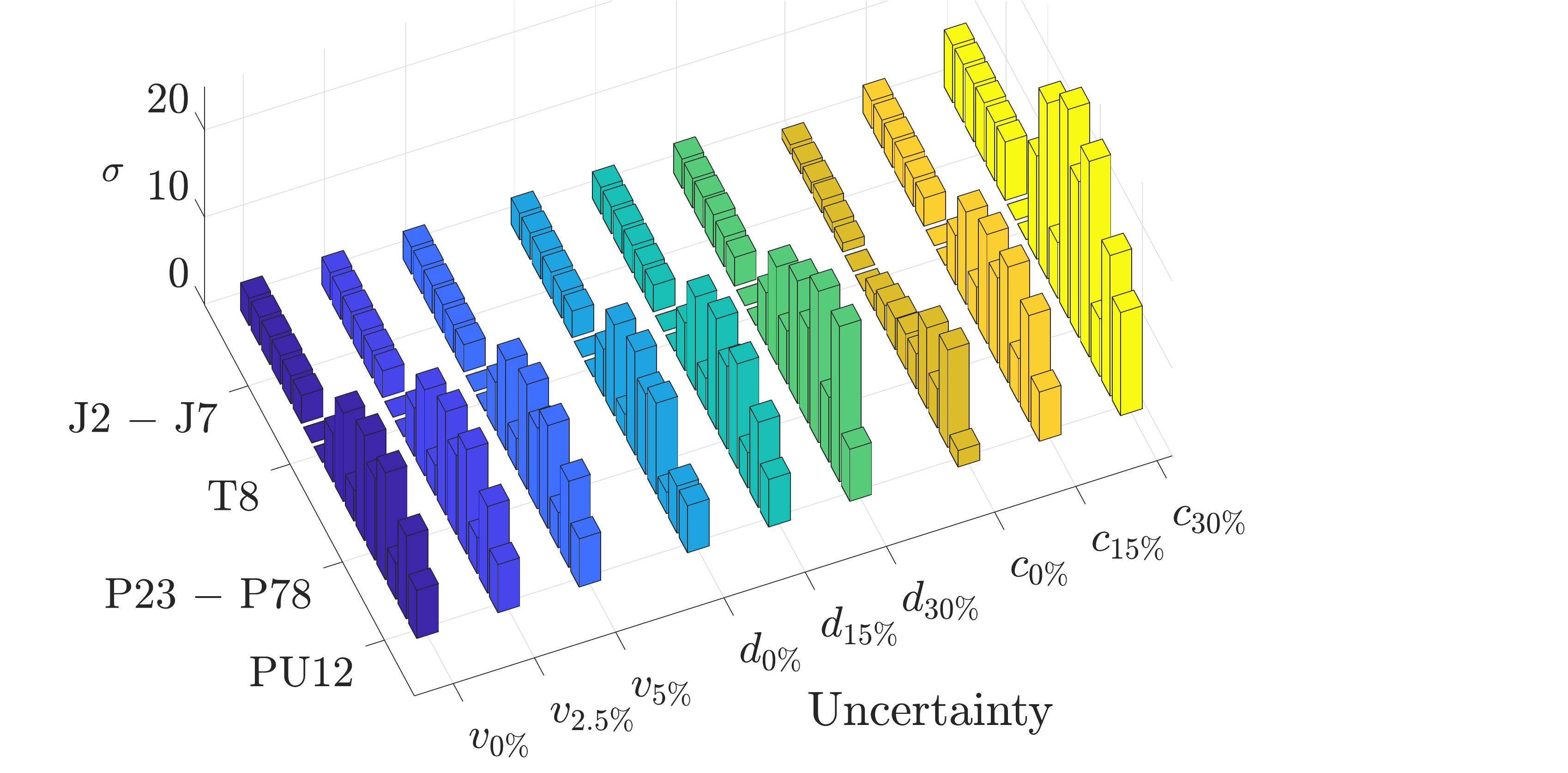}
	\caption{Effects of 3 different uncertainty sources on the standard deviations of components' heads (J2-J7, and T8) and flows (P23-P78, and PU12), see Table~\ref{tab:8noderesultcompare} for details of ID.}
	\label{fig:effect}
	% \vspace{-2em}
\end{figure}
\subsubsection{Effects of different uncertainty sources}\label{sec:effects}
We test the individual effects of the uncertainty sources on the standard deviation of $\m x$ (head and flow) when head at Reservoir 1 $h_1$, Tank 8 $h_8$ are measured and demands at junctions are available.  Fig.~\ref{fig:effect} shows the standard deviations when the margin of errors of uncertainty sources changes. Note that the measurement noise uncertainty $\m v$ is the smallest among all three uncertainty source because the modern sensors are accurate. Most sensors are in $v_{2\%}$\footnote{Quantity $v_{a\%}$ defined in Section~\ref{sec:3-node} means the 99\% of estimated error fall into the range of $\pm$a\% around estimated value. Particularly, $v_{0\%}$ indicates an ideal situation, that is, there is no measurement uncertainty.}, and we assume the worst case is $v_{5\%}$. The pipe roughness coefficient uncertainty  $\m c$ is the largest  among all sources since  $\m c$ is difficult to estimate accurately nowadays, but $\m c$ has a certain range in practice (see Section~\ref{sec:uncertainty}) which can be fully covered by $c_{30\%}$ under 99\% confidence interval, and the worst case for demand are assumed as $d_{30\%}$.

The first three  bar graphs (marked with different shades of blue) show the impact of measurement noise uncertainty $\m v$ (for $h_1$ and $h_8$) on the standard deviation $\sigma$ of  each component with fixed  $d_{15\%}$ and $c_{15\%}$. As the margin of error of $v$ increases from $0\%$, $2.5\%$, to $5\%$, the standard deviation $\sigma$ of some components (P34, P46, P76, and P78) has relatively larger changes compared with other components, but the overall changes for all components are not significant. This indicates that the impact of noise uncertainty is relatively small (see Remark~\ref{remark:impact}).

The second group of three bar graphs present the impact of demand uncertainty under fixed $c_{15\%}$ and $v_{1\%}$. Only the standard deviation of pipes ($\mathrm{P}{23}$--$\mathrm{P}{78}$) changes obviously as the margin of error of demand $d$ goes from $0\%$, $15\%$, up to $30\%$. However, the head at each node still remains relatively unchanged, and this result reflects the demand uncertainty has direct impact on flows instead of heads shown in Fig.~\ref{fig:poof}.

The third group three bar graphs show the impact of pipe roughness uncertainty under fixed $d_{15\%}$ and $v_{2.5\%}$. Note that all standard deviations $\sigma$ of all components change significantly as range of the roughness parameter $c$ goes from $0\%$, $15\%$, to $30\%$. We also test other networks, including the BAK and PESCARA, and the results are similar which verifies the Remark~\ref{remark:impact}.  

The result also indicates the main impact of uncertainty is from pipe roughness coefficients and demand. This guides the network operator in selecting a larger weight $\m W$ in the modeling of PSE and improve the accuracy of state estimation. Specifically, the weight on measurements should be relatively large as the measurements are reliable, whereas the weight on mass balance equation where demands are encoded should be relatively small. The uncertainty from pipe roughness coefficients  have  the greatest impact on state estimation in WDN. An ideal way to improve the performance of state estimation is to curb the pipe roughness coefficients uncertainty propagation through adding more sensors.

Based on the results of the three- and eight-node networks, the Theorem~\ref{thm:looped} has been verified using Kolmogorov-Smirnov test~\cite{daniel1990kolmogorov}, and the plots comparing the CDF based on the results with the standard CDF of normal distribution are omitted due to space limitation.

\subsection{Extended Eight-node network with valves} 

%	% \vspace{-0.25cm}

\begin{figure}[t]
	% \vspace{-0.75cm}
	\centering
	\includegraphics[width=0.75\linewidth]{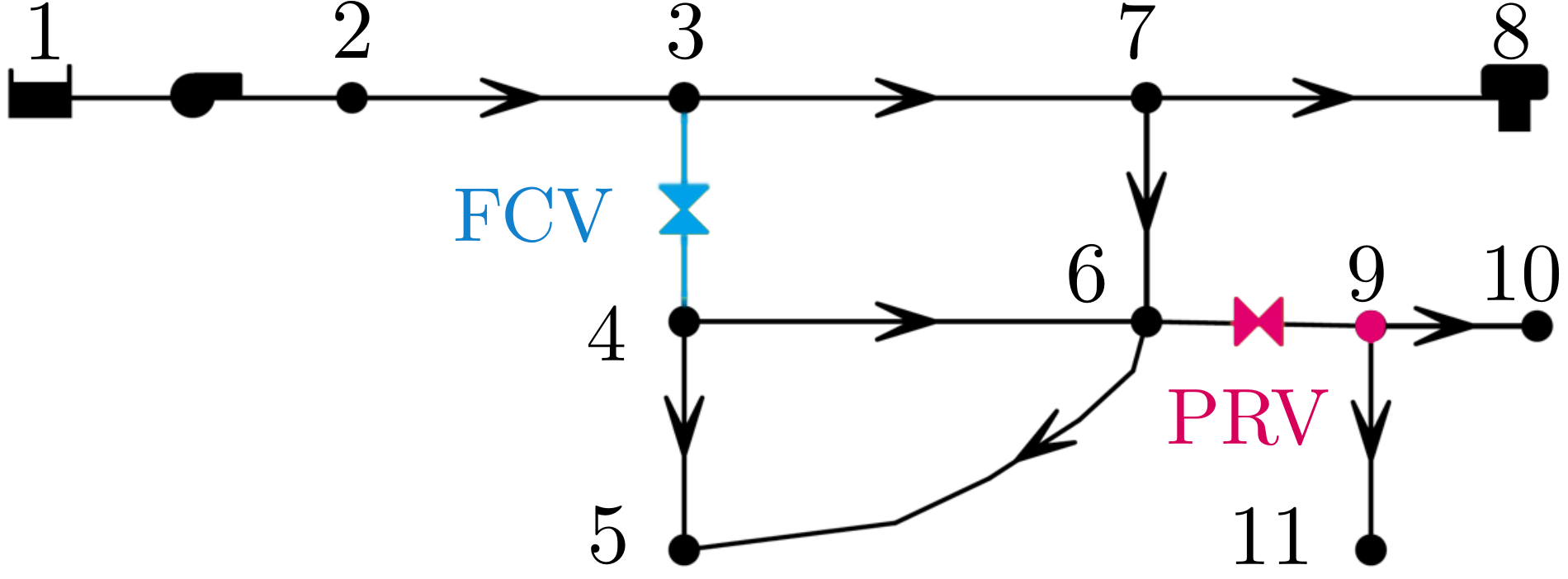}
	\caption{Exteneded 8-node network with FCV and PRV.}
	\label{fig:8nodefcvprv}
	% \vspace{-0.5cm}
\end{figure}

\begin{figure}[t]
	%	% \vspace{-0.75cm}
	\centering
	\includegraphics[width=0.9\linewidth]{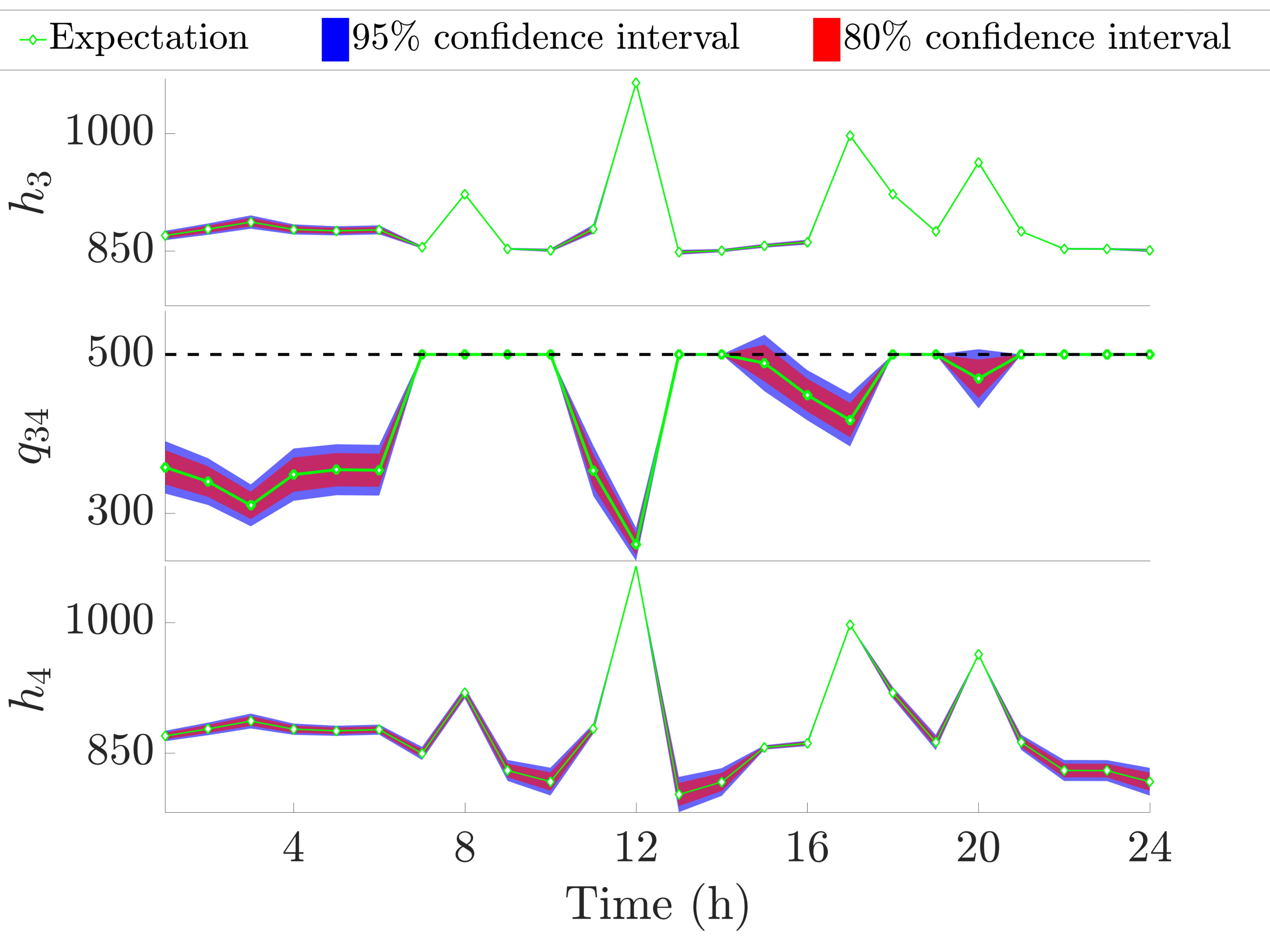}
	\caption{Head and flow related to FCV.}
	\label{fig:FCV}
	% \vspace{-0.3cm}
\end{figure}

\begin{figure}[t]
	\centering
	\includegraphics[width=0.9\linewidth]{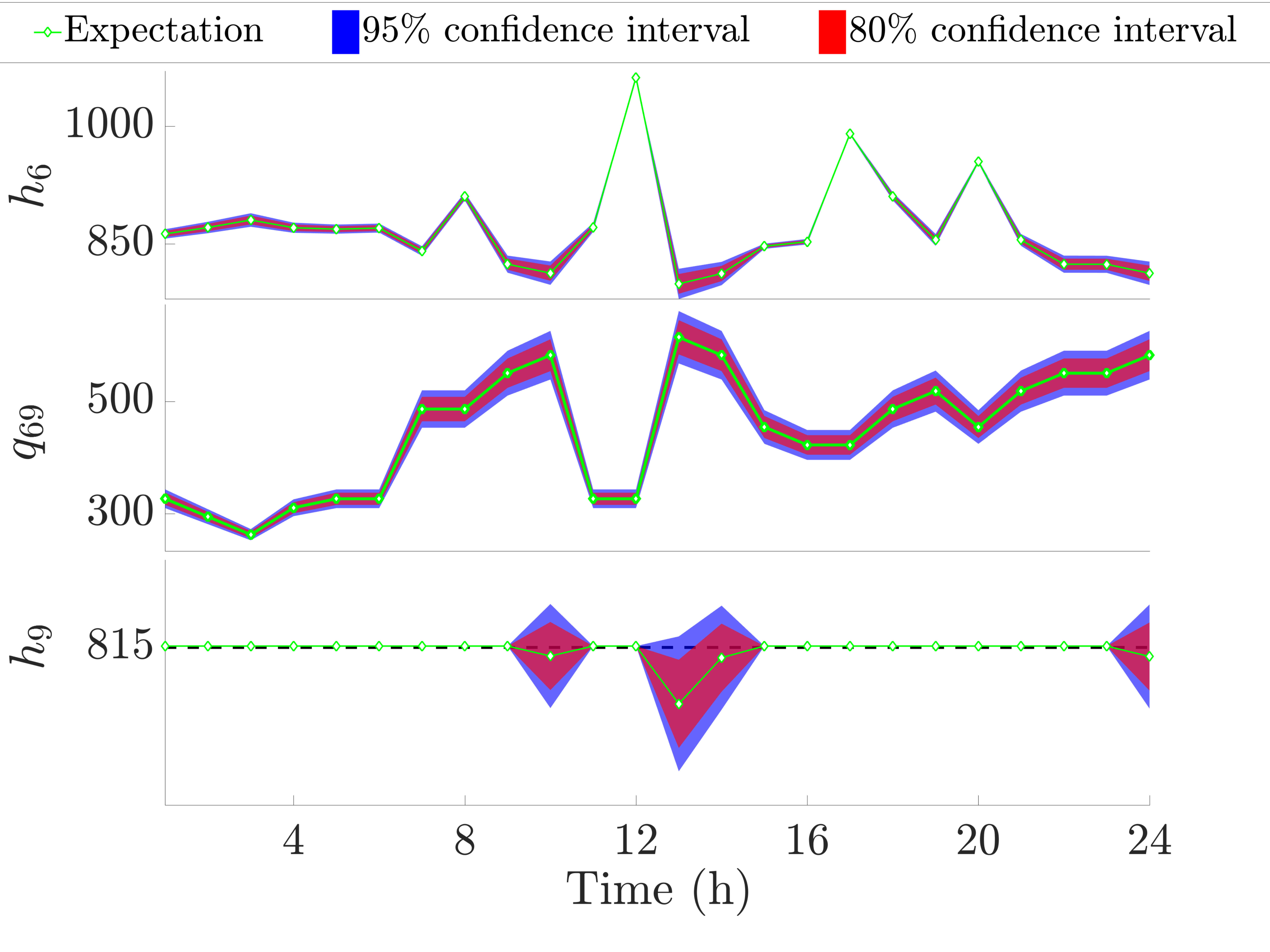}
	\caption{Head and flow related to PRV.}
	\label{fig:PRV}
	%% \vspace{-0.5cm}
\end{figure}

In order to validate our method with water networks with valves, a flow control valve (FCV) and a pressure reducing valve (PRV) are modeled in extended eight-node network shown in Fig.~\ref{fig:8nodefcvprv}. We assume that the FCV is installed between J3 and J4 to limit the flow through the link, and the PRV is installed to maintain the pressure at J9, J10, and J11. 
For the FCV, no head loss exists in $\mathrm{OPEN}$ status and its flow is set to $q^{\mathrm{L}_\mathrm{set}}$ in $\mathrm{ACTIVE}$ status  according to~\eqref{equ:head-fcv-valve}. In particular, the $q^{\mathrm{L}_\mathrm{set}} = 500\ \mathrm{GPM}$ for this test case. The extended period simulation (EPS) results for the FCV are presented in Fig.~\ref{fig:FCV}. {The variables instead of component ID are used as Y labels for convenience in Fig.~\ref{fig:FCV} and Fig.~\ref{fig:PRV}.} We can see that the valve is in $\mathrm{OPEN}$ status during time period $1-6$, $11-12$, and $15-17$, the expectation and variance of $h_3$ and $h_4$ are exactly the same, and the flow through it $q_{34}$ is not controlled by FCV. For the rest of time period in 24 hours, it is in $\mathrm{ACTIVE}$ status, and the expectation of $q_{34}$ is under  $q^{\mathrm{L}_\mathrm{set}}$, and the variance of  $q_{34}$ is zero.

For the PRV, its head $h_9$ is set to $h^{\mathrm{L}_\mathrm{set}}=815\ \mathrm{ft}$ in $\mathrm{ACTIVE}$ status according to~\eqref{equ:head-prv-valve}, and corresponding EPS results for the PRV are presented in Fig.~\ref{fig:PRV}. We can see that it is in $\mathrm{OPEN}$ status during time period $10$, $13-14$, and $24$, the expectation and variance of $h_9$ varies since it is not under control. For the rest of time period in 24 hours, it is in $\mathrm{ACTIVE}$ status, and the expectation of $h_9$ is equal to  $h^{\mathrm{L}_\mathrm{set}}$ with zero variance. Note that once $h_9$ is under control and  equal to less than the settings  $h^{\mathrm{L}_\mathrm{set}}$, the pressures at J10 and J11 are also less than settings. We note that the upper limit of confidence interval may be incorrect when the FCV or PRV are in $\mathrm{OPEN}$ status. The reason is that a valve is considered as a pipe with zero head loss, and the result is not limited by the setting value.
\subsection{Testing scalability and various uncertainty distributions }
In order to verify the scalability of our proposed PSE approach, we {test} different networks {with varying size and complexity} including Anytown, BAK, PESCARA, OBCL, and D-Town networks~\cite{ostfeld2011battle}. The details of each network (number of components), problem size, average relative error of standard deviation (absolute error of expected value), and computational time of each method and network are all summarized in~Table~\ref{tab:computational-time}.  For example,  the Anytown network has $19$ junctions, $3$ reservoirs, $40$ pipes, and one pump. The total number of components is $63$.  The results showcase that the PSE formulation is indeed scalable for networks with hundreds of components due to the highly sparse system matrix $\m A[k]$. We have also verified that $\m A[k]$ for each network is full column rank using \textit{sprank} command from Matlab---thereby ensuring the existence of solution. 

\begin{table}[t]
	% \vspace{-0.25cm}
	\centering
	\caption{Tested networks and their corresponding average $\m {\mathrm{RE}}$ of $\sigma$ ($\m {\mathrm{AE}}$ of $\mu$) and computational time with pipe roughness and measurement uncertainty following normal distributions while demand uncertainty following normal, uniform, or laplace distribution.}
	\label{tab:computational-time}
	\setlength\arraycolsep{0pt}
	\setlength{\tabcolsep}{1.5pt}
	\renewcommand{\arraystretch}{1.9}
	%\resizebox{\linewidth}{!}{%
	\begin{tabular}{c|c|c|c|c|c|c}
		\hline
		\multicolumn{2}{c|}{\textit{Network}} &  \textit{Anytown} & \textit{BAK}& \textit{PESCARA}  & \textit{OBCL}  & \textit{D-Town} \\ \hline
		\multicolumn{2}{c|}{\textit{\makecell{Number of each\\ component$^*$}}} &   \makecell{\{19,3,0,\\40,1,0\}} & \makecell{\{35,1,0,\\58,0,0\}}   & \makecell{\{68,3,0,\\99,0,0\}}  &  \makecell{\{262,1,0,\\288,1,0\}}  &  \makecell{\{364,1,7,\\405,11,4\}}         \\ \hline
		\multicolumn{2}{c|}{\textit{\makecell{ Dimension of $\operatorname {K}_{\m x \m x } $}}} & $63\times63 $  &   $94 \times 94$  & $170 \times 170$  &  $552 \times 552$  &  $783 \times 783$   \\ \hline
		\multirow{3}{*}{\textit{\makecell{Average  \\$\m {\mathrm{RE}}$ \\ of $\sigma$ \%  \\ (Average \\ $\m {\mathrm{AE}}$ \\ of $\mu$ )}}} 
		& \textit{Normal } & \makecell{2.96\%\\ (0.5944)}  & \makecell{1.78\% \\ (1.070)}   & \makecell{1.72\% \\ (0.3249)} & \makecell{0.37\% \\ (0.0072)} & \makecell{1.61\% \\ (0.0547)}    \\  \cline{2-7} 
		&\textit{Uniform} & \makecell{2.64\% \\ (1.152)}  & \makecell{2.39\%\\ (1.374)}  & \makecell{1.96\% \\ (0.5337)} & \makecell{0.29\% \\ (0.0053)}& \makecell{1.51\%\\(0.0643)}   \\ \cline{2-7} 
		&\textit{Laplace} & \makecell{1.36\%\\ (0.1992)}  & \makecell{1.49\%\\ (0.2219)}  & \makecell{1.69\%\\(0.0729)} & \makecell{0.49\%\\(0.0068)} & \makecell{2.05\% \\ (0.1031)} \\ 
		\hline 
		\multirow{2}{*}{\textit{Time (sec)}$^\dagger$} & \textit{MCS }& 1.94 & 2.75  & 4.85 & 16.17 & 21.45 \\ 
		&\textit{PSE} & 0.02  & 0.041  & 0.063 & 0.11 & 0.56  \\ 
		\hline 
		\hline
		\multicolumn{7}{l}{\footnotesize{ \makecell{$^*$\{\# Junctions, \# Reservoirs, \# Tanks, \# Pipes, \# Pumps, \# Valves\}}
		}}\\
		\multicolumn{7}{l}{\footnotesize{ \makecell{$^\dagger$Time under all uncertainty sources follow normal distribution.}
		}}
	\end{tabular}%
	%}
	% \vspace{-0.5cm}
\end{table}
\normalcolor
Besides that, to verify the generality of our method (uncertainty sources can follow different types of distributions, or even do not have to follow specific distributions as long as their expectation and covariance exist), we use normal distributions for pipe roughness and measurement uncertainty, while the demand uncertainty varies from the normal distribution, uniform distribution, to Laplace distribution. The small average relative error of all states is presented as the fourth row in Table~\ref{tab:computational-time}. The results indicates that our method still maintains good performance even with different uncertainty distributions.

The last row list the average computational time for MCS and PSE~\eqref{equ:weighted-hy-matrix} when all uncertainty sources follow the normal distribution. The results of computational time for the mixture of distributions are similar and thus omitted but can be recovered via running the provided codes on Github.

\section{Paper Summary and Future Work}~\label{sec:Future}
The paper's {{objective}} is to develop probabilistic modeling for water network state estimation through scalable computational algorithms considering various sources of uncertainty. To this end, a probabilistic state estimation  algorithm is proposed. We analytically show that the covariance matrix of unknown system states can be linearly expressed by the covariance matrix of uncertainty from measurement noise, network parameters, and demand at arbitrary operating points. Case studies demonstrate the applicability of the proposed method in bounding unmeasured WDN states. Future work will focus on the following relevant problems that are not addressed in this work: \textit{(i)} probabilistic modeling and estimation of network parameters such as demands and pipe roughness coefficients; \textit{(ii)} designing sensor placement methods to maximize the observability of the nonlinear WDN model. 

%	% \vspace{-0.5cm}

\bibliographystyle{IEEEtran}
\bibliography{IEEEabrv,bibfile4}

%	% \vspace{-0.5cm}
\appendices
\section{Equivalent covariance operation }\label{app:covform}
For random vectors $\m a$ and $\m b$,
let $\m a =\m A \m x$, and assume that $\m a= \m b $. We show that $$\Cov((\m a - \m b),(\m a - \m b)) = \m 0 \Leftrightarrow\Cov(\m a, \m a) =  \Cov(\m b, \m b).$$
According to linear combinations covariance property,
\begin{align*}
	\Cov((\m a - \m b),(\m a - \m b)) &= \Cov(\m a, \m a) - 2 \Cov(\m a, \m b) + \Cov(\m b, \m b) \\
	& = \Cov(\m a, \m a) - 2 \Cov(\m b, \m b) + \Cov(\m b, \m b) \\
	&= \Cov(\m a, \m a) -  \Cov(\m b, \m b).
\end{align*}
Note that the fact $\m a = \m b$ is used from the first equality to the second one in above proof. 
Thus,    $$\Cov((\m a - \m b),(\m a - \m b)) = \m 0 \Leftrightarrow \Cov(\m a, \m a) = \Cov(\m b, \m b).$$
% \vspace{-1cm}
\section{Modeling WDN \& Background}~\label{sec:Modeling}
We model  the WDN  by a directed graph $\mathcal{G} = (\mathcal{V},\mathcal{E})$.  Set $\mathcal{V}$ defines the nodes and is partitioned as $\mathcal{V} = \mathcal{J} \cup \mathcal{T} \cup \mathcal{R}$ where $\mathcal{J}$, $\mathcal{T}$, and $\mathcal{R}$ stand for the collection of junctions, tanks, and reservoirs, respectively. Let $\mathcal{E} \subseteq \mathcal{V} \times \mathcal{V}$ be the set of links, and define the partition $\mathcal{E} = \mathcal{P} \cup \mathcal{M} \cup \mathcal{L}$, where $\mathcal{P}$, $\mathcal{M}$, and $\mathcal{L}$ stand for the collection of pipes, pumps, and valves, respectively. The directed graph  $\mathcal{G}$ can be expressed by its incidence matrix $\m E_\mathcal{G}$~\eqref{equ:WFP-Incidence-matrix} which stands for the connection relationship between vertices and edges. For the $i^\mathrm{th}$ node, set $\mathcal{N}_i$ collects its neighboring nodes and is partitioned as $\mathcal{N}_i = \mathcal{N}_i^\mathrm{in} \cup \mathcal{N}_i^\mathrm{out}$, where  			$\mathcal{N}_i^\mathrm{in}$ and $\mathcal{N}_i^\mathrm{in}$ are the sets of  inflow and outflow neighbors of the $i^\mathrm{th}$ node. Notice that the assignment of direction to each link (and the resulting inflow/outflow node classification) is arbitrary.  Thus, $\m E_\mathcal{G}$ is comprised of $-1$, $0$, and $1$ representing negative connection, no connection, and positive connection. Besides, $\m E_\mathcal{G}$ has row partitioned form $[{\m E_{\m h}^\mathrm{P}}^\top {\m E_{\m h}^\mathrm{M}}^\top {\m E_{\m h}^\mathrm{L}}^\top]$  and column partitioned form ${[\m E_{\m q}^\mathrm{J}; \m E_{\m q}^\mathrm{R}; \m E_{\m q}^\mathrm{TK}]}$  from  different perspectives. 

% \vspace{-1.3em}
{ 	\setlength\extrarowheight{3pt}
	\begin{align} ~\label{equ:WFP-Incidence-matrix}
	{\large \m E_\mathcal{G}} = \begin{array}{*{4}{cccc}@{}c}
	&\textit{Pipe}&\textit{Pump}&\textit{Valve}\\
	\cline{2-4}
	\textit{Junction}&\multicolumn{1}{|c|}{\m E^{\mathrm{J}}_\mathrm{P}}&\multicolumn{1}{|c|}{{\m E_{\mathrm{M}}^\mathrm{J}}}&\multicolumn{1}{|c|}{{\m E_{\mathrm{L}}^\mathrm{J}}}& \multirow{1}{*}{$\left.\rule[1ex]{0pt}{1ex}\right\}\m E_{\m q}^\mathrm{J}$}\\
	\cline{2-4}
	\textit{Reservoir}&\multicolumn{1}{|c|}{\m E^{\mathrm{R}}_\mathrm{P}}&\multicolumn{1}{|c|}{{\m E_{\mathrm{M}}^\mathrm{R}}}&\multicolumn{1}{|c|}{{\m E_{\mathrm{L}}^\mathrm{R}}}& \multirow{1}{*}{$\left.\rule[1ex]{0pt}{1ex}\right\}\m E_{\m q}^\mathrm{R}$}\\
	\cline{2-4}
	\textit{Tank}&\multicolumn{1}{|c|}{\m E^{\mathrm{TK}}_\mathrm{P}}&\multicolumn{1}{|c|}{{\m E_{\mathrm{M}}^\mathrm{TK}}}&\multicolumn{1}{|c|}{{\m E_{\mathrm{L}}^\mathrm{TK}}}& \multirow{1}{*}{$\left.\rule[1ex]{0pt}{1ex}\right\}\m E_{\m q}^\mathrm{TK}$}\\
	\cline{2-4}
	\noalign{ \vspace{-6pt}}
	\multicolumn{1}{c}{}	& \multicolumn{1}{@{}c@{}}{\underbrace{\hspace*{2\tabcolsep}\hphantom{.......}}_{{\textstyle {\m E_{\m h}^\mathrm{P}}^\top}}}  & \multicolumn{1}{@{}c@{}}{\underbrace{\hspace*{2\tabcolsep}\hphantom{.......}}_{{\textstyle {\m E_{\m h}^\mathrm{M}}^\top}}} & \multicolumn{1}{@{}c@{}}{\underbrace{\hspace*{2\tabcolsep}\hphantom{.......}}_{{\textstyle \ {\m E_{\m h}^\mathrm{L}}^\top}}}
	\end{array}
	\end{align}
}

According to the principles of \textit{conservation of mass} and \textit{energy}, we present the corresponding modeling of each component of a WDN. %WDN. 
Models of network nodes and links  are given below. %junctions, tanks, reservoirs, pipes, pumps, and valves

\subsubsection{Conservation of mass at junctions, tanks, and reservoirs} Junctions are the points where water flow merges or splits. The expression of mass conservation of the $i^\mathrm{th}$ junction at time $k$  can be written as
\begin{align}~\label{equ:nodes}
\sum_{j \in \mathcal{N}_i^\mathrm{in}} q_{ji} (k)- \sum_{j \in \mathcal{N}_i^\mathrm{out}} q_{ij} (k)= d_i(k),
\end{align}
where  $q_{ji}(k),\;j \in \mathcal{N}_i^\mathrm{in} $ is the inflow from the $j^\mathrm{th}$ neighbor, $q_{ij}(k),\;j \in \mathcal{N}_i^\mathrm{out} $ is the outflow to the $j^\mathrm{th}$ neighbor, and $d_i(k)$ is  the demand extracted from node $i$.  {Here, we assume there is are either no leaks or that leaks uncertainty is encoded in the $\m d(k)$ demand uncertainty in~\eqref{equ:nodes}.} 

The water hydraulic dynamics in the $i^\mathrm{th}$ tank can be expressed by a discrete-time difference equation
% \vspace{-0.5em}
\begin{align}~\label{equ:tankhead}
\hspace{-9pt}h_{i}^{\mathrm{TK}}(k\hspace{-2pt}+\hspace{-2pt} 1)\hspace{-2pt} =\hspace{-2pt} h_{i}^{\mathrm{TK}}(k) \hspace{-2pt}+\hspace{-2pt} \frac{\Delta t}{A_i^{\mathrm{TK}}}\hspace{-3pt}\left(\hspace{-1pt}\sum_{j \in \mathcal{N}_i^\mathrm{in}}\hspace{-3pt}q_{ji}(k)\hspace{-2pt}-\hspace{-9pt}\sum_{j \in \mathcal{N}_i^\mathrm{out}} \hspace{-3pt}q_{ij}(k)\hspace{-3pt}\right)\hspace{-3pt},
\end{align}
where  $h_i^{\mathrm{TK}}$, $A_i^{\mathrm{TK}}$ respectively stand for the head, cross-sectional area  of the $i^\mathrm{th}$ tank, and $\Delta t$ is the discretization time.
We also assume that reservoirs have infinite water supply and the head of the $i^\mathrm{th}$ reservoir is fixed~\cite[Chapter 3.1]{rossman2000epanet}.
\subsubsection{Conservation of energy at pipes, pumps, and valves} The major head loss of a pipe is determined by Hazen-Williams,  and can be expressed as~\eqref{equ:head-flow-pipe} in Table~\ref{tab:models}, where resistance coefficient $R_{ij} = 4.727 L^{\mathrm{P}} (C^{\mathrm{HW}})^{-1.852} (D^{\mathrm{P}})^{-4.871}$ is a function of the corresponding roughness coefficient $C^{\mathrm{HW}}$, pipe diameter $D^{\mathrm{P}}$, and pipe length $L^{\mathrm{P}}$. The minor head loss of a pipe is not considered.
\begin{table*}[t]
	% \vspace{-0.5cm}
	\def\arraystretch{1.0}
	\caption{Hydraulic modeling of pipes and pumps (time index $k$ is ignored for each variable for simplicity).}
	\centering
	\makegapedcells
	\setcellgapes{0.2pt}
	\begin{tabular}{c|c|c}
		\hline
		& {\textit{Pipes}}  & \textit{Pumps} \\ \hline
		{\textit{Original Hydraulic Model}}
		&  \parbox{7cm}{
			% \vspace{-0.7em}
			\begin{align}~\label{equ:head-flow-pipe}
			\Delta h_{ij}^\mathrm{P}  = h_{i} - h_{j} = R_{ij} {q_{ij}}|q_{ij}|^{\alpha-1}
			\end{align}
			% \vspace{-1em}
		}
		& 
		\parbox{6cm}{
			% \vspace{-0.7em}
			\begin{align} \label{equ:head-flow-pump}
			\hspace{-10pt} \Delta h_{ij}^\mathrm{\mathrm{M}} = h_{i} - h_{j} = -{s_{ij}^2}(h_0 - r  (q_{ij} s_{ij}^{-1})^\beta )
			\end{align}
			% \vspace{-1em}
		} 
		\\ 
		\hline
		{\textit{First Order Taylor Series Form}} 
		& \parbox{7cm}{
			% \vspace{-0.7em}
			\begin{align} \label{equ:head-flow-pipe-linear}
			\Delta h_{ij}^\mathrm{P}  = h_{i} - h_{j} = k_{q_{ij}}^{\mathrm{P}} q_{ij} + k_{c_{ij}}^{\mathrm{P}}c_{ij}+ {b_{ij}^{\mathrm{P}}}
			\end{align}
			% \vspace{-1em}
		} 
		&\parbox{6cm}{
			% \vspace{-0.7em}
			\begin{align} \label{equ:head-flow-pump-linear}
			\hspace{-10pt}\Delta h_{ij}^\mathrm{\mathrm{M}} = h_{i} - h_{j} = -k_{ij}^{\mathrm{M}} q_{ij}+ {b_{ij}^{\mathrm{M}}}
			\end{align}
			% \vspace{-1em}
		}
		\\ \hline \hline
	\end{tabular}
	\label{tab:models}%
	%		}
	%	\end{subtable}
	% \vspace{-0.5cm}
\end{table*}

The head gain generated by a pump from suction node $i$ to delivery node $j$ is determined by the pump curve and can be expressed as~\eqref{equ:head-flow-pump} in Table~\ref{tab:models}, where $q_{ij}$ and $s_{ij}$ are the flow and  speed of a pump; quantities $h_0$, $r$ and $\beta$ are the pump curve coefficients. Flow control valve (FCV) and pressure reducing valve (PRV) are commonly used valves to regulate flow or pressures, respectively, and are controlled through valve openness or set points. The valve models used in our paper is based on~\cite{rossman2000epanet}.
FCV  limits the flow to a specified setting  $q^{\mathrm{L}_\mathrm{set}}$, when the head $h_{i}$ at upstream node~$i$ is greater than the head $h_{j}$ at downstream node~$j$; otherwise, FCV is treated as an open pipe with minor head loss. In short, FCV can be modeled as
\begin{align}~\label{equ:head-fcv-valve}
\begin{dcases}
h_{ij}^\mathrm{L} = h_{i} - h_{j} = l_{ij} {q_{ij}}|q_{ij}| = 0,\mathrm{OPEN}\\
q_{ij} = q^{\mathrm{L}_\mathrm{set}}  ,\;\text{ACTIVE}
\end{dcases}
\end{align} 
where  $l_{ij} = \frac{k_{ij}}{2g A^2_{ij}}$, and $k_{ij}$ is the minor head loss coefficient, $A_{ij}$ is the corresponding cross-sectional area of the FCV, $g$ is the acceleration of gravity, and all three parameters are constant which make $l_{ij}$ also as a constant. $ q^{\mathrm{L}_\mathrm{set}}$ is the setting value. 

PRV limits the pressure at a specific location (reverse flow is not allowed) and set the pressure to $P^\mathrm{set}$ on its downstream side $j$ when the upstream side $i$ pressure is higher than $P^\mathrm{set}$~\cite[Chapter 3.1]{rossman2000epanet}. Given the status of a PRV, it can be modeled as 
% \vspace{-0.5em}
\begin{align}~\label{equ:head-prv-valve}
\begin{dcases}
h_{ij}^\mathrm{L} = h_{i} - h_{j} = l_{ij} {q_{ij}}|q_{ij}| = 0,\mathrm{OPEN}\\
h_{j}  = h^{\mathrm{L}_\mathrm{set}}  ,\;\text{ACTIVE}
\end{dcases}
\end{align} 
where $l_{ij}$ is the same constant as the one of PRV;  $h^{\mathrm{L}_\mathrm{set}}$ is a constant head converted from the constant pressure setting $P^\mathrm{set}$. The fact that reverse flow is not allowed in PRVs can be expressed as a limit $q_{ij} \geq 0$ and included in~\eqref{equ:inequality}. Besides that, minor head losses in PRVs and FCVs are not considered in this paper, and that $l_{ij}$ is set as $0$ results in the linear model for PRVs and FCVs.

% \vspace{-0.5cm}

\section{Derivation of  Linearized  DAE} \label{sec:LinearDAE}
The linearized pipe and pump modeling~\eqref{equ:linearPipesPumps}  is included in Table V and can be rewritten as (time index $k$ is ignored during the following derivation)
\begin{align}~\label{equ:linearPipesPumps-matrix}\hspace{-10pt}
\setlength\arraycolsep{0pt}
\begin{bmatrix}
\setlength\arraycolsep{0pt}
\m \Delta \m h^\mathrm{P}  \\
\m \Delta \m h^\mathrm{\mathrm{M}} 
\end{bmatrix} \hspace{-3pt}=  \hspace{-3pt} \begin{bmatrix}
\m K_q^{\mathrm{P}}  & \\
& -\m K_q^{\mathrm{M}} 
\end{bmatrix}  \hspace{-4pt} \begin{bmatrix}
\m q^{\mathrm{P}} \\
\m q^{\mathrm{M}} 
\end{bmatrix}   \hspace{-3pt}+ \hspace{-3pt} \begin{bmatrix}
\m K_c^{\mathrm{P}}  & \\ 
& \m 0 
\end{bmatrix}  \begin{bmatrix}
\m c \\
\m 0
\end{bmatrix}  \hspace{-3pt}+ \hspace{-3pt} \begin{bmatrix}
\m b^{\mathrm{P}} \\
\m b^{\mathrm{M}}
\end{bmatrix}.
\end{align}
Head difference $\m \Delta \m h^\mathrm{P} $ and $\m \Delta \m h^\mathrm{\mathrm{M}} $ across any pipe or pump can also be represented as $\m h_i - \m h_j$, and  $\m h_i$ or $\m h_j$ must be in set $\m h^{\mathrm{J}}$, $\m h^{\mathrm{R}}$, and 
$\m h^{\mathrm{TK}}$. Hence, we have 
\begin{align}
\setlength\arraycolsep{0pt}
\begin{bmatrix}
\setlength\arraycolsep{0pt}
\m \Delta \m h^\mathrm{P}  \\
\m \Delta \m h^\mathrm{\mathrm{M}} 
\end{bmatrix} = \m E_{h'} 	\begin{bmatrix}
\m h^{\mathrm{J}}\\
\m h^{\mathrm{R}}\\
\m h^{\mathrm{TK}}
\end{bmatrix} = \m E_{h'} \m h,\notag
\end{align}
where $\m E_{h'}$ is a submatrix of $\m E_{\mathcal{G}}$~\eqref{equ:WFP-Incidence-matrix} defined as $\m E_{h'} = [{\m E_{\m h}^\mathrm{P}}\ {\m E_{\m h}^\mathrm{M}}]^\top$. Substituting it to~\eqref{equ:linearPipesPumps-matrix}, the linearized hydraulic equations are obtained as
\begin{align}~\label{equ:linearPipesPumps-matrix2}
\setlength\arraycolsep{0pt}
\m E_{h'}  \m h\hspace{-2pt}=  \hspace{-2pt}\m E_{k}   \begin{bmatrix}
\m q^{\mathrm{P}} \\
\m q^{\mathrm{M}} 
\end{bmatrix}   \hspace{-3pt}+ \hspace{-3pt} \begin{bmatrix}
\m K_c^{\mathrm{P}}  & \\ 
& \m 0 
\end{bmatrix}  \begin{bmatrix}
\m c \\
\m 0
\end{bmatrix}  \hspace{-3pt}+ \hspace{-3pt} \begin{bmatrix}
\m b^{\mathrm{P}} \\
\m b^{\mathrm{M}}
\end{bmatrix},
\end{align}
where $\m E_{k} = \diag\{\m K_q^{\mathrm{P}},-\m K_q^{\mathrm{M}} \}$. We can reorganize \eqref{equ:linearPipesPumps-matrix2} into the following form
\begin{align} \label{equ:linearizeenergy}
\setlength\arraycolsep{0pt}
\begin{bmatrix}
\m E_{h'} \  -\hspace{-2pt} \m E_k
\end{bmatrix} \begin{bmatrix}
\m h \\
\m q^{\mathrm{P}} \\
\m q^{\mathrm{M}} 
\end{bmatrix}  =    \begin{bmatrix}
\m K_c^{\mathrm{P}}  & \\ 
& \m 0 
\end{bmatrix}  \begin{bmatrix}
\m c \\
\m 0
\end{bmatrix}  \hspace{-3pt}+ \hspace{-3pt} \begin{bmatrix}
\m b^{\mathrm{P}} \\
\m b^{\mathrm{M}}
\end{bmatrix}.
\end{align}
Recall that both FCV~\eqref{equ:head-fcv-valve} and PRV~\eqref{equ:head-prv-valve} have two different statuses, and the total number of combinations is four. However, regardless of whether valves are open or active, it is linear. Hence, we can model them as
\begin{align} \label{equ:valve}
\begingroup
\setlength\arraycolsep{0pt}
\renewcommand{\arraystretch}{1.2}
[\m E_{l} \ | \ \m 0_{(n_p + n_m) \times n_q} \  | \ \m E_{l'}] \begin{bmatrix}
\m h \\ \cline{1-1}
\m q^{\mathrm{P}} \\
\m q^{\mathrm{M}} \\ \cline{1-1}
\m q^{\mathrm{L}} 
\end{bmatrix}   =  \begin{bmatrix}
\m 0\\
\m 0\\
\m l
\end{bmatrix} ,
\endgroup
\end{align}
where matrix $\m E_{l}$ includes rows of matrix ${\m E_{\m h}^\mathrm{L}}^\top$ in $\m E_{\mathcal{G}}$ [see~\eqref{equ:WFP-Incidence-matrix}] for FCVs and PRVs in open status; or selects the appropriate entry of $\m h$ for a PRV in active status; matrix $\m E_{l'} $ is  comprised of 0 or 1 corresponding to the active status in~\eqref{equ:head-fcv-valve}; $\m l$ is a $n_l \times 1$ vector consisting of the combination of either $\m 0$ or setting of FCV or PRV. With the compact form of conversion of mass~\eqref{equ:nodes-abcstracted}, linearized form of conservation of energy~\eqref{equ:linearizeenergy}, and linear valve modeling~\eqref{equ:valve}, the overall  linearized Hydraulic Modeling can be rewritten as
\setlength\arraycolsep{1pt}
\begin{align} ~\label{equ:LinearHyd}
\renewcommand{\arraystretch}{1.15}
\hspace{-1em}\underbrace{ \left[ \begin{array}{c|c|c}
	\m 0 & \multicolumn{2}{c}{\m E_{\m q}^\mathrm{J}}  \\
	\cline{2-3}
	\m E_{h'}  & -\m E_k & \m 0  \\ %_{(n_p+n_m)\times n_l}
	{\m E_l} & \m 0 & \m E_{l'}
	\end{array}\right] }_{\text{ \normalsize $\m E(k)$}}   
&\begingroup
\renewcommand{\arraystretch}{1.15}
\underbrace{\begin{bmatrix}
	\m h \\ \cline{1-1}
	\m q^{\mathrm{P}} \\
	\m q^{\mathrm{M}} \\ \cline{1-1}
	\m q^{\mathrm{L}} 
	\end{bmatrix}}_{\text{ \normalsize $\m x$}} 
\endgroup
=   \begin{bmatrix}
\m 0 & & &\\
& \m K_c^{\mathrm{P}}  & &\\ 
& & \m 0 &\\
& & & \m 0
\end{bmatrix} \hspace{-4pt}  \begin{bmatrix}
\m 0 \\ 
\m c \\
\m 0 \\
\m 0
\end{bmatrix}  \hspace{-3pt}+ \hspace{-3pt} \begin{bmatrix}
\m 0 \\
\m b^{\mathrm{P}} \\
\m b^{\mathrm{M}}\\
\m 0
\end{bmatrix} \notag \\
& +     \begin{bmatrix}
-\m d\\
\m 0\\
\m 0\\
\m 0 
\end{bmatrix}  +      \begin{bmatrix}
\m 0\\
\m 0\\
\m 0\\
\m l
\end{bmatrix} =
\underbrace{ \setlength\arraycolsep{1pt}
	\begin{bmatrix}
	-\m d\\
	\m  K_c^{\mathrm{P}} \m c + \m b^{\mathrm{P}}  \\
	\m b^{\mathrm{M}}\\
	\m l
	\end{bmatrix}.}_{\text{ \normalsize $\m z$}}
\end{align}

Here we redefine the matrices in the left and right hand as $\m E(k)$ and $\m z$ for simplicity, the corresponding dimension are $(n_j + n_{q})\times n_{x}$ and $(n_j + n_{q})\times 1$,  and the final linearized DAE can be rewritten as
\setlength\arraycolsep{3pt}
\begin{subequations} 
	\begin{align*}
	\mathrm{DAE_{linear}}:\;\;\;	\m h^{\mathrm{TK}}(k + 1) &= \m A_{\m h} \m h^{\mathrm{TK}}(k) + \m B_{\m q} \m q(k)   \\
	\hspace{-10pt}	\m 0&=\m E(k) \m x(k) - \m z(k).
	\end{align*} 
\end{subequations}
This shows the compact, linearized DAE model.

\section{Proof of Theorem~\ref{thm:inver} and Corollary~\ref{cor:inver}}\label{app:inver}

\begin{proof}(of Theorem~\ref{thm:inver})
	The proof consists of two parts. The first part shows that $\m A[k]$ is indeed full rank.  Actually if the submatrix $\m A^\mathrm{s}[k]$ is full rank, then $\m A[k]$ is full rank.  This is because $\m A[k]$ has more rows than columns, and it is clear that each row in $\m A^\mathrm{TK}[k]$ is linearly independent from the rows in $\m A^\mathrm{s}[k]$. Next we show the $\m A^\mathrm{s}[k]$ is full rank and it is given by
	\begin{equation*}
	\m A^\mathrm{s}[k]= \diag\{\m A^\mathrm{s}(k),\ldots, \m A^\mathrm{s}(k+T) \}.
	\end{equation*}
	where $\m A^\mathrm{s}(k)$ is given in~\eqref{equ:combin}. After substituting $\m E(k)$~\eqref{equ:LinearHyd} into $\m A^\mathrm{s}(k)$, the $\rank(\m A^\mathrm{s}(k))$ can be expressed as the left hand side in~\eqref{equ:rank}. Matrix row operations do not impact the matrix rank, thus, $\rank(\m A^\mathrm{s}(k)) = \rank(\m{\tilde{A}^\mathrm{s}}(k))$ where $\m{\tilde{A}^\mathrm{s}}(k)$ is the new matrix after applying row operations on $\m A^\mathrm{s}(k)$.
	\setlength\arraycolsep{0pt}
	% \vspace{-0.0em}
	\begin{equation}~\label{equ:rank}
	\renewcommand{\arraystretch}{1.15}
	\hspace{-1em}			{\small \rank \hspace{-2pt} \underbrace{\left( \hspace{-2pt} \left[ \begin{array}{c|cc}
			\m 0 & \multicolumn{2}{c}{\m E_{\m q}^\mathrm{J}(k)}  \\
			\cline{2-3}
			\m E_{h'}(k)  & -\m E_k(k) & \m 0  \\ %_{(n_p+n_m)\times n_l}
			\cline{2-3}
			{\m E_l(k)} & \m 0 & \m E_{l'}(k) \\
			\cline{1-3}
			\multicolumn{3}{c}{\m C}
			\end{array}\right] \hspace{-2pt} \right)}_{\textstyle \m A^\mathrm{s}(k)}\hspace{-3pt}  = \hspace{-2pt} 
		\rank \hspace{-2pt} \underbrace{\left( \hspace{-2pt} \left[ \begin{array}{c|cc}
			\m 0 & \multicolumn{2}{c}{\m E_{\m q}^\mathrm{J}(k)}  \\ \hline
			\multicolumn{3}{c}{\m C}\\ \hline
			\m E_{h'}(k)  & -\m E_k(k) & \m 0  \\ %_{(n_p+n_m)\times n_l}
			\cline{2-3}
			{\m E_l(k)} & \m 0 &  \m E_{l'}(k) \\
			\end{array}\right] \hspace{-2pt} \right)}_{\textstyle \m{\tilde{A}^\mathrm{s}}(k)}}
	\end{equation}
	 We consider the simple case (sufficient scenario) and show $\m{\tilde{A}^\mathrm{s}}(k)$ is full column rank first, that is, $\m C$ only selects heads at tanks and reservoirs. Because the linear models of PRVs and FCVs are used, and after expanding ${\m E_{\m q}^\mathrm{J}}$, $\m E_{h'}$, $\m E_{k}$, $\m E_{l}$ and $\m E_{l'}$ which are defined in Appendix~\ref{sec:LinearDAE}, the $\m{\tilde{A}^\mathrm{s}}(k)$ can be expressed as
	\begin{equation}~\label{equ:LinearMatrix}
	%\resizebox{.7\hsize}{!}{\input{Amatrix1.tikz}}
	%\includegraphics[width=0.8\linewidth,valign=c]{Amatrix1.pdf},
	\begin{tikzpicture}[
	baseline={([yshift=-.5ex]current bounding box.center)},
	style0/.style={
		matrix of math nodes,
		row sep=-\pgflinewidth,
		column sep=-\pgflinewidth,
		every node/.append style={text width=#1,anchor=center,minimum height=3.5ex},
		anchor=base,
		nodes in empty cells,
		left delimiter=.,
		right delimiter=.
	},
	style1/.style={
		matrix of math nodes,
		row sep=-\pgflinewidth,
		column sep=-\pgflinewidth,
		every node/.append style={text width=#1,anchor=center,minimum height=3.5ex},
		anchor=base,
		nodes in empty cells,
		left delimiter={[},
		right delimiter={]}
	},
	style2/.style={
		matrix of math nodes,
		row sep=-\pgflinewidth,
		column sep=-\pgflinewidth,
		every node/.append style={text width=#1,anchor=center,minimum height=3.5ex},
		anchor=base,
		nodes in empty cells,
		left delimiter={[},
		right delimiter={]}
	}
	]
	%\matrix[style0=1.5cm] (0mat)
	%{  \m{\tilde{A}} = \\
	%};
	\matrix[style1=0.9cm] (1mat)
	{
		\m 0_{n_j \times n_j} & & & {\m E^{\mathrm{J}}_\mathrm{P}} & {{\m E_{\mathrm{M}}^\mathrm{J}}}  & {{\m E_{\mathrm{L}}^\mathrm{J}}} \\
		& \m I_{n_r \times n_r}& & & \\
		& & {\m I_{n_t \times n_t}}& & & \\
		{\m E_{\mathrm{J}}^\mathrm{P}}&{\m E_{\mathrm{R}}^\mathrm{P}}  & {\m E_{\mathrm{TK}}^\mathrm{P}} &  -\m K_q^{\mathrm{P}} & &\\
		{{\m E^{\mathrm{M}}_\mathrm{J}}} &{{\m E^{\mathrm{M}}_\mathrm{R}}} & {{\m E^{\mathrm{M}}_\mathrm{TK}}}  & & \m K_q^{\mathrm{M}}&\\
		{{\m E^{\mathrm{L}}_\mathrm{J}}} &{{\m E^{\mathrm{L}}_\mathrm{R}}} &{{\m E^{\mathrm{L}}_\mathrm{TK}}}& & & {\m E_{l'}} \\
	};
	\draw [opacity=0.4] (1mat-1-1.south west) -- (1mat-1-6.south east);
	\draw [opacity=0.4] (1mat-3-1.south west) -- (1mat-3-6.south east);
	%\draw [dashed] (1mat-5-1.south west) -- (1mat-5-6.south east);
	\draw [opacity=0.4]  (1mat-5-4.south west) -- (1mat-5-6.south east);
	\draw  [opacity=0.4]  (1mat-1-3.north east) -- (1mat-1-3.south east);
	\draw  [opacity=0.4]  (1mat-4-3.north east) -- (1mat-6-3.south east);
	\end{tikzpicture}.
	\end{equation}
Note that \textit{(i)} $\m{\tilde{A}^\mathrm{s}}$ is a  square matrix with size of $n_x \times n_x$; \textit{(ii)} $\m{\tilde{A}^\mathrm{s}}$ is split into six block-columns as shown in~\eqref{equ:LinearMatrix}, that is  $\m{\tilde{A}^\mathrm{s}} = [\m{\tilde{A}^\mathrm{s}}_{1} \  \m{\tilde{A}^\mathrm{s}}_{2} \  \m{\tilde{A}^\mathrm{s}}_{3} \  \m{\tilde{A}^\mathrm{s}}_{4} \  \m{\tilde{A}^\mathrm{s}}_{5} \  \m{\tilde{A}^\mathrm{s}}_{6}]$;  \textit{(iii)} $\m K_q^{\mathrm{P}}$ and $\m K_q^{\mathrm{M}}$ are diagonal matrices; \textit{(iv)} $\m E_{\m q}^\mathrm{J} = [{\m E^{\mathrm{J}}_\mathrm{P}} \ {{\m E_{\mathrm{M}}^\mathrm{J}}}  \ {\m E_{\mathrm{L}}^\mathrm{J}} ] = \{\m E_{\mathrm{J}}^\mathrm{P}, {{\m E^{\mathrm{M}}_\mathrm{J}}}, \m E^{\mathrm{L}}_\mathrm{J}\}^\top = {\m E_\mathrm{J}^{\m q}}^\top$, and  $\m E_{\m q}^\mathrm{J}$ is the reduced incidence matrix of $\m E_\mathcal{G}$ because reservoirs and tanks are not included in $\m E_{\m q}^\mathrm{J}$; the last diagonal block in  $\m{\tilde{A}^\mathrm{s}}$ is $\m E_{l'}$  and it covers all possible situations when PRVs and/or FCVs are open or active. 

In order to prove $\m{\tilde{A}^\mathrm{s}}$ is full column rank, we only need to prove each column is linearly independent from each other, meaning that each block $\m{\tilde{A}^\mathrm{s}}_i$ is full column rank and it is linearly independent from $\m{\tilde{A}^\mathrm{s}}_j$ when $i \neq j \in \{1,\ldots,6\}$. It is clear that the four columns ($ \m{\tilde{A}^\mathrm{s}}_{2}$, $\m{\tilde{A}^\mathrm{s}}_{3}$, $\m{\tilde{A}^\mathrm{s}}_{4}$, and $\m{\tilde{A}^\mathrm{s}}_{5}$) are linearly independent of each other due to the fact that the $\m I_{n_r \times n_r}$, $\m I_{n_t \times n_t}$, $\m K_q^{\mathrm{P}}$, and $\m K_q^{\mathrm{M}}$ are diagonal matrices. Next we need to show $ \m{\tilde{A}^\mathrm{s}}_{1}$ ($\m E_{\m q}^\mathrm{J}$) and $ \m{\tilde{A}^\mathrm{s}}_{6}$ ($\m E_{\mathrm{L}}^\mathrm{J}$) are full column rank and linearly independent from the rest of columns.

We note that $\m E_{\m q}^\mathrm{J}$ is encoded in~\eqref{equ:nodes}. Examining the water flow in the mass balance equation~\eqref{equ:nodes} and using network- and electric circuit-theoretic results from~\cite[Theorem 3.2]{muthuswamy2018introduction},  matrix~$\m E_{\m q}^\mathrm{J}$ is linearly row independent. In other words, the $\m E_{\m q}^\mathrm{J}$ ($\m E^{\m q}_\mathrm{J}$) is linearly row (column) independent. Hence, $ \m{\tilde{A}^\mathrm{s}}_{1}$ is full column rank. We next discuss the full column rank property of $\m{\tilde{A}^\mathrm{s}}_{6}$.
% An assumption about the topology of valve connections is given first.
% \begin{assumption}\label{assumtion:valve}
We now assume that the  links with open valves (PRVs and FCVs) \textit{do not} form a loop in water network.
% \end{assumption}
This assumption is practical due to the way valves are installed in water networks, and is corroborated by examining tens of water network templates. Hence, ${\m E_{\mathrm{L}}^\mathrm{J}}$ (${\m E_{\mathrm{J}}^\mathrm{L}}$) has linearly independent columns.  Hence, $ \m{\tilde{A}^\mathrm{s}}_{6}$  is also full column (row) rank.
 
 Note that $ \m{\tilde{A}^\mathrm{s}}_{1}$ is linearly column independent of the rest five column matrices, because each  submatrix in $ \m{\tilde{A}^\mathrm{s}}_{1}$ can not be eliminated with the corresponding $\m I_{n_r \times n_r}$, $\m I_{n_t \times n_t}$, $\m K_q^{\mathrm{P}}$, $\m K_q^{\mathrm{M}}$, or $\m E_{l'}$. Similarly, this also holds true for $ \m{\tilde{A}^\mathrm{s}}_{6}$. Hence, $ \m{\tilde{A}^\mathrm{s}}$ is full column rank. In fact, matrix $ \m{\tilde{A}^\mathrm{s}}$ is full row rank since  we have shown $\m E_{\m q}^\mathrm{J}$ and ${\m E_{\mathrm{J}}^\mathrm{L}}$ are full row rank, and  the row blocks are linearly independent of each other after splitting $ \m{\tilde{A}^\mathrm{s}}$ into six row-blocks. Thus, the original matrix $\m A^\mathrm{s}(k)$ is full  rank which results in a full rank $\m A[k] = \{\m A^\mathrm{s}[k], \m A^\mathrm{TK}[k]\}$.
 
Then, we consider $\m C$ under over-determined scenario, that is, $\m C$ selects extra measurements besides heads at tanks or reservoirs. In this case, the indices new added of measurements can be reformed in diagonal matrix, which cannot be eliminated with any other columns. Matrix $\m{\tilde{A}^\mathrm{s}}$ remains full rank under over-determined scenario. Hence, $\m A$  is full rank.
 
The second part of the proof derives~\eqref{equ:pse-solution} from~\eqref{equ:pse-matrix}. After multiplying $\m A^\top$ and $\m A$ on the left and right side of~\eqref{equ:pse-matrix}, we have
$(\m A^\top  \m A) { \operatorname {K}_{\m x \m x }}( \m A^\top  \m A )= \m A^\top { \operatorname {K}_{\m b \m b }} \m A$ (index $k$ is removed to make the equation clear), and the solution ${ \operatorname {K}_{\m x \m x }}$ can be obtained by inverting the full rank matrix $\m A^\top  \m A$. 
	\end{proof}
\begin{proof}(of Corollary~\ref{cor:inver})
	The proof of Corollary~\ref{cor:inver} is similar to the proof of Theorem~\ref{thm:inver} and hence omitted for brevity. 
\end{proof}

\section{Proof of Lemma~\ref{lemma:treeNonlinear} and Theorem~\ref{thm:looped}  }\label{sec:proof}
Before presenting the following proofs, we note that \textit{(i)} the results are applicable when the distribution is uniform. As a result, we stick to normally distributed uncertainty. \textit{(ii)} In the proof of Lemma~\ref{lemma:treeNonlinear} or Theorem~\ref{thm:looped}, we consider the tree network at first, then consider the looped/grid network. The valve models are linear in this paper, hence, are not considered in this proof.
\begin{proof} (of Lemma~\ref{lemma:treeNonlinear} for networks with nonlinear hydraulic models)
	
\noindent \textbf{Tree network:} According to  Definition 2.2.2 in~\cite{gertsbakh2016models}, a tree $\mathcal{T}$ is a graph in which any two vertices are connected by exactly one path, or equivalently a connected acyclic undirected graph.  Mathematical induction is used to prove that the flows follow the normal distribution in  a tree network  for normally distributed demand.
	
\noindent	\textbf{Base case}: We show that the statement holds for a tree $\mathcal{T}_{ij}$ with only two nodes with demand $d_i$ and $d_j$ and one link with flow $q_{ij}$. We assume that flow direction is from $i$ to $j$ which means $q_{ij} = d_{j}$.  When demand $d_j \sim \mathcal{N}(\mu_{d_j},\,\sigma_{d_j}^2)$, flow $q_{ij} \sim \mathcal{N}(\mu_{d_j},\,\sigma_{d_j}^2)$  based on Assumption~\ref{assumtion:all2}. This base case can be viewed as a subtree of another high level tree.
	
\noindent	\textbf{Inductive step}: We show that if base case (subtree) holds, then a tree $\mathcal{T}$ consisting of two base cases (subtrees) also holds; see Fig.~\ref{fig:tree} (left) for the notation. Note that two subtrees $\mathcal{T}_{jm}, \mathcal{T}_{jn}$ intersect at node $j$ connecting anther node $i$. According to the conservation of mass, $q_{ij}$ is also normally distributed since $q_{jm}$, $q_{jn}$ and $d_j$  follow a normal distribution:
	\begin{subequations} ~\label{equ:tree1}  
		\begin{align}
		&q_{ij} - q_{jm} - q_{jn} =  d_j  \\
		&q_{jm} = d_{m} \\
		&q_{jn} = d_{n} .
		\end{align}
	\end{subequations}
	Since both the base and inductive step have been performed, by mathematical induction the statement for flows holds for $\forall \ \mathcal{T}$ in  WDN.\\
	Second, we prove the statement for heads  in a tree network with nonlinear modeling of hydraulics.  The nonlinearities from a pump and a pipe are different, and we consider the pump first. \\ 
	\textbf{Case 1}: Link $ij$ is installed with a pump and the corresponding PDF can be noted as $f_{Q_{ij}}(q_{ij})$. We derive the PDF for nonlinear head increase of a pump $\Delta h_{ij}^{\mathrm{M}}$ next.
	We need to rewrite the cumulative distribution function (CDF) of $\Delta h_{ij}^{\mathrm{M}}$ in terms of the CDF of  $q_{ij}$.
	\begin{align*}
	&\F_{\Delta H_{ij}^{\mathrm{M}}}(\Delta h_{ij}^{\mathrm{M}}) = \prob(\Delta H_{ij}^{\mathrm{M}} \hspace{-2pt}<\hspace{-2pt} \Delta h_{ij}^{\mathrm{M}})=\prob(-(h_0 - r  q_{ij}^\beta )\hspace{-2pt} <\hspace{-2pt} \Delta h_{ij}^{\mathrm{M}}) \\
	&=\prob\left(q_{ij} \hspace{-2pt} <\hspace{-2pt} \left( \frac{h_0 + \Delta h_{ij}^{\mathrm{M}}}{r}\right)^{\frac{1}{\beta}}\right)  = \F_{Q_{ij}}\left(\hspace{-1pt}\left( \frac{h_0 + \Delta h_{ij}^{\mathrm{M}}}{r}\right)^{\frac{1}{\beta}}\right).
	\end{align*}
	The PDF equals to the derivative of CDF, hence, we obtain
	\begin{align}
	&f_{\Delta H_{ij}^{\mathrm{M}}}(\Delta h_{ij}^{\mathrm{M}}) = \F_{\Delta H_{ij}^{\mathrm{M}}}'\left(\Delta h_{ij}^{\mathrm{M}}\right) \label{equ:pdf_pump_new} \\ 
	&=\frac{1}{r \beta} \left(\left( \frac{h_0 + \Delta h_{ij}^{\mathrm{M}}}{r}\right)^{\frac{1}{\beta}-1}\right) f_{Q_{ij}}\left(\left( \frac{h_0 + \Delta h_{ij}^{\mathrm{M}}}{r}\right)^{\frac{1}{\beta}}\right). \notag
	\end{align}
	
	Usually, the suction side of a pump is connected with another node which is either a reservoir or a junction, see Fig.~\ref{fig:tree} (left).  If it is a reservoir, then the head is fixed as the elevation which a constant. If it is a junction, the head can be measured and the measurement noise is normally distributioned. In either way, $h_i$ is normally distributed. Now we consider a simple case, which is the node $i$ is a reservoir. At the delivery side of a pump, the head $h_j = h_i+ \Delta h_{ij}^{\mathrm{M}}$ is not normally distributed, since $f_{H_j}(h_j)$ is not a normal function. Next, we consider the pipe.\\
	\textbf{Case 2}: Link $ij$ is  a pipe. Similarly, the PDF of nonlinear head loss $\Delta h_{ij}^{\mathrm{P}}$ derived as follows.
	As for the nonlinear head loss model, we have
	\begin{equation*}
	\begin{aligned}
	\F_{\Delta H_{ij}^{\mathrm{P}}}(\Delta h_{ij}^{\mathrm{P}}) &= \prob(\Delta H_{ij}^{\mathrm{P}} < \Delta h_{ij}^{\mathrm{P}})=\prob\left(q_{ij} < \left( \frac{ \Delta h_{ij}^{\mathrm{P}}}{R_{ij}}\right)^{\frac{1}{\alpha}}\right)  \\
	&= \F_{Q_{ij}}\left( ( \frac{ \Delta h_{ij}^{\mathrm{P}}}{R_{ij}})^{\frac{1}{\alpha}}\right),
	\end{aligned}
	\end{equation*}
	and the corresponding PDF is 
	\begin{equation}\label{equ:pdf_pipe_new}   
	\begin{aligned}
	&f_{\Delta H_{ij}^{\mathrm{P}}}(\Delta h_{ij}^{\mathrm{P}}) = \F_{\Delta H_{ij}^{\mathrm{P}}}'(\Delta h_{ij}^{\mathrm{P}})  =\F_{Q_{ij}}'\left( \left( \frac{ \Delta h_{ij}^{\mathrm{P}}}{R_{ij}}\right)^{\frac{1}{\alpha}}\right)  \\ 
	&=\frac{1}{R_{ij}} \left( \left( \frac{ \Delta h_{ij}^{\mathrm{P}}}{R_{ij}}\right)^{\frac{1}{\alpha}-1}\right) f_{Q_{ij}}\left( \left( \frac{ \Delta h_{ij}^{\mathrm{P}}}{R_{ij}}\right)^{\frac{1}{\alpha}}\right).
	\end{aligned}
	\end{equation}
	Unfortunately, the PDF for $\Delta h_{ij}^{\mathrm{M}}$~\eqref{equ:pdf_pump_new} and $\Delta h_{ij}^{\mathrm{P}}$~\eqref{equ:pdf_pipe_new} are not normal.  
	
	From the above, we note that the head at arbitrary node is not normally distributed  in nonlinear tree network as mentioned in statement.\\
	
\noindent \textbf{Looped/grid network:}
A looped network is different from the tree one where flow variables can be only presented by the linear combination of demands. However, the random variables, flow  and demand, could be dependent on each other in a looped network making  it difficult to find the exact distribution even after linearization. 

Now we simply connect Junction $m$ and $n$ in Fig.~\ref{fig:tree} (left) together, and a looped network is constructed and shown in Fig.~\ref{fig:tree} (right). If we list all conservation of mass and energy for this looped network, we can obtain~\eqref{equ:loop1} and~\eqref{equ:loop2}.
Compared with~\eqref{equ:tree1}, the flow $q_{jm}$ and $q_{jm}$ are related now, and cannot be only presented by demands in~\eqref{equ:loop1}. Besides that, we have one more head loss equation in~\eqref{equ:loop2} due to the newly-added pipe. %Honestly, the PDF of $h_m$ and $h_n$ can still be obtained via~\eqref{equ:pipe1-loop2} and~\eqref{equ:pipe2-loop2} when $h_j$ is known. But at the same time, they also have to be constrained by~\eqref{equ:pipe3-loop2}.
% \vspace{-3em}
\begin{multicols}{2}
	\noindent \begin{subequations} ~\label{equ:loop1}  
		\begin{align}
		&q_{ij} - q_{jm} - q_{jn} =  d_j  \\
		&q_{jm} - q_{mn}= d_{m} \\
		&q_{jn} + q_{mn} = d_{n} 
		\end{align}
	\end{subequations}
	\begin{subequations} ~\label{equ:loop2}  
		\begin{align}
		h_i - h_j &= \Delta h_{ij}^{\mathrm{M}} \label{equ:pump-loop2} \\
		h_j - h_m &= \Delta h_{jm}^{\mathrm{P}} \label{equ:pipe1-loop2} \\
		h_j - h_n &= \Delta h_{jn}^{\mathrm{P}}  \label{equ:pipe2-loop2} \\
		h_m - h_n &= \Delta h_{mn}^{\mathrm{P}}  \label{equ:pipe3-loop2} 
		\end{align}
	\end{subequations}
\end{multicols}
From mentioned tree network case, we know that the PDFs for pump curve $\Delta h^{\mathrm{M}}$ and head loss function $\Delta h^{\mathrm{P}}$ in~\eqref{equ:loop2}  are not normally distributed.   After solving~\eqref{equ:loop1} and~\eqref{equ:loop2} together, it is impossible that  heads and flows follow the normal distribution.
\end{proof}
\begin{proof} (of Theorem~\ref{thm:looped} for networks with linearized hydraulic models )\\
Note that the proof of Theorem~\ref{thm:looped}  is based on the proof of Lemma~\ref{lemma:treeNonlinear}, please refer to it before reading this proof.\\
\noindent \textbf{Tree network:}
	Considering the same tree $\mathcal{T}_{ij}$ with two nodes with demand $d_i$ and $d_j$ and one link with flow $q_{ij}$.  The flow $q_{ij}  \sim \mathcal{N}(\mu_{q_{ij}}, \sigma_{q_{ij}}^2)$ based on Lemma~\ref{lemma:treeNonlinear}. Thus, only the heads need to be disscused for linearized tree network.\\
Similarly,  we still need to discuss the pump case and pipe case, but the derivation of PDF of linearized pump and pipe  modeling are omitted, since the linear transformation of normal distribution would remain as normal distribution. The pump is considered next.\\
\textbf{Case 1}: Link $ij$ is a pump and the corresponding linearized modeling is~\eqref{equ:head-flow-pump-linear}. According to the rule of expectation and variance of linear combination for random variables, $\Delta h_{ij}^{\mathrm{M}} \sim \mathcal{N}(k_{ij}^{\mathrm{M}} \mu_{q_{ij}}+b_{ij}^{\mathrm{M}},\,(k_{ij}^{\mathrm{M}}\sigma_{q_{ij}})^2)$. Similarly, $h_j = h_i + \Delta h_{ij}^{\mathrm{M}}$. Hence, we have  $ h_{j} \sim \mathcal{N}(k_{ij}^{\mathrm{M}} \mu_{q_{ij}}+b_{ij}^{\mathrm{M}} + h_i,\,(k_{ij}^{\mathrm{M}}\sigma_{q_{ij}})^2)$. For the pipe case, we have\\
\textbf{Case 2}: Link $ij$ is a pipe and the corresponding linearized modeling is~\eqref{equ:head-flow-pipe-linear}. , $\Delta h_{ij}^{\mathrm{P}} \sim \mathcal{N}(k_{q_{ij}}^{\mathrm{P}} \mu_{q_{ij}}+k_{c_{ij}}^{\mathrm{P}} \mu_{c_{ij}}+b_{ij}^{\mathrm{P}},\,(k_{q_{ij}}^{\mathrm{P}}\sigma_{q_{ij}})^2 + (k_{c_{ij}}^{\mathrm{P}}\sigma_{q_{ij}})^2)$. According to, $h_j = h_i + \Delta h_{ij}^{\mathrm{P}}$, we have  $ h_{j} \sim \mathcal{N}(k_{q_{ij}}^{\mathrm{P}} \mu_{q_{ij}}+k_{c_{ij}}^{\mathrm{P}} \mu_{c_{ij}}+b_{ij}^{\mathrm{P}} + h_i,\,(k_{q_{ij}}^{\mathrm{P}}\sigma_{q_{ij}})^2 + (k_{c_{ij}}^{\mathrm{P}}\sigma_{q_{ij}})^2)$, because we know $h_i$ is normally distributed based on the discussion of nonlinear case, so $h_j$ also follows the normal distribution. 

\noindent \textbf{Looped/grid network:}
In WDN, after linearization for~\eqref{equ:loop2} around operating points, we obtain the first order Taylor series model~\eqref{equ:loop3}.
\begin{subequations} ~\label{equ:loop3}  
	\begin{align}
	h_i - h_j + k_{ij}^{\mathrm{M}} q_{ij} &=  b_{ij}^{\mathrm{M}}  \label{equ:pump-loop3} \\
	h_j - h_m - k_{jm}^{\mathrm{P}} q_{jm}  &= b_{jm}^{\mathrm{P}} + k_{jm}^{\mathrm{C}} c_{jm} \label{equ:pipe1-loop3} \\
	h_j - h_n - k_{jn}^{\mathrm{P}} q_{jn}  &=  b_{jn}^{\mathrm{P}} + k_{jn}^{\mathrm{C}} c_{jn}  \label{equ:pipe2-loop3} \\
	h_m - h_n - k_{mn}^{\mathrm{P}} q_{mn} &=  b_{mn}^{\mathrm{P}}  + k_{mn}^{\mathrm{C}} c_{mn}   \label{equ:pipe3-loop3} 
	\end{align}
\end{subequations}
From~\eqref{equ:loop1} and~\eqref{equ:loop3}, we can see that \textit{(i)} flow $\m q$ is related with $\m d$, and the head $\m h$ is determined by linearized model with $\m q$ and $\m c$. Hence, head $\m h$ is also related with $\m d$. \textit{(ii)} 
Each random variable in $\m d$ are assumed to be independent of each other, but we note that they are correlated with flow $\m q$ and head $\m h$. \textit{(iii)} in order to find each head and flow, we need to solve the linear equations consisting of~\eqref{equ:loop1} and~\eqref{equ:loop3}, and this is the different part comparing with a tree network where the expectation and variance of heads and flows can be solved directly.  Now we can summarize the PSE for any type of networks in WDN. Next, we show the proof of Theorem~\ref{thm:looped} for general looped network.

For a general looped network, we need to solve the linear equation~\eqref{equ:combin}, and if we consider it in sufficient scenario, then $\m x = \m A^{-1}  \m b$. The proof on invertibility of $\m A$ is in Appendix~\ref{app:inver}. According to the rule of variance of linear combination, the head is also distributed normally when  $\m b$ follow a normal distribution. The above process can be repeated using~\eqref{equ:weighted-hy-matrix} under over-determined scenario. Thus, the statement holds true for  any scenario in WDN.
	
	Based on the above two cases, the statement holds true for any network.
\end{proof}
\vspace{-5em}
\begin{IEEEbiography}[{\includegraphics[width=1in,height=1.25in,clip,keepaspectratio]{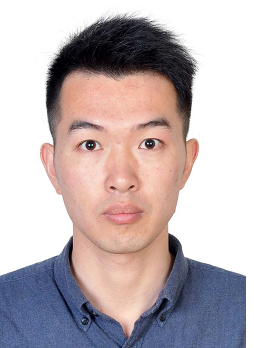}}]{Shen Wang} received the B.E. degree in automation from Anhui Jianzhu University, Hefei, China, in 2013 and the M.E. degree in control science and engineering from the University of Science and Technology of China, Hefei, China, in 2016. He is currently pursuing a Ph.D. degree in Electrical Engineering at the University of Texas at San Antonio, Texas, USA.

His current research interests include optimal control in cyber-physical systems with special focus on power systems and drinking water distribution networks.
\end{IEEEbiography}

\vspace{-5em}

\begin{IEEEbiography}
	[{\includegraphics[width=1in,height=1.25in,clip,keepaspectratio]{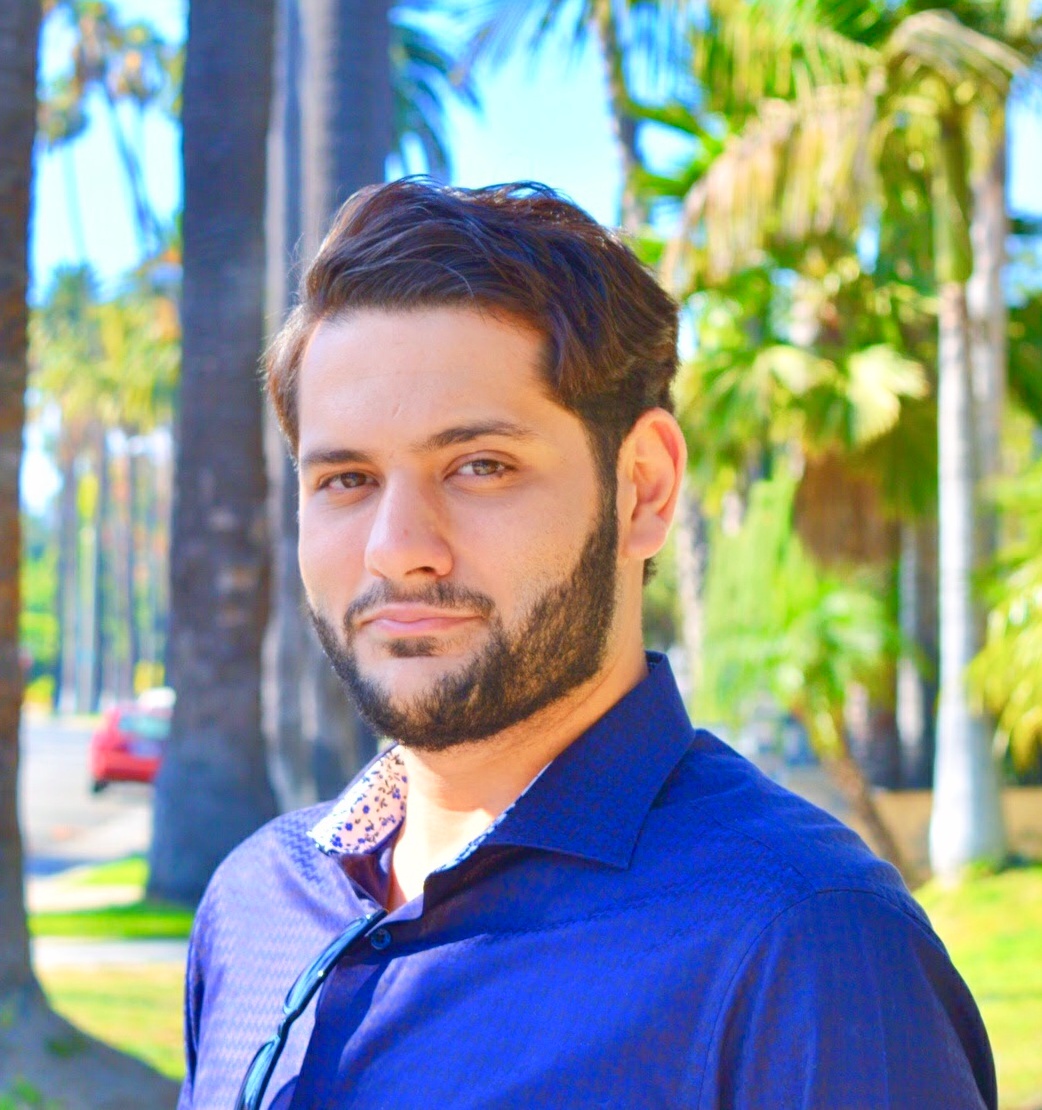}}]
	{Ahmad F. Taha} is  an assistant professor with the Department of Electrical and Computer Engineering at the University of Texas, San Antonio. He received the B.E. and Ph.D. degrees in Electrical and Computer Engineering from the American University of Beirut, Lebanon in 2011 and Purdue University, West Lafayette, Indiana in 2015.  Dr. Taha is interested in understanding how complex cyber-physical systems (CPS) operate, behave, and \textit{misbehave}. His research focus includes optimization, control, and security of CPSs with applications to power, water, and transportation networks.  Dr. Taha is an editor of IEEE Transactions on
	Smart Grid and the editor of the IEEE Control Systems Society Electronic
	Letter (E-Letter).
\end{IEEEbiography}

\vspace{-5em}

\begin{IEEEbiography}
	[{\includegraphics[width=1in,height=1.25in,clip,keepaspectratio]{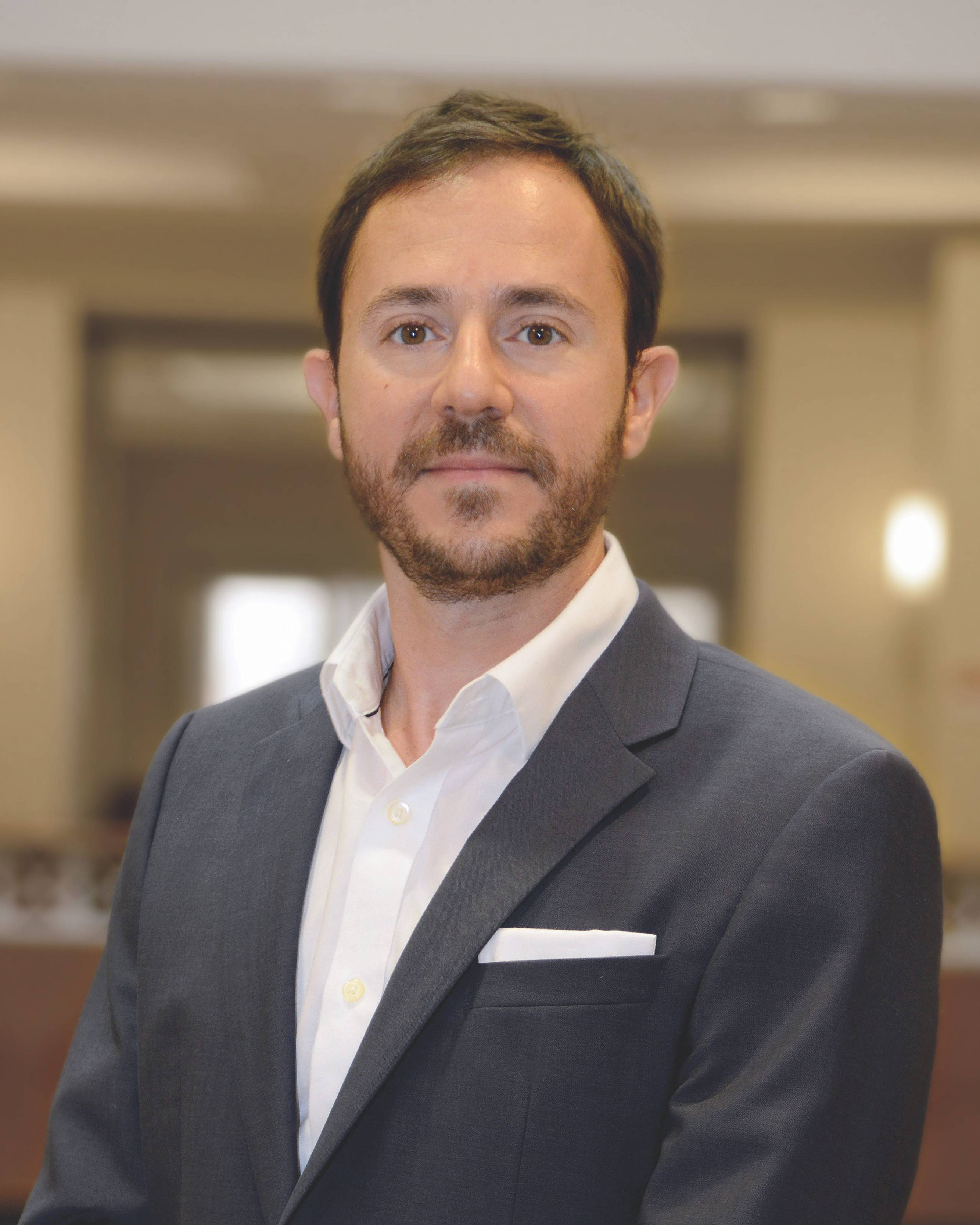}}]
{Nikolaos Gatsis} (Member, IEEE) received the Diploma degree (Hons.) in electrical and computer engineering from the University of Patras, Patras, Greece, in 2005, and the M.Sc. degree in electrical engineering, in 2010, and the Ph.D. degree in electrical engineering with minor in mathematics, in 2012, from the University of Minnesota, Minneapolis, MN, USA. He is currently an Associate Professor with the Department of Electrical and Computer Engineering, University of Texas at San Antonio, San Antonio, TX, USA and a Lutcher Brown Professorship Endowed Fellow.  His research focuses on optimal and secure operation of smart power grids and other critical infrastructures, including water distribution networks and the Global Positioning System. Dr. Gatsis is a recipient of the NSF CAREER award. He was a Co-Guest Editor for the IEEE Journal on Selected Topics in Signal Processing (Special Issue on Signal and Information Processing for Critical Infrastructures). He also presented to the 2020 NSF Engineering CAREER Proposal Writing Workshop. 
\end{IEEEbiography}
\vspace{-11em}

\begin{IEEEbiography}
	[{\includegraphics[width=1in,height=1.25in,clip,keepaspectratio]{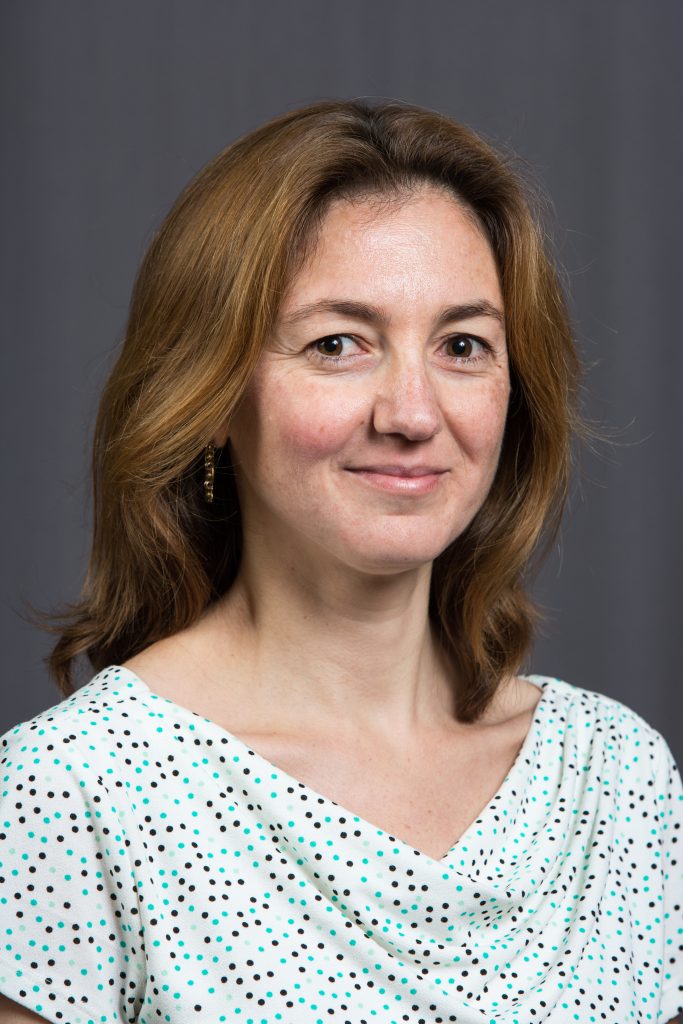}}]
	{Lina Sela} is an Assistant Professor in Environmental Water Resources Engineering, in the Civil, Architectural and Environment Engineering department at the University of Texas at Austin. Dr. Sela’s research focuses on improving the efficiency of water distribution systems facing challenges related to finite water sources, aging infrastructure, and population growth. Her work relies on integrating the increasingly available digital information from distributed sensing devices with physical‐based models to improve operations and management of urban water systems by integrating the cyber components in decision making process. Her current research projects include a network of pressure and acoustic sensors distributed in the City of Austin continuously monitoring the state of the water distribution system. Her collaboration with industry and water utilities provides a mechanism for technology development and implementation, realization of research works, and enables end technology users’ involvement in research. 
\end{IEEEbiography}

\vspace{-11em}

\begin{IEEEbiography}
	[{\includegraphics[width=1in,height=1.25in,clip,keepaspectratio]{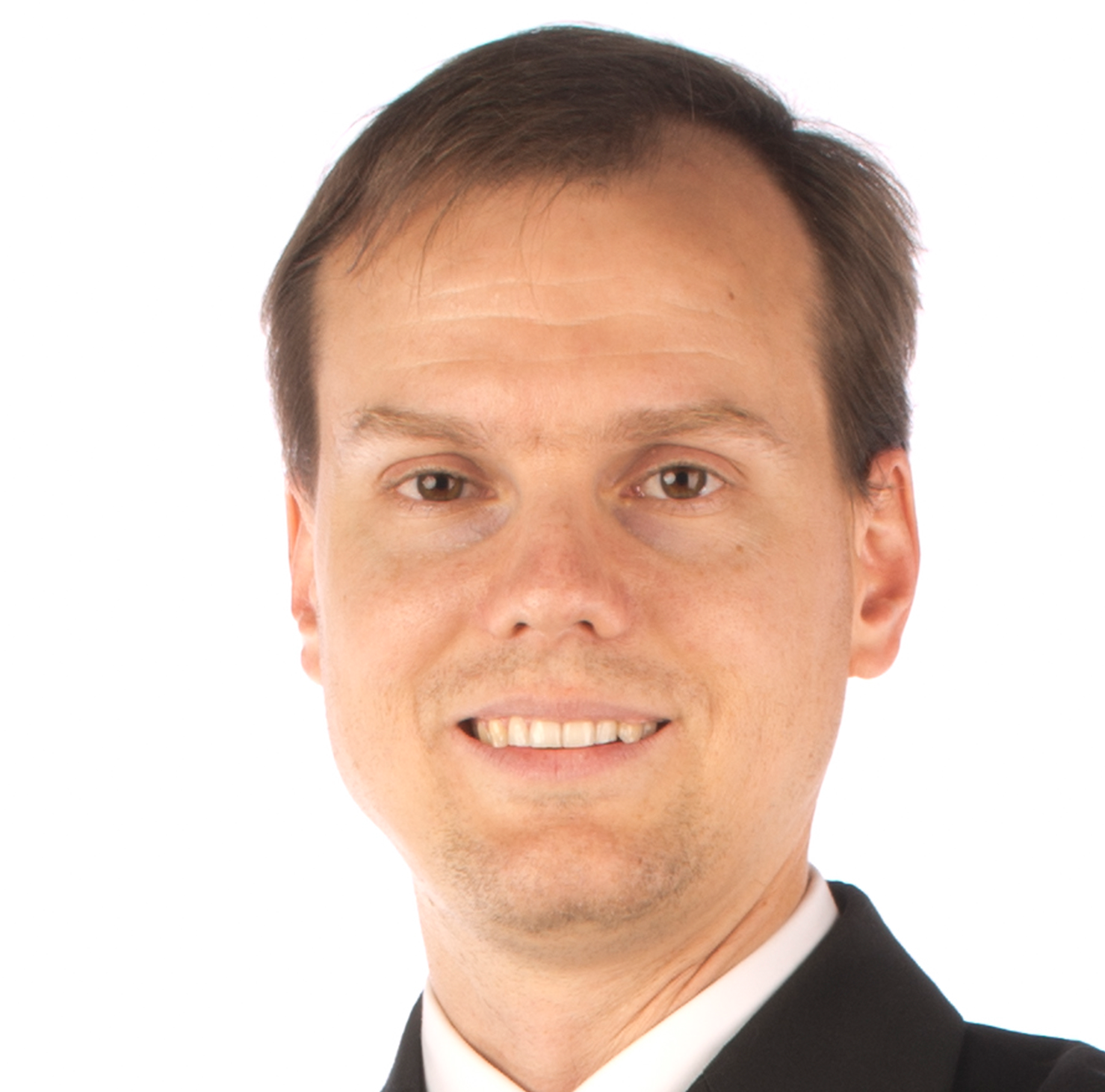}}]
	{Marcio H. Giacomoni} received a B.S. degree in civil engineering from University of Brasilia, Brazil, in 2002, a M.S. degree in water resources engineering from Federal University of Rio Grande do Sul, Porto Alege, Brazil, in 2005, and a PhD. Degree in civil engineering from Texas A\&M University, College Station, Texas, USA, in 2012. He is currently an Associate Professor in the Department of Civil and Environmental Engineering at the University of Texas at San Antonio, Texas, USA. His research interests include water resources systems analysis, sustainability of the built and natural environments, stormwater management and green infrastructure, and resilience and security of cyber-physical systems.
\end{IEEEbiography}
\end{document}